\documentclass[11pt,a4paper]{article}

\usepackage[T1]{fontenc}
\usepackage[utf8]{inputenc}
 
\usepackage{a4wide}

\usepackage{xcolor}
\definecolor{lgray}{rgb}{0.95,0.95,0.95}
\definecolor{yel}{rgb}{1,0.993,0.95}
\definecolor{mydarkblue}{rgb}{0,0.0,0.8}
\newcommand\bl[1]{{\color{mydarkblue}#1}}
\newcommand\h{\bl{h}}
\definecolor{mydarkred}{rgb}{0.65,0.0,0.0}
\newcommand\dr[1]{{\color{mydarkred}#1}}
\newcommand\g{\dr{g}}
\definecolor{mydarkgreen}{rgb}{0,0.42,0}
\newcommand\dg[1]{{\color{mydarkgreen}#1}}
\newcommand\ff{\dg{f}}
\definecolor{mybrown}{rgb}{0.75,0.45,0} 
\newcommand\br[1]{#1}

\usepackage{amsmath, amsfonts}
\usepackage{cite}
\DeclareMathOperator*{\argmin}{arg\,min}
\DeclareMathOperator*{\minimize}{minimize}

\usepackage[amsmath,thmmarks,hyperref,framed]{ntheorem}

\definecolor{mydarkgreen2}{rgb}{0,0.5,0}
\usepackage[dvips]{hyperref}
\hypersetup{
  colorlinks=true,
  allcolors = mydarkgreen2,
  urlcolor = mydarkred,
  frenchlinks=false,
  pdfborder={0 0 0},
  naturalnames=false,
  hypertexnames=false,
  breaklinks
}

\usepackage[capitalize,nameinlink]{cleveref}

\crefname{section}{section}{sections}
\crefname{subsection}{subsection}{subsections}
\Crefname{section}{Section}{Sections}
\Crefname{subsection}{Subsection}{Subsections}

\numberwithin{equation}{section}

\Crefname{figure}{Figure}{Figures}
\crefformat{equation}{\textup{#2(#1)#3}}
\crefrangeformat{equation}{\textup{#3(#1)#4--#5(#2)#6}}
\crefmultiformat{equation}{\textup{#2(#1)#3}}{ and \textup{#2(#1)#3}}
{, \textup{#2(#1)#3}}{, and \textup{#2(#1)#3}}
\crefrangemultiformat{equation}{\textup{#3(#1)#4--#5(#2)#6}}%
{ and \textup{#3(#1)#4--#5(#2)#6}}{, \textup{#3(#1)#4--#5(#2)#6}}{, and \textup{#3(#1)#4--#5(#2)#6}}

\Crefformat{equation}{#2Equation~\textup{(#1)}#3}
\Crefrangeformat{equation}{Equations~\textup{#3(#1)#4--#5(#2)#6}}
\Crefmultiformat{equation}{Equations~\textup{#2(#1)#3}}{ and \textup{#2(#1)#3}}
{, \textup{#2(#1)#3}}{, and \textup{#2(#1)#3}}
\Crefrangemultiformat{equation}{Equations~\textup{#3(#1)#4--#5(#2)#6}}%
{ and \textup{#3(#1)#4--#5(#2)#6}}{, \textup{#3(#1)#4--#5(#2)#6}}{, and \textup{#3(#1)#4--#5(#2)#6}}

\crefdefaultlabelformat{#2\textup{#1}#3}

\newtheoremstyle{plainnn}%
{\item[\hskip\labelsep \theorem@headerfont ##1\ ##2\theorem@separator]}%
{\item[\hskip\labelsep \theorem@headerfont ##1\ ##2\ (##3)\theorem@separator]}%

\usepackage{framed,pstricks}
\shadecolor{lgray}

\makeatletter
\renewtheoremstyle{plain}
 {\item[\hskip\labelsep \theorem@headerfont ##1\ \textup{##2\theorem@separator}]}%
 {\item[\hskip\labelsep \theorem@headerfont ##1\ \textup{##2}]\textup{(##3)\theorem@separator}\ }
\renewtheoremstyle{nonumberplain}%
 {\item[\theorem@headerfont\hskip\labelsep ##1\theorem@separator]}%
{\item[\theorem@headerfont\hskip\labelsep ##3\theorem@separator]}
\makeatother

\newcommand{\proofbox}{\vbox{\hrule height0.6pt\hbox{\vrule height1.3ex width0.6pt\hskip0.8ex\vrule width0.6pt}\hrule height0.6pt}}

\theorempreskip{0em}
\theorempostskip{0mm}
\theoremframepreskip{0cm}
\theoremframepostskip{0cm}
\theoreminframepreskip{0cm}
\theoreminframepostskip{0cm}

\theoremstyle{nonumberplain}
\theoremheaderfont{\normalfont\itshape}
\theorembodyfont{\normalfont}
\theoremseparator{.}
\theoremsymbol{\proofbox}
\newtheorem{proof}{Proof}

\theoremstyle{plain}
\theoremheaderfont{\normalfont\sc}
\theorembodyfont{\normalfont\itshape}
\theoremseparator{.}
\theoremsymbol{}
\newshadedtheorem{theorem}{Theorem}[section]

\newcommand{\newsiamthm}[2]{
  \theoremstyle{plain}
  \theoremheaderfont{\normalfont\sc}
  \theorembodyfont{\normalfont\itshape}
  \theoremseparator{.}
  \theoremsymbol{}
  \newshadedtheorem{#1}[theorem]{#2}
}

\newsiamthm{lemma}{Lemma}
\newsiamthm{corollary}{Corollary}
\newsiamthm{proposition}{Proposition}
\newsiamthm{definition}{Definition}

\usepackage[ntheorem]{empheq}

\usepackage{authblk}
 
\title{\textbf{Proximal Splitting Algorithms\\ for Convex Optimization:\\A Tour of Recent Advances, with New Twists}} 
\date{\today}
\author[1]{Laurent Condat\thanks{Corresponding author. Contact: see his webpage \url{https://lcondat.github.io}}\thanks{The first author did  part of this work during a stay at Ritsumeikan University in 2018, hosted by the second and fourth authors, thanks to a fellowship from the Japanese Society for the Promotion of Science (JSPS), no. L17565.}}
\author[2]{Daichi Kitahara}
\author[3,4]{Andr\'es Contreras\thanks{The third author contributed to this work during his visit to the first author at GIPSA-lab, Grenoble, France, in Spring 2019. He was supported by the CMM-ANID PIA Grant AFB170001.}}
\author[5]{Akira Hirabayashi}

\affil[1]{Visual Computing Center, King Abdullah University of Science and Technology (KAUST), Thuwal, Saudi Arabia}
\affil[2]{Division of Electrical, Electronic and Infocommunications Engineering, Osaka University, Suita, Japan}
\affil[3]{Centro de Modelamiento Matemático, Santiago, Chile}
\affil[4]{Department of Mathematical and Physical Sciences, Faculty of Engineering, UC-Temuco, Temuco, La Araucanía, Chile}
\affil[5]{College of Information Science and Engineering, Ritsumeikan University, Kusatsu, Japan}

\date{February 2023\\
Authors' final version. To appear in \emph{SIAM Review}}

\newcommand*\mycolouredbox[1]{%
\setlength{\fboxsep}{3pt}\colorbox{lgray}{\ #1\ }}
{\endEmphEqMainEnv}

\newcommand{\ceq}[2]{
\begin{empheq}[box=\colorbox{yel}
]{#1}#2\end{empheq}}

\begin{document}

\maketitle

\begin{center}\large\textbf{Abstract}\end{center}
Convex nonsmooth optimization problems, whose solutions live in very high dimensional spaces, have become ubiquitous. To solve them, the class of first-order algorithms known as proximal splitting algorithms is particularly adequate: they consist of simple operations, handling the terms in the objective function separately. In this overview, we demystify a selection of recent proximal splitting algorithms: we present them within a unified framework, which consists in applying splitting methods for monotone inclusions 
in  primal-dual product spaces, with well-chosen metrics. Along the way, we easily derive new variants of the algorithms and revisit existing convergence results, extending the parameter ranges in several cases.  In particular, we emphasize that  when the smooth term in the objective function is quadratic, e.g., for least-squares problems, convergence is guaranteed with larger values of the relaxation parameter than previously known. Such larger values are usually beneficial for the convergence speed in practice.\bigskip

\medskip


{
  \hypersetup{linkcolor=black}
  \tableofcontents
}
\newpage

\section{Introduction}

Optimization, also called mathematical programming, is the search for a best object in a set with respect to some criterion. Mathematically, this amounts to characterizing a point where a function attains its minimum value, and computationally, this often consists in exhibiting an iterative algorithm in which a variable 
converges to such an optimal point. When a function is convex, 
its minimizers are global: a point which attains a minimum value 
in its neighborhood is actually optimal in the whole search space. Thus, it is important to know about the principles behind
convex optimization algorithms, before considering more general problems. 
This paper is a constructive and self-contained introduction to the class of proximal splitting algorithms, which are efficient  for large-scale convex optimization.

Many problems in statistics, machine learning, 
signal and image processing, control, and many other fields can be formulated as convex optimization problems~\cite{pal09, sta10, bac12, pol15, cha16, sta16} \br{of the form:
\begin{equation}
    \minimize_{x\in\mathcal{X}}\, \sum_{m=1}^M g_m(L_m x),
\end{equation}
for some $M\geq 1$, where $\mathcal{X}$ is a real Hilbert space, the $g_m$ are convex, possibly nonsmooth,  functions, and the $L_m$ are linear operators. We will see that this generic formulation encompasses the possible presence of constraints.} In the age of ``big data,'' with the explosion in size and complexity of the data to process, it is increasingly important to be able to solve such optimization problems, whose solutions live in very high dimensional spaces~\cite{ela10,sra11,bub15,emr16,pol19}. 
For instance, in
 image processing and computer vision, 
 there are one or several variables to estimate for each pixel, with typically several millions of pixels. 
There is an extensive literature on \emph{proximal splitting} algorithms for solving convex optimization problems, with applications in various fields~\cite{com10,sra11,bot14,par14,kom15,ryu16,glo16,cse16,bec17,com21}. They consist of simple, easy-to-compute steps that can 
deal with the terms in the 
objective function separately: 
if a function is smooth, its gradient can be used, whereas for a nonsmooth function, the proximity operator will be called instead (hence the wording \emph{proximal algorithm}). Since the proximity operator of a sum of functions does not have a computable form, in general, it is necessary to handle individually the proximity operators of the functions that appear in the problem; this is what \emph{splitting} refers to, a notion reminiscent to the divide-and-conquer principle, which is ubiquitous in computational sciences.

It has been known for several decades that the \emph{dual} problem,\footnote{More rigorously, \emph{some} dual problem, since there is no unique way to define the dual of an optimization problem, in general.} associated to the (primal) optimization problem under consideration, can be easier to solve and, if used, may simplify the computation of the primal solution. Yet, the development of \emph{primal-dual} algorithms, which solve the primal and dual problems simultaneously, in an intertwined way, is more recent. 
In this tutorial survey, we introduce a selection 
of primal-dual proximal splitting algorithms, which have been developed over the last decade. 
We present them in a unified way, by solving a monotone inclusion expressed in a well-chosen primal-dual product space. Along the way, some new variants are naturally obtained. An important idea for the development of primal-dual algorithms is \emph{preconditioning}, or equivalently changing the metric of the ambient Hilbert space~\cite{poc11,com142,bre15}. We will make extensive use of this notion throughout the paper.

In our presentation, we put the emphasis on the potential to \emph{overrelax} the algorithms: let us consider an algorithm producing a new estimate $z^{(i+\frac{1}{2})}$ of some variable at iteration number $i+1$, given the previous estimate $z^{(i)}$. Instead of setting $z^{(i+1)}=z^{(i+\frac{1}{2})}$ before continuing with the next iteration, overrelaxation consists in the extrapolated update $z^{(i+1)}=z^{(i)}+\rho^{(i)} (z^{(i+\frac{1}{2})}-z^{(i)})$, for some $\rho^{(i)}>1$. This step is easy to implement and can  improve the convergence speed of the algorithm. Thus, it is important to study the conditions on the parameters $\rho^{(i)}$ that guarantee convergence. In particular, we show that if the optimization problem involves a smooth term that is quadratic, convergence is guaranteed with larger values of the relaxation parameter than previously known. Regularized least-squares problems are an important class of problems to which this extension applies.

Also, a desirable property of the proximal splitting algorithms is that they can be parallelized. Thus, when minimizing the sum of a large number of functions on a parallel computing architecture, each processor can handle one of the functions. Mathematical splitting gracefully blends into the physical splitting of the computing hardware in that case. The need for such parallel optimization algorithms is rapidly growing in distributed or federated learning \cite{kon16,kai21,mal20} and artificial intelligence. The last part of the paper is dedicated to parallel versions of the algorithms.\medskip

The paper is organized as follows. \br{In \cref{secadd}, we present the most important notions and principles underlying the class of iterative algorithms for large-scale convex nonsmooth optimization under study. }
In \cref{sec2}, we present relaxed versions of the forward-backward splitting algorithm to solve monotone inclusions. This analysis is subsequently applied in \cref{sec3} to the Loris--Verhoeven algorithm, which is a primal-dual forward-backward algorithm. In \cref{sec4}, we analyze the Chambolle--Pock algorithm and its particular case, the Douglas--Rachford algorithm, which is  equivalent to the Alternating Direction Method of Multipliers (ADMM). For these too, we derive convergence results for relaxed versions, which extend previously known results. In \cref{sec5}, we study an algorithm that we call the Generalized Chambolle--Pock algorithm, since it is applicable to the minimization of two functions composed with two linear operators and reverts back to the Chambolle--Pock algorithm when one of the linear operators is the identity. The algorithm has been presented in the literature as a linearized or preconditioned version of the ADMM, but without relaxation. In \cref{sec6}, we study an algorithm proposed independently by the first author and by B.\ C.\ V\~u that is also a primal-dual forward-backward algorithm. \Cref{sec7} is devoted to recently proposed algorithms based on a three-operator splitting scheme, which can be viewed as the fusion of the forward-backward and Douglas--Rachford two-operator splitting schemes. Along the way, new algorithm variants are proposed, for example, in \cref{sec61,secnaq}.
Finally, in \cref{sec8}, we propose parallel versions of all these algorithms by applying them in product spaces.

\section{\br{A Brief Introduction to Convex Analysis and Fixed-Point Algorithms}}\label{secadd}

\subsection{Notions of convex analysis}

In this section, we introduce some notions and notations that will be used throughout the paper. The reader can find a much more complete account  of convex analysis and operator theory in textbooks, e.g.,~\cite{boy04,bau17}, or in papers like~\cite{ryu16,com18}. 

We consider optimization over real Hilbert spaces: a real Hilbert space is a vector space equipped with a real inner product $\langle\cdot\,,\cdot\rangle$. When we make use of a norm without further specification, this is the norm $\|\cdot\|=\langle\cdot\,,\cdot\rangle^{1/2}$ induced by the inner product. This general formalism allows us to optimize with respect to
vectors, matrices, or tensors~\cite{gan11},  which are real-valued or complex-valued with some Hermitian properties, but also more complicated objects like quaternions~\cite{miz19}. The Hilbert spaces can be of infinite dimension; this can be useful in control~\cite{you97} or when solving PDEs, and we will see that most results hold in this general setting. However, there are a few technicalities to keep in mind:

1) Let $L:\mathcal{X}\rightarrow \mathcal{U}$ be a linear operator, with $\mathcal{X}$ and $\mathcal{U}$ two real Hilbert spaces. We define the operator norm of $L$ 
as $\|L\|=\sup\{\|Lx\| : x\in\mathcal{X}$ and $\|x\|\leq 1\}$. 
If $\|L\|<+\infty$, $L$ is said to be \textbf{bounded}. $L$ is bounded if and only if it is continuous. In the paper, we will assume all linear operators to be bounded. If $\mathcal{X}$ is of finite dimension, $L$ is necessarily bounded.  

2) A sequence of points $(x^{(i)})_{i\in\mathbb{N}}$ in a real Hilbert space $\mathcal{X}$ is said to converge \textbf{weakly} to a point $x^\star\in\mathcal{X}$ if, for every $y\in\mathcal{X}$, $\langle x^{(i)},y\rangle \rightarrow \langle x^\star,y\rangle$. On the other hand, $(x^{(i)})_{i\in\mathbb{N}}$ is said to converge \textbf{strongly} to $x^\star\in\mathcal{X}$ if $\|x^{(i)}-x^\star\|\rightarrow 0$. Strong convergence implies weak convergence.  If $\mathcal{X}$ is of finite dimension, both notions are equivalent and  
we just say that  $(x^{(i)})_{i\in\mathbb{N}}$ converges to $x^\star$. 
We will see that most convergence results state weak convergence; this is not a weakness of the proof techniques, and strong convergence in infinite-dimensional spaces cannot be guaranteed without further assumptions~\cite{bui202}. In particular, given a continuous operator $T$ on $\mathcal{X}$, if a sequence $(x^{(i)})_{i\in\mathbb{N}}$ converges weakly to  $x^\star\in\mathcal{X}$, one cannot deduce that $\big(T(x^{(i)})\big)_{i\in\mathbb{N}}$ converges weakly to $T(x^\star)$. This makes the convergence proofs of algorithms more difficult to derive.\medskip

For the rest of \cref{secadd}, let $\mathcal{X}$ be a real Hilbert space. A subset $\Omega$ of $\mathcal{X}$ is said to be \textbf{convex} if, for every $(x,y)\in \Omega^2$ and $a\in (0,1)$, $ax+(1-a)y \in \Omega$. Note that a convex set can be open, like the interval $(0,1)$ in $\mathbb{R}$, or closed, like $[0,1]$; 
it can be bounded or unbounded, like $[0,+\infty)$. 
The intersection of convex sets is convex. 
The union of convex sets is not convex,  in general.

A function $f:\mathcal{X}\rightarrow\mathbb{R}\cup\{+\infty\}$ is said to be \textbf{convex} if, for every $(x,y)\in \mathcal{X}^2$ and $a\in (0,1)$, 
\begin{equation}
f\big(ax+(1-a)y\big) \leq a f(x) + (1-a)f(y).\label{eqconv}
\end{equation}
We note that $f$ is allowed to take the value $+\infty$, a distinctive and very useful feature of convex analysis, compared to classical analysis. For calculus, we only need to adopt the following rules: for every $t>0$, $t(+\infty)=+\infty$, and for every $t\in \mathbb{R}$, $t + (+\infty)=+\infty$. 
Note that for every $t\in\mathbb{R}$, the intervals $(t,+\infty]$ are considered open in $\mathbb{R}\cup\{+\infty\}$, so that, for instance,  the function $t\in\mathbb{R}\mapsto (1/t$ if $t>0$, $+\infty$ otherwise$)$ is continuous on $\mathbb{R}$~\cite[Example 1.22]{bau17}. The domain of $f$ is the convex 
set  $\mathrm{dom}\,f=\{x\in\mathcal{X}:f(x)\neq +\infty\}$. $f$ is said to be \textbf{proper} if its domain is nonempty.

Indicator functions relate the two notions of convexity of a set and of a function: given a set $\Omega\subset \mathcal{X}$, we define the \textbf{indicator function} of $\Omega$, denoted by $\imath_\Omega$, as:
\begin{equation}
\imath_\Omega: x\in\mathcal{X}\mapsto \left\{\begin{array}{ll}
0&\mbox{if }x\in\Omega,\\
+\infty&\mbox{otherwise}.
\end{array}\right.
\end{equation}
Then, if $\Omega$ is convex, $\imath_\Omega$ is convex; if $\Omega$ is nonempty,  
$\imath_\Omega$ is proper. Also, given two convex sets $\Omega_1\subset \mathcal{X}$ and $\Omega_2\subset \mathcal{X}$,
\begin{equation}
\imath_{\Omega_1 \cap \Omega_2} = \imath_{\Omega_1} + \imath_{\Omega_2}.
\end{equation}
Indicator functions are important, since they allow us to remove constraints and integrate them into the objective function to minimize. Indeed, given a nonempty convex set $\Omega\subset \mathcal{X}$ and a convex function $f$ on $\mathcal{X}$,
\begin{align}
\minimize_{x\in\Omega}\, f(x)&\ \equiv\  \minimize_{x\in\mathcal{X}}\, f(x)\ \ \mbox{s.t.}\  \ x\in\Omega,\\
&\ \equiv\  \minimize_{x\in\mathcal{X}}\, f(x) + \imath_\Omega(x).\label{eqi24}
\end{align}
We see here the interest in allowing the functions to take the value $+\infty$: it is used to exclude some parts of $\mathcal{X}$ from the set of possible solutions. 
Thus, we will consider the problem of minimizing a sum of convex functions, knowing 
that, thanks to indicator functions, this covers the possible presence of constraints.
\medskip

A convex optimization problem consists in estimating a \textbf{minimizer} $x^\star\in\mathcal{X}$, supposed to exist, 
of a proper convex function $f:\mathcal{X}\rightarrow\mathbb{R}\cup\{+\infty\}$, 
that is, a point $x^\star$ at which $f$ 
attains its minimum value $\min_{x\in\mathcal{X}} f(x)$. The set of minimizers of $f$, which is convex, is denoted by $\argmin_{x\in\mathcal{X}} f(x)$. Thus, $x^\star \in \argmin_{x\in\mathcal{X}} f(x)$ if and only if, for every $x\in\mathcal{X}$, $f(x^\star)\leq f(x)$. Note that the function $t\in\mathbb{R}\mapsto (1/t$ if $t>0$, $+\infty$ otherwise$)$  is proper, convex, and bounded from below, but does not have any minimizer.

An optimization algorithm to minimize a proper convex function $f:\mathcal{X}\rightarrow\mathbb{R}\cup\{+\infty\}$, which has a minimizer, constructs a sequence $(x^{(i)})_{i\in \mathbb{N}}$ of points such that $f(x^{(i)})$ decreases overall. Thus, for the process to be fruitful, if $(x^{(i)})_{i\in \mathbb{N}}$ converges to $x^\star$, we must have 
$\lim\inf_{i\in \mathbb{N}} f(x^{(i)})\geq f(x^\star)$. Otherwise, even if $f(x^{(i)})$ converges to $\min f$, we cannot deduce that $x^\star$ is a minimizer of $f$. So, we will only consider functions having this lower semicontinuity property. Equivalently, a function $f$ is said to be \textbf{lower semicontinuous} if the convex set $\{x\in\mathcal{X}: f(x)\leq t\}$ is closed for every $t\in\mathbb{R}$. A continuous function is lower semicontinuous. 
A proper convex lower semicontinuous function is continuous on the interior of its domain~\cite[Corollary 8.39]{bau17}, so that lower semicontinuity concerns the 
boundary of the domain, where the function jumps from real values to $+\infty$.

We denote by $\Gamma_0(\mathcal{X})$ the set of \textbf{convex, proper, lower semicontinuous} functions  from $\mathcal{X}$ to $\mathbb{R}\cup \{+\infty\}$. We will only consider such functions in optimization problems. If $\Omega\subset\mathcal{X}$ is nonempty and closed, $\imath_\Omega\in \Gamma_0(\mathcal{X})$. 
If $f\in\Gamma_0(\mathcal{X})$ and $L$ is a bounded linear operator with values in $\mathcal{X}$, $f\circ L $ is 
convex and lower semicontinuous~\cite[Proposition 9.5]{bau17}, where $\circ$ denotes the composition of functions; 
it is proper if and only if $\mathrm{dom}\,f \cap \mathrm{ran}\,L\neq \emptyset$, where  $\mathrm{ran}\, L$ denotes the range of $L$. Also, if $f\in\Gamma_0(\mathcal{X})$ and $g\in\Gamma_0(\mathcal{X})$, $f+g$ is convex and lower semicontinuous~\cite[Corollary 9.4]{bau17}; it is proper if and only if $\mathrm{dom}\,f \cap \mathrm{dom}\,g\neq \emptyset$. 

The difficulty of convex optimization stems from the fact that even if $f$ and $g$ are two simple convex functions with known minimizers, $f+g$ is difficult to minimize, in general. A well-known example is the LASSO problem~\cite{tib96}, which consists in minimizing over $x\in\mathcal{X}$ the function $\|Lx-b\|^2 + \|x\|_1$, where $L$ is a linear operator, $b$ an element, and $\|\cdot\|_1$ is the $l_1$ norm; note that every norm on $\mathcal{X}$ is convex. The minimizer of the least-squares term can be obtained efficiently by solving a linear system, the minimizer of the $l_1$ norm is the zero element, but there is no straightforward way of finding a minimizer of the sum of the two terms, in general.  Algorithms able to solve such problems are precisely the purpose of this paper.\medskip

When dealing with vectors, a linear operator can be represented by a matrix, so that  the application of a linear operator to a vector can be viewed as a matrix-vector product.  But since we are placing ourselves in the general setting of Hilbert spaces, some  definitions 
are in order. Let $L:\mathcal{X}\rightarrow\mathcal{U}$ be a bounded linear operator for some real Hilbert space $\mathcal{U}$.
The \textbf{adjoint} operator $L^*:\mathcal{U}\rightarrow\mathcal{X}$ of $L$ is the only bounded linear operator such that $\langle Lx,u\rangle_\mathcal{U}=\langle x,L^*u\rangle_\mathcal{X}$ for every $x\in\mathcal{X}$ and  $u\in\mathcal{U}$, where we have denoted by $\langle \cdot,\cdot\rangle_\mathcal{X}$ and $\langle \cdot,\cdot\rangle_\mathcal{U}$ the inner products in $\mathcal{X}$ and  $\mathcal{U}$, respectively, to differentiate them. $L$ is said to be \textbf{self-adjoint} if $L^*=L$. We have $\|L\|=\|L^*\|=\sqrt{\|LL^*\|}=\sqrt{\|L^*L\|}$. $\mathrm{Id}$ denotes the identity operator. Let $P$ be a self-adjoint bounded linear operator on $\mathcal{X}$.  $P$ is said to be \textbf{positive} if $\langle x,Px\rangle \geq 0$ for every $x\in\mathcal{X}$, and strongly positive if $P-a\mathrm{Id}$ is positive for some real $a>0$. If $P$ is positive, for every $y\in\mathcal{X}$, the function $x\mapsto \langle Px+y,x\rangle $ belongs to $\Gamma_0(\mathcal{X})$; 
for instance, with $P=\mathrm{Id}$ and $y=0$, we obtain that  $\|\cdot\|^2\in\Gamma_0(\mathcal{X})$. \medskip

In this paper, we consider first-order algorithms; that is, which exploit information about the gradient or subdifferential of the functions. 
Let us define these notions. A  continuous real-valued function $f\in\Gamma_0(\mathcal{X})$ is smooth, or differentiable, if at every $x\in\mathcal{X}$, there exists an element of $\mathcal{X}$, called the \textbf{gradient} of $f$ at $x$ and denoted by $\nabla f(x)$, such that for every  $e\in \mathcal{X}$,
\begin{equation}
f(x+e) -f(x) - \langle e, \nabla f(x)\rangle = o(\|e\|), \quad\mbox{as } \|e\|\rightarrow 0.
\end{equation}
That is, the affine function $y\mapsto f(x) + \langle y-x, \nabla f(x)\rangle$ is a first-order approximation of $f$ around $x$. Note that, as a consequence of convexity, this affine function is a minorant of $f$: for every  $y\in \mathcal{X}$, 
\begin{equation}
f(y) \geq f(x) + \langle y-x, \nabla f(x)\rangle.
\end{equation}
Thus, $\nabla f$, the gradient of $f$, is the operator on $\mathcal{X}$: $x\mapsto \nabla f(x)$. It is continuous on $\mathcal{X}$~\cite[Corollary 17.43]{bau17}. 
 For smooth functions $f$ and $g$, we have $\nabla(f+g)=\nabla f +\nabla g$, and for a linear operator $L$, $\nabla (f\circ L)=L^*(\nabla f)L$, where $AB$ denotes the composition of the two operators $A$ and $B$. 
We note that the gradient of the quadratic convex function $x\mapsto \langle Px+y,x\rangle $, for some positive self-adjoint bounded linear operator $P$ on $\mathcal{X}$ and some $y\in\mathcal{X}$, is $x\mapsto 2Px+y$.

We want to deal with nonsmooth functions; 
for instance, the absolute value function $f:\mathbb{R}\rightarrow\mathbb{R}:t\mapsto |t|$, and by extension the $l_1$ norm, are convex but not differentiable at $0$. Moreover, we allow convex functions to take the value $+\infty$. Thus, we need a more general notion than the gradient, which captures first-order information. This is where the subdifferential comes to the rescue. Let $f\in \Gamma_0(\mathcal{X})$. For every $x\in\mathcal{X}$, we define the \textbf{subdifferential} of $f$ at $x$, denoted by $\partial f(x)$, as the closed convex  set~\cite[Proposition 16.4]{bau17}
\begin{equation}
\partial f(x) =  \big\{u\in\mathcal{X}\ :\ \forall y\in\mathcal{X},\ f(x)+\langle y-x, u\rangle\leq f(y)\big\}.
\end{equation}
The elements of $\partial f(x)$ are called subgradients of $f$ at $x$. A subgradient $u\in \partial f(x)$ can be viewed as the gradient of the affine minorant $y\mapsto  f(x)+\langle y-x, u\rangle$ of $f$. 
If $f$ is differentiable at $x$, there is only one subgradient, which is the gradient: $\partial f(x)=\{\nabla f(x)\}$. On the other hand, if $f(x)=+\infty$, $\partial f(x)=\emptyset$. A classical example is the function absolute value, for which $\partial |\cdot| (0)=[-1,1]$. Another example is $f:t\in\mathbb{R}\mapsto (0$ if $t\leq 0$, $+\infty$ otherwise$)$, for which $\partial f (0) = [0,+\infty)$. It is important to note that for any $g\in \Gamma_0(\mathcal{X})$ and $x\in\mathcal{X}$~\cite[Proposition 16.6]{bau17},
\begin{equation}
\partial f(x) +\partial g(x) \subset \partial (f+g)(x),\label{eqi29}
\end{equation}
where the sum of the two sets is the Minkowski sum, but the inclusion may be strict~\cite[Remark 16.8]{bau17}. 
Similarly, given a bounded linear operator $L$, for every $x\in\mathcal{X}$, $L^*\partial f(Lx) \subset \partial (f \circ L)(x)$. We denote by $2^\Omega$ the set of all subsets of a set $\Omega$. Then the subdifferential of $f$ is the \textbf{set-valued operator}  $\partial f:\mathcal{X}\rightarrow 2^\mathcal{X}, x\mapsto \partial f(x)$. We assimilate a set-valued operator $M:\mathcal{X}\rightarrow 2^\mathcal{X}$ with exactly one element in each set $M(x)$ as a single-valued operator $M:\mathcal{X}\rightarrow\mathcal{X}$, and conversely. For instance, for a smooth convex function $f$, we have $\partial f = \nabla f$.\medskip

We will consider 
convex optimization problems of the form:
\begin{equation}
    \minimize_{x\in\mathcal{X}}\, \sum_{m=1}^M g_m(L_m x),\label{eq00}
\end{equation}
where $M\geq 1$, each $L_m$ is a  bounded linear operator from $\mathcal{X}$ to some real Hilbert space $\mathcal{U}_m$, possibly the identity 
if $\mathcal{U}_m=\mathcal{X}$, and each $g_m\in\Gamma_0(\mathcal{U}_m)$. 
The well-known \textbf{Fermat's rule} 
states that $x^\star\in\mathcal{X}$ is a solution to 
\eqref{eq00} if and only if
\begin{equation}
0\in \partial \Big( \sum_{m=1}^M g_m \circ L_m \Big)(x^\star).\label{eqfer2}
\end{equation}
This fundamental property is easy to demonstrate: let $f\in\Gamma_0(\mathcal{X})$ and $x^\star\in\mathcal{X}$. Then 
\begin{align}
x^\star\in \argmin_{x\in\mathcal{X}} f(x)&\Leftrightarrow \forall x \in \mathcal{X},\ f(x^\star)\leq f(x) \\
&\Leftrightarrow\forall x \in \mathcal{X},\ f(x^\star)+\langle x-x^\star,0\rangle \leq f(x)\\
&\Leftrightarrow 0 \in \partial f(x^\star).
\end{align}
As already mentioned in \eqref{eqi29}, for every $x\in\mathcal{X}$ we have
\begin{equation}
 \Big(\sum_{m=1}^M L_m^* (\partial g_m) L_m\Big)(x)  \subset \partial \Big( \sum_{m=1}^M g_m \circ L_m \Big)(x),\label{eq01}
\end{equation}
but the inclusion may be strict. All algorithms considered in this paper actually find a solution $x^\star\in \mathcal{X}$ satisfying
\begin{equation}
0\in\sum_{m=1}^M L_m^*  \partial g_m(L_m x^\star),\label{eq02}
\end{equation}
if such a solution exists. But in view of \eqref{eq01}, the existence of a minimizer of $\sum_m g_m\circ L_m$ is not sufficient to guarantee the existence of a solution to \eqref{eq02}. 
Therefore, throughout the paper, we assume that for every optimization problem considered, \textbf{the solution set of the corresponding monotone inclusion \eqref{eq02} is nonempty}. Given \eqref{eq01} and Fermat's rule, a solution $x^\star$ of \eqref{eq02} is a solution to the optimization problem \eqref{eq00}. If the functions satisfy a so-called \textbf{qualification constraint}, we have an equality instead of an inclusion in \eqref{eq01}; in that case, it is sufficient to show that a solution to \eqref{eq00} exists. 
A qualification constraint is~\cite[Proposition 4.3]{com12}:
\begin{equation}
\bigcap_{m=1}^M L_m^{-1} (\mathrm{sri}\;\mathrm{dom}\;g_m)\neq\emptyset,\label{eqqualif0}
\end{equation}
where, for any linear operator $L:\mathcal{X}\rightarrow \mathcal{U}$ and $\Omega\subset\mathcal{U}$, $L^{-1}(\Omega)=\{x\in\mathcal{X}:Lx\in\Omega\}$, and $\mathrm{sri}\;\Omega$ denotes the strong relative interior of a convex set $\Omega$, which 
is the set of points $u\in\Omega$ such that the cone generated by $-u+\Omega$ is a closed linear subspace.
 Other qualification constraints can be found in \cite[Proposition 4.3]{com12}. Note that for $\sum_{m=1}^M g_m\circ L_m$ to be proper, it is necessary to have $\cap_{m=1}^M L_m^{-1} (\mathrm{dom}\;g_m)\neq\emptyset$. Since the strong relative interior is even larger than the interior of a set, the condition \eqref{eqqualif0} is really not a lot to ask.\medskip
 
 We have seen that the subdifferential is an important notion, in particular with Fermat's rule. But the algorithms considered in this paper to solve optimization problems involving nonsmooth functions will not make use of subgradients directly. Instead, they will make calls to the proximity operators of the functions: we define the \textbf{proximity operator}~\cite{mor62} of any function $f\in \Gamma_0(\mathcal{X})$
as
\begin{equation}
    \mathrm{prox}_f:\mathcal{X}\rightarrow\mathcal{X},x\mapsto \argmin_{x'\in\mathcal{X}}\Big(f(x')+{\textstyle\frac{1}{2}}\|x'-x\|^2\Big).\label{eqdefprox}
\end{equation}
By strong convexity of $f+\frac{1}{2}\|\cdot-x\|^2$, its minimizer always exists and is unique, so that the proximity operator is well defined. For any nonempty closed convex set $\Omega\subset\mathcal{X}$ and $\tau>0$,
\begin{equation}
    \mathrm{prox}_{\tau \imath_\Omega} = \mathrm{proj}_\Omega: x\mapsto \argmin_{x'\in\Omega}\|x'-x\|
\end{equation}
is the projection onto $\Omega$, which does not depend on $\tau$. Thus, proximal algorithms can be viewed as generalizing algorithms to find points in intersections of convex sets~\cite{com93}. 
Although the definition of the proximity operator is implicit, it has a closed form for many functions of practical interest. For instance, for the absolute value (and by extension for the $l_1$ norm by elementwise application), this is soft-thresholding: we have, for any $\tau>0$, 
\begin{equation}
    \mathrm{prox}_{\tau |\cdot|}:t\in\mathbb{R}\mapsto  \mathrm{sign}(t)\max(|t|-\tau,0).
\end{equation}
There are fast and exact methods to compute the proximity operators for a large class of functions~\cite{cha092,com10,con132,con16,pus17,elg18}; see also the website \url{http://proximity-operator.net}. In this paper, we consider \textbf{proximal splitting algorithms}, which are iterative algorithms producing an estimate of a solution to a convex optimization problem of the form \eqref{eq00}, by calling at each iteration the proximity operator or the gradient of each function $g_m$, as well as the operators $L_m$ and $L_m^*$.

An optimization problem can often be written in the form \eqref{eq00}
in several equivalent but different ways, with each formulation having  
appropriate algorithms to solve it. For instance, when a function $\imath_\Omega$ for some convex set $\Omega$ appears in  the problem, it is implicitly assumed, in general, that one knows how to project onto $\Omega$. For instance, to minimize $f\in \Gamma_0(\mathcal{X})$ subject to the constraint $L\,\cdot=b$ for some linear operator $L$ and element $b$, 
the two following formulations are equivalent:
\begin{align}
     \minimize_{x\in\mathcal{X}}&\, f(x)+\imath_{\{x'\in\mathcal{X}\ :\ Lx'=b\}}(x),\\
      \minimize_{x\in\mathcal{X}}&\, f(x)+\imath_{\{b\}}(Lx).
\end{align}
The first formulation is appropriate if one knows how to project on the affine constraint set efficiently, whereas the second formulation is preferable if one does not and asks for an algorithm making calls to $L$ and $L^*$. 
\medskip

For every $f\in \Gamma_0(\mathcal{X})$, we denote by $f^*$ the conjugate  of $f$ defined by 
\begin{equation}
f^*:x\mapsto \sup_{x'\in\mathcal{X}} \big(\langle x,x'\rangle -f(x')\big),\label{eqcon}
\end{equation}
which belongs to $\Gamma_0(\mathcal{X})$ as well. We have $f^{**}=f$. An important property is 
\textbf{the Moreau identity}, which 
 allows us to compute the proximity operator of $f^*$ from that of $f$, and vice versa: for every $\tau>0$,
\begin{equation}
\mathrm{prox}_{\tau f}(x)=x-\tau\,\mathrm{prox}_{f^*/\tau}(x/\tau).\label{eqmoreau}
\end{equation}
We do not define here further notions related to duality, which will be introduced progressively during the presentation of primal-dual algorithms.

\subsection{Notions of operator theory and fixed-point algorithms}

Let  $M:\mathcal{X}\rightarrow 2^\mathcal{X}$ be a set-valued operator; as mentioned above, single-valued operators can be viewed as set-valued operators, so that the subsequent notions apply to them as well.
 The graph of $M$ is the set $\{(x,v)\in\mathcal{X}^2 : v\in Mx\}$ (we may use the shortened notation $Mx$ for $M(x)$, usually adopted for linear operators, with any operator). The \textbf{inverse} operator of $M$, denoted by $M^{-1}$, is the set-valued operator with graph $\{(v,x)\in\mathcal{X}^2 : v\in Mx\}$. That is, $v\in Mx$ if and only if $x\in M^{-1}v$. 
 For instance, for any bounded linear operator $L$, $L^{-1}$ is well defined everywhere, in this sense. But $L^{-1}x$ can be the empty set if $x\notin \mathrm{ran}\,L$, or can contain several elements: $L^{-1} 0$ is the kernel of $L$.  The inverse relates the subdifferentials of a function $f\in\Gamma_0(\mathcal{X})$ and its conjugate $f^*\in\Gamma_0(\mathcal{X})$ defined in \eqref{eqcon} via the simple relation~\cite[Corollary 16.30]{bau17}:
 \begin{equation}
 \partial f^*=(\partial f)^{-1}.
 \end{equation}

$M$  is said to be \textbf{monotone} if, for every
$(x,x')\in\mathcal{X}^2$, $v \in Mx$, $v'\in Mx'$, 
 \begin{equation}
\langle x-x',v-v'\rangle\geq 0.
 \end{equation}
 We note that for a self-adjoint bounded linear operator, being positive is equivalent to being monotone. 
 A monotone operator $M$ is said to be \textbf{maximally monotone} if there exists no monotone operator whose graph strictly contains the graph of $M$. For every $f\in\Gamma_0(\mathcal{X})$, 
 \begin{equation}
  \partial f \mbox{ is maximally monotone}
  \end{equation}\cite[Theorem 20.25]{bau17}. The first derivative of a smooth convex function $f:\mathbb{R}\rightarrow \mathbb{R}$ is nondecreasing, so monotonicity can be viewed as a generalization of the nondecreasing property to multidimensional and set-valued operators.  
 The problem of finding $x^\star$ such that $0\in Mx^\star$ is called a monotone inclusion. Thus, as we have seen with Fermat's rule \eqref{eq00}--\eqref{eqfer2}, convex optimization problems consist in solving monotone inclusions. This is why monotone operator theory plays an essential role in optimization.
 
If $M$ is monotone, we define its \textbf{resolvent}, denoted by $J_M$, as the operator
 \begin{equation}
J_M=(M+\mathrm{Id})^{-1}.
 \end{equation}
$J_M$ is a well-defined single-valued operator on $\mathcal{X}$. Now, let us apply Fermat's rule to the definition of the proximity operator in \eqref{eqdefprox}: if $y=\mathrm{prox}_f(x)$ for any $f\in\Gamma_0(\mathcal{X})$ and $x\in\mathcal{X}$, then $0 \in \partial f(y) + y-x$. Equivalently, $y = (\partial f + \mathrm{Id})^{-1}(x)$. That is,
\begin{equation}
\mathrm{prox}_f=J_{\partial f}.
\end{equation}
Thus, proximal algorithms handle nonsmooth functions by calls to their proximity operators; this can be viewed as a regularized and implicit way of dealing with their set-valued subdifferentials. But this comes at a price: $\mathrm{prox}_{f+g}$ cannot be expressed from $\mathrm{prox}_{f}$ and $\mathrm{prox}_{g}$, in general; in rare cases, $\mathrm{prox}_{f+g}=\mathrm{prox}_{f}\circ \mathrm{prox}_{g}$ \cite{pus17}. Also $\mathrm{prox}_{f\circ L}$ does not have a computable form, in general. For instance, the proximity operator of the least-squares term $\frac{\tau}{2}\|L\cdot-b\|^2$ for some bounded linear operator $L:\mathcal{X}\rightarrow \mathcal{U}$, $b\in\mathcal{U}$ and $\tau>0$, is
 \begin{equation}
x\mapsto    (\tau L^*L +\mathrm{Id})^{-1} (x+\tau L^* b),
\end{equation}
which has a closed form, but this amounts in practice to solving a linear system, which is generally too costly. 
Thus, even the apparently simple problem of minimizing the sum of two nonsmooth functions is in fact difficult. Thus, \textbf{splitting} is necessary: we have to call the proximity operator of each simple-enough function individually when minimizing a sum of functions, and  
deal with the linear operators separately.\medskip

Let $T:\mathcal{X}\rightarrow\mathcal{X}$ be a single-valued operator. A \textbf{fixed point} of $T$ is any $x^\star\in\mathcal{X}$ satisfying 
\begin{equation}
x^\star = Tx^\star.\label{eqfp1}
\end{equation}
$T$ is said to be $\beta$-\textbf{Lipschitz continuous} for some real $\beta\geq 0$ if, for every $(x,x')\in\mathcal{X}^2$,
\begin{equation}
\|Tx-Tx'\|\leq \beta \|x-x'\|. \label{eqlip}
\end{equation}
We note that a bounded linear operator $L$ is $\|L\|$-Lipschitz continuous. 
$T$ is said to be \textbf{nonexpansive} if it is 1-Lipschitz continuous. $T$ is said to be $\alpha$-\textbf{averaged} for some $\alpha\in (0,1)$ if  there exists a nonexpansive operator $T':\mathcal{X}\rightarrow\mathcal{X}$ such that 
\begin{equation}
T = (1-\alpha)\mathrm{Id}+\alpha T'.
\end{equation}
If $T$ is $\alpha$-averaged for some $\alpha\in (0,1)$, it is $\alpha'$-averaged for every $\alpha' \in (\alpha,1)$ and nonexpansive; in addition, 
for every $\rho\in (0,\frac{1}{\alpha})$, the \textbf{relaxed} operator
\begin{equation}
\mathrm{Id} + \rho(T-\mathrm{Id}) \label{eqav1}
\end{equation}
is $\rho\alpha$-averaged. 
\textbf{Firmly nonexpansive} means $(1/2)$-averaged. Finally, $T$ is said to be  $\xi$-\textbf{cocoercive} for some $\xi>0$ if $\xi T$ is firmly nonexpansive or, equivalently~\cite[Proposition 4.4]{bau17}, if, for every $(x,x')\in\mathcal{X}^2$,
\begin{equation}
\xi \| Tx - Tx'\|^2\leq \langle x-x', Tx - T x'\rangle.
\end{equation}

Let $\mathcal{Z}$ be a real Hilbert space and $T$ be a nonexpansive operator on $\mathcal{Z}$ with nonempty set of fixed points.
A \textbf{fixed-point algorithm}, initialized with some $z^{(0)} \in \mathcal{Z}$, consists in iterating 
\begin{equation}
z^{(i+1)}=T z^{(i)},\label{eqfp2}
\end{equation}
where $i$ is the iteration counter and $z^{(i)}$ denotes the version of the variable $z$ produced during the $i$th iteration. 
Note that we have changed the name of the variable from $x$ to $z$ to emphasize that  $z^{(i)}$  may not be the variable converging to a solution of the optimization problem at hand: it can be a concatenation of variables, like $z^{(i)}=(x^{(i)},u^{(i)})$, with $x^{(i)}$ converging to a desired solution and $u^{(i)}$ an auxiliary variable,  the value of which is not of primary concern. Alternatively, 
 it might be that there is an internal variable $x^{(i)}$ among the set of operations modeled by the abstract operator $T$, which converges to a solution, whereas $z^{(i)}$ is not interesting in itself.

Then, by combining \eqref{eqfp1} and \eqref{eqlip}, we have, for every fixed point $z^\star$ of $T$ and $i\geq 0$,
\begin{equation}
\|z^{(i+1)} - z^\star\| \leq \|z^{(i)}-z^\star\|. \label{eqfp3}
\end{equation}
Thus, the variable $z^{(i)}$ grows closer and closer at every iteration to the set of fixed points of $T$.  But this is not sufficient to obtain convergence of the sequence $(z^{(i)})_{i\in\mathbb{N}}$ to some fixed point. 
A counterexample is  the rotation operator of angle $\theta\in [0,2\pi)$ in $\mathbb{R}^2$, represented by the matrix
\begin{equation}
T= \left[\begin{array}{cc} \cos(\theta)&-\sin(\theta)\\
\sin(\theta)&\cos(\theta)\end{array}\right].
\end{equation}
In that case the inequality in \eqref{eqfp3} is an equality, so that  $(z^{(i)})_{i\in\mathbb{N}}$ generated by \eqref{eqfp2} turns over and over again in the plane and does not converge to the fixed point $(0,0)$ of $T$. 

Thus, we have to be slightly more strict in our assumptions: we will suppose that $T$ in \eqref{eqfp2} is $\alpha$-averaged for some $\alpha\in(0,1)$, and not only nonexpansive. And in view of \eqref{eqav1}, given some sequence $(\rho^{(i)})_{i\in\mathbb{N}}$ of relaxation parameters, 
we will consider the \textbf{relaxed iteration}
\begin{align}
    &\left\lfloor
   \begin{array}{l}
    z^{(i+\frac{1}{2})}= T z^{(i)}\\
     z^{(i+1)}=z^{(i)}+\rho^{(i)} (z^{(i+\frac{1}{2})}-z^{(i)}).
      \end{array}\right. 
    \label{eq0}
\end{align}
Note that with $\rho^{(i)}=1$, there is no relaxation and we recover the iteration \eqref{eqfp2}. While there is no interest in doing underrelaxation with $\rho^{(i)}$ less than 1, \textbf{overrelaxation} with $\rho^{(i)}$ larger than 1 may be beneficial to the convergence speed, as is observed in practice in many cases. The idea of using overrelaxation comes from the following intuition: if the sequence $(z^{(i)})_{i\in\mathbb{N}}$   generated by \eqref{eqfp2} converges to some fixed point $z^\star$, this means that $z^{(i+1)}=T z^{(i)} $ tends to be closer to $z^\star$ than $z^{(i)} $, on average. Hence, we may want, starting at $z^{(i)}$, to move further in the direction $T z^{(i)}-z^{(i)}$, which improves the estimate. This is exactly what the relaxed iteration \eqref{eq0} does if $\rho^{(i)}>1$.

Now, we can state the fundamental result on which we will rely to establish convergence of optimization algorithms:\medskip

\noindent\colorbox{lgray}{\parbox{0.985\textwidth}{\textsc{Krasnosel'ski\u{\i}--Mann theorem} (algorithm \eqref{eq0})~\cite[Proposition 5.16]{bau17}. \emph{
 Let $\mathcal{Z}$ be a real Hilbert space and let $T:\mathcal{Z}\rightarrow \mathcal{Z}$ be an $\alpha$-averaged operator for some $\alpha\in(0,1)$, with a nonempty set of fixed points. Let $(\rho^{(i)})_{i\in\mathbb{N}}$ be a sequence in $[0,1/\alpha]$, such that $\sum_{i\in \mathbb{N}} \rho^{(i)} (1-\alpha\rho^{(i)})=+\infty$. Let $z^{(0)}\in\mathcal{Z}$. 
Then the sequence $(z^{(i)})_{i\in \mathbb{N}}$ defined by the iteration \eqref{eq0} 
converges weakly to some fixed point $z^\star$ of $T$, and $(z^{(i+\frac{1}{2})}-z^{(i)})_{i\in \mathbb{N}}$ converges strongly to 0.
}}}
\medskip

We note that as a direct implication of this theorem, $(z^{(i+\frac{1}{2})})_{i\in \mathbb{N}}$ converges weakly to the same fixed point $z^\star$. In practice, one might prefer to consider $z^{(i+\frac{1}{2})}$ instead of $z^{(i+1)}$ as an estimate of the solution, because it is in the range of $T$.

In this theorem, the assumption that $T$ has a nonempty set of fixed points is important: consider, for instance, $T=\mathrm{Id}+c$ for some $c\neq 0$, which is $\alpha$-averaged for every $\alpha\in(0,1)$. Then $T$ has no fixed point and the iteration \eqref{eq0} diverges.

Finally, we note that $\rho^{(i)}=1$ is allowed in the theorem. In practice, we recommend setting $\rho^{(i)}=1$ when running an algorithm for the first time, before trying larger values.

All algorithms presented in this paper are instances of 
 the simple form \eqref{eq0}, and their convergence analysis relies on the Krasnosel'ski\u{\i}--Mann theorem.\medskip

It remains to see how to construct averaged operators whose fixed points are solutions to optimization problems. To that aim, we will combine gradients and proximity operators of functions like Lego bricks. Indeed, the set of averaged operators is stable by convex combinations (a convex combination is a linear combination with nonnegative weights whose sum is one) and compositions:\medskip

\cite[Proposition 4.42]{bau17}  (convex combination of averaged operators). \emph{Let $(T_m)_{m=1}^M$  be a  family of operators for some $M\geq 1$, such that every $T_m:\mathcal{X}\rightarrow\mathcal{X}$ is $\alpha_m$-averaged for some $\alpha_m\in(0,1)$. Let $(\omega_m)_{m=1}^M$ be real numbers in $(0,1]$, such that $\sum_{m=1}^M \omega_m=1$.  Set $\alpha=\sum_{m=1}^M \omega_m \alpha_m\in (0,1)$ and $T=\sum_{m=1}^M \omega_m T_m$. Then $T$ is $\alpha$-averaged. }\medskip

\cite{ogu02}\cite[Proposition 4.44]{bau17}  (composition of averaged operators). \emph{Let $T_1:\mathcal{X}\rightarrow\mathcal{X}$ be an $\alpha_1$-averaged operator and $T_2:\mathcal{X}\rightarrow\mathcal{X}$ be an $\alpha_2$-averaged operator for some $\alpha_1\in(0,1)$ and $\alpha_2\in(0,1)$. Set 
\begin{equation}
\alpha = \frac{\alpha_1+\alpha_2-2\alpha_1\alpha_2}{1-\alpha_1\alpha_2}\label{eqcompav}
\end{equation}
and $T=T_2T_1$. Then $\alpha\in (0,1)$ and $T$ is $\alpha$-averaged. 
}\medskip

Moreover, if $M:\mathcal{X}\rightarrow2^\mathcal{X}$ is a maximally monotone operator,
\begin{equation}
J_M \mbox{ is firmly nonexpansive}
\end{equation}\cite[Proposition 23.8]{bau17}. And since the proximity operator is the resolvent of the subdifferential, it is firmly nonexpansive as well. Finally, when we deal with a smooth function using calls to its gradient, we will assume that the gradient is Lipschitz continuous. Indeed, the Baillon--Haddad theorem~\cite[Corollary 18.17]{bau17} states the following: let $f:\mathcal{X}\rightarrow\mathbb{R}$ be a differentiable convex function whose gradient $\nabla f$ is $\beta$-Lipschitz continuous for some $\beta>0$. Then 
\begin{equation}
\nabla f \mbox{ is }\textstyle \frac{1}{\beta}\mbox{-cocoercive}.
\end{equation}

We now have all the ingredients to start designing algorithms. Let us start with the minimization of a differentiable convex function $f:\mathcal{X}\rightarrow\mathbb{R}$, whose gradient $\nabla f$ is $\beta$-Lipschitz continuous for some $\beta>0$, and whose set of minimizers is nonempty. \textbf{Gradient descent}, initialized at some $x^{(0)}\in\mathcal{X}$, consists in iterating, for $i\geq 0$,
\begin{equation}
x^{(i+1)}= x^{(i)} - \gamma \nabla f (x^{(i)})
\end{equation}
for some parameter $\gamma \in (0,\frac{2}{\beta})$. It is called a descent method because one can show that $f(x^{(i+1)})\leq f(x^{(i)})$, so that the algorithm ``goes down'' along the function values. To analyze this algorithm, we first observe that $x^{(i+1)}= T x^{(i)}$ for the operator $T=\mathrm{Id} -\gamma \nabla f$. Since $\frac{1}{\beta}\nabla f$ is firmly nonexpansive, it is easy to show that $T$ is $\alpha$-averaged with $\alpha=\frac{\gamma\beta}{2}\in (0,1)$. Moreover, $x^\star\in\mathcal{X}$ is a fixed point of $T$ if and only if $\nabla f(x^\star)=0$, which is equivalent to $x^\star$ being a minimizer of $f$. Hence, we can invoke the Krasnosel'ski\u{\i}--Mann theorem (without relaxation) to obtain weak convergence of the sequence $(x^{(i)})_{i\in\mathbb{N}}$ to a minimizer $x^\star$ of $f$.

We note that the relaxed version \eqref{eq0} with a fixed parameter $\rho$ of gradient descent simplifies to $x^{(i+1)}= x^{(i)} - \rho \gamma \nabla f (x^{(i)})$ with $\rho\gamma \in (0,\frac{2}{\beta})$, so that we can just set $\rho=1$ and keep $\gamma$ as the only parameter. How to choose the value of $\gamma$, then? As $\gamma$ tends to $\frac{2}{\beta}$, the operator $\mathrm{Id} -\gamma \nabla f$ is $\alpha$-averaged with $\alpha$ tending to 1, and we have seen that convergence is not guaranteed with merely nonexpansive operators. So, we can expect convergence to become slower as $\gamma$ tends to $\frac{2}{\beta}$; but this is true \emph{in the worst case}: if $f=\frac{1}{2}\|\cdot\|^2$, whose gradient $\nabla f=\mathrm{Id}$ is $\beta=1$-Lipschitz continuous, we have $x^{(i+1)}= (1-\gamma)x^{(i)}$, and $\gamma=\frac{1}{\beta}=1$ is the best choice, whereas convergence becomes slower and slower as $\gamma\rightarrow \frac{2}{\beta}=2$. However, a more insightful and realistic example is the least-squares function $f:x\mapsto \frac{1}{2}\|Ax-b\|^2$ for some bounded linear operator $A$ with $\|A\|=1$ and element $b$. Gradient descent is then equivalent to the classical Richardson iteration that solves the linear system $A^*Ax=A^*b$:
\begin{equation}
x^{(i+1)}=x^{(i)}-\gamma A^*(A x^{(i)}-b), 
\end{equation}
for some  $\gamma \in (0,2)$. The operator mapping $x^{(i)}$ to $x^{(i+1)}$ is $\gamma/2$-averaged, but it is also affine and quasicontractive, 
and the convergence rate is fully characterized by the eigenvalues of $A^*A$. Let $\lambda_{\min}>0$ be its smallest positive eigenvalue; we have assumed that the largest eigenvalue is 1. Then the best convergence rate is obtained for $\gamma=2/(1+\lambda_{\min})$, which, for $\lambda_{\min}$ small, is close to 2. In many applications,  the spectrum of $A^*A$ is  concentrated around 0, so that $\lambda_{\min}$, even if unknown, is likely to be small. 
In any case, there is no  interest in choosing $\gamma<1$. So, as a rule of thumb, in gradient descent one should try the two values $\gamma=\frac{1}{\beta}$ and $\gamma=\frac{1.9}{\beta}$. Similarly, in the Krasnosel'ski\u{\i}--Mann theorem, we recommend tuning the parameters in the algorithms with $\rho=1$ and then trying overrelaxation with $\rho=\frac{1.9}{2\alpha}$.\medskip

Let us now consider the minimization of a function $f\in\Gamma_0(\mathcal{X})$ whose set of minimizers is nonempty. The \textbf{proximal point algorithm}~\cite{roc76},  initialized at some $x^{(0)}\in\mathcal{X}$, consists in iterating, for $i\geq 0$,
\begin{equation}
x^{(i+1)}= \mathrm{prox}_{\gamma f} (x^{(i)})
\end{equation}
for some $\gamma>0$. Since, by definition of the proximity operator, $f(x^{(i+1)})\leq f(x^{(i)})-\frac{1}{2\gamma}\|x^{(i+1)}- x^{(i)}\|^2$, this is also a descent algorithm. We can easily show that 
\begin{equation}
x^{(i+1)}\in  x^{(i)} - \gamma \partial f (x^{(i+1)}),
\end{equation}
so that the proximal point algorithm looks similar to gradient descent, but with an implicit call to the subdifferential via the proximity operator. $\mathrm{prox}_{\gamma f}$ is firmly nonexpansive, and $x^\star\in\mathcal{X}$ is a fixed point of $\mathrm{prox}_{\gamma f}$ if and only if $0\in \partial f(x^\star)$, which is equivalent to $x^\star$ being a minimizer of $f$. So, we can again invoke the Krasnosel'ski\u{\i}--Mann theorem (without relaxation) to obtain weak convergence of the sequence $(x^{(i)})_{i\in\mathbb{N}}$ to a minimizer $x^\star$ of $f$. Note that if the proximal point algorithm can be applied with  every $\gamma>0$, by taking $\gamma$ very large, it essentially converges in one single iteration. So, it is rare to be able to minimize a function, whose minimizer is unknown, using the proximal point algorithm as is. However, we will see in the rest of the paper that proximal splitting algorithms are able to minimize a sum of functions using calls to their gradients or proximity operators.

\section{The Forward-Backward Algorithm}\label{sec2}

Let $\mathcal{Z}$ be a real Hilbert space.
Let $M:\mathcal{Z}\rightarrow 2^\mathcal{Z}$ be a set-valued maximally monotone operator and let $\bl{C}:\mathcal{Z}\rightarrow \mathcal{Z}$ be a $\xi$-cocoercive operator for some real $\xi>0$. Let us consider
the monotone inclusion 
\begin{equation}
    0\in Mz+\bl{C}z,\label{eq1}
\end{equation}
whose solution set is supposed nonempty.
Let $z^{(0)}\in\mathcal{Z}$ be some initial estimate of a solution and let $\gamma>0$ be some real parameter. The classical forward-backward iteration to find a solution of \eqref{eq1} consists in iterating
\begin{equation}
     z^{(i+1)}= J_{\gamma M}(z^{(i)}-\gamma \bl{C} z^{(i)}). \label{eq2}
\end{equation}
This method, proposed by Mercier~\cite{mer792}, was further developed by many authors~\cite{lio79, pas79,gab83,tse91,che941,com04, com05}.

A not-so-well-known extension of the forward-backward iteration consists in \emph{relaxing} it. 
Let $(\rho^{(i)})_{i\in\mathbb{N}}$ be a sequence of relaxation parameters. The iteration becomes:\\[-2mm]
\ceq{align}{
 &\,\mbox{\textbf{Forward-backward iteration for \eqref{eq1}}: for }i=0,1,\ldots,\notag\\[-1mm]
    &\left\lfloor
    \begin{array}{l}
    z^{(i+\frac{1}{2})}= J_{\gamma M}(z^{(i)}-\gamma \bl{C} z^{(i)})\\
     z^{(i+1)}=z^{(i)}+\rho^{(i)} (z^{(i+\frac{1}{2})}-z^{(i)}).
      \end{array}\right.\label{eq3}}%
To the best of our knowledge, Theorem 25.8 of \cite{bau11} is the first convergence result about the relaxed forward-backward algorithm, though with a smaller relaxation range than in \cref{lemma21} below.

We remark that the explicit mapping from $z^{(i)}$ to $z^{(i+\frac{1}{2})}$ in \eqref{eq3} can be equivalently written in the implicit form
\begin{equation}
    0\in M z^{(i+\frac{1}{2})} + \bl{C} z^{(i)} + \frac{1}{\gamma} (z^{(i+\frac{1}{2})} - z^{(i)}).\label{eq4}
\end{equation}

The now standard convergence result for the forward-backward iteration is the following~\cite[Lemma 4.4]{con13} \cite[Theorem 26.14]{bau17}:

\begin{lemma}[forward-backward algorithm \eqref{eq3}]\label{lemma21}
Let $z^{(0)}\in\mathcal{Z}$. Let $\gamma \in (0, 2\xi)$. Set $\delta=2-\gamma/(2\xi)$. Let $(\rho^{(i)})_{i\in\mathbb{N}}$ be a sequence in $[0,\delta]$ such that $\sum_{i\in \mathbb{N}} \rho^{(i)} (\delta-\rho^{(i)})=+\infty$. 
Then the sequence $(z^{(i)})_{i\in \mathbb{N}}$ defined by the iteration \eqref{eq3} converges weakly to a solution of \eqref{eq1}.
\end{lemma}

\begin{proof}
This result follows from the Krasnosel'ski\u{\i}--Mann theorem, since the operator $T=J_{\gamma M}(\mathrm{Id}-\gamma \bl{C})$, which maps $z^{(i)}$ to $z^{(i+\frac{1}{2})}$, is $(1/\delta)$-averaged as the composition of a firmly nonexpansive operator and a $\gamma/(2\xi)$-averaged operator, according to \eqref{eqcompav}.  Moreover, the fixed points of $T$ are the solutions to \eqref{eq1}, as is clearly visible in \eqref{eq4}.
\end{proof}\smallskip

We can allow $\gamma=2\xi$ in \cref{lemma21}, but in this case $\delta=1$ and we cannot set $\rho^{(i)}=1$. Since we are interested in overrelaxation, we write all theorems such that the default choice $\rho^{(i)}=1$ is allowed.

We note that in this algorithm, like in every algorithm of the paper, the sequence $(z^{(i+\frac{1}{2})})_{i\in \mathbb{N}}$ converges weakly to the same solution as $(z^{(i)})_{i\in \mathbb{N}}$. For any $i$, $z^{(i+\frac{1}{2})}$ may be a better estimate of the solution than $z^{(i)}$, since it is in the domain of $M$; that is, $Mz^{(i+\frac{1}{2})}\neq \emptyset$.

We can also remark that if $\gamma$ is close to $2\xi$, $\delta$ is close to 1, so that overrelaxation cannot be used.   This explains why the relaxed forward-backward iteration is not so well known.

Now, let $P$ be a self-adjoint, strongly positive, bounded linear operator on $\mathcal{Z}$. Clearly, solving \eqref{eq1} is equivalent to solving
\begin{equation}
0\in P^{-1}Mz + P^{-1}\bl{C}z.\label{eqfpb}
\end{equation}
Let $\mathcal{Z}_P$ be the Hilbert space obtained by endowing $\mathcal{Z}$ with the inner product $(x, x') \mapsto  \langle x, x'\rangle_P = \langle x, Px'\rangle$. Then $P^{-1}M$ is maximally monotone in $\mathcal{Z}_P$~\cite[Proposition 20.24]{bau17}.
However, the cocoercivity of $P^{-1}\bl{C}$ in $\mathcal{Z}_P$ has to be checked on a case-by-case basis. The \emph{preconditioned} forward-backward iteration to solve \eqref{eqfpb} is:
\ceq{align}{
&\,\mbox{\textbf{Preconditioned forward-backward iteration for \eqref{eq1}}: for }i=0,1,\ldots,\notag\\[-1mm]
    &\left\lfloor
    \begin{array}{l}
    z^{(i+\frac{1}{2})}= ( P^{-1}M+\mathrm{Id})^{-1}(z^{(i)}-P^{-1}\bl{C} z^{(i)})\\
     z^{(i+1)}=z^{(i)}+\rho^{(i)} (z^{(i+\frac{1}{2})}-z^{(i)}).
      \end{array}\right.
    \label{eq3p}}%
The corresponding convergence result follows.

\begin{theorem}[preconditioned forward-backward algorithm \eqref{eq3p}]\label{thm22}
Suppose that $P^{-1}\bl{C}$ is $\chi$-cocoercive in $\mathcal{Z}_P$ for some $\chi>\frac{1}{2}$. 
Set $\delta=2-1/(2\chi)$. Let $z^{(0)}\in\mathcal{Z}$ and let $(\rho^{(i)})_{i\in\mathbb{N}}$ be a sequence in $[0,\delta]$ such that $\sum_{i\in \mathbb{N}} \rho^{(i)} (\delta-\rho^{(i)})=+\infty$. 
Then the sequence $(z^{(i)})_{i\in \mathbb{N}}$ defined by the iteration \eqref{eq3p} converges weakly to a solution of \eqref{eq1}.
\end{theorem}

\begin{proof}This is a straightforward application of \cref{lemma21} in $\mathcal{Z}_P$ instead of $\mathcal{Z}$, with $\gamma=1$. Weak convergence in $\mathcal{Z}_P$ is equivalent to weak convergence in $\mathcal{Z}$, and the solution sets of \eqref{eq1} and \eqref{eqfpb} are the same.
\end{proof}\smallskip

We remark that the explicit mapping from $z^{(i)}$ to $z^{(i+\frac{1}{2})}$ in \eqref{eq3p} can be equivalently written in the implicit form:
\begin{equation}
    0\in M z^{(i+\frac{1}{2})} + \bl{C} z^{(i)} + P (z^{(i+\frac{1}{2})} - z^{(i)}).\label{eq4p}
\end{equation}
Also, we can write the iteration as $z^{(i+\frac{1}{2})}= ( M+P)^{-1} (P- \bl{C}) z^{(i)}$, an iteration studied with nonlinear kernels $P$ in \cite{gis21} and \cite{bui20}.
\bigskip

If $\bl{C}=0$, the forward-backward iteration reduces to the proximal point algorithm~\cite{roc76,bau17}; \br{we have discussed this algorithm in \cref{secadd} with proximity operators, but it generalizes to any resolvent.} 
 Let us give a formal convergence statement for the preconditioned proximal point algorithm. 
Let $M: \mathcal{Z} \rightarrow  2^\mathcal{Z} $ be a maximally monotone operator.
The problem is to solve
\begin{equation}
    0\in Mz,\label{eq7}
\end{equation}
whose solution set is supposed nonempty. 
Let $z^{(0)}\in\mathcal{Z}$ and let $(\rho^{(i)})_{i\in\mathbb{N}}$ be a sequence of relaxation parameters. The relaxed and
preconditioned proximal point algorithm consists in the iteration:
\ceq{align}{
&\,\mbox{\textbf{Proximal point iteration for \eqref{eq7}}: for }i=0,1,\ldots,\notag\\[-1mm]
    &\left\lfloor
    \begin{array}{l}
    z^{(i+\frac{1}{2})}= (P^{-1} M+\mathrm{Id})^{-1} z^{(i)}\\
     z^{(i+1)}=z^{(i)}+\rho^{(i)} (z^{(i+\frac{1}{2})}-z^{(i)}).
      \end{array}\right.
    \label{eq5}}%
The mapping from $z^{(i)}$ to $z^{(i+\frac{1}{2})}$ in \eqref{eq5} can be equivalently written as
\begin{equation}
    0\in M z^{(i+\frac{1}{2})}  + P (z^{(i+\frac{1}{2})} - z^{(i)}).\label{eq6}
\end{equation}
The convergence of the proximal point algorithm can be stated as 
 follows:
 
 \begin{theorem}[proximal point algorithm \eqref{eq5}]\label{thm23}
 Let $z^{(0)}\in\mathcal{Z}$ and let  $(\rho^{(i)})_{i\in\mathbb{N}}$ be a sequence in $[0,2]$ such that $\sum_{i\in \mathbb{N}} \rho^{(i)} (2-\rho^{(i)})=+\infty$. Then the sequence $(z^{(i)})_{i\in \mathbb{N}}$ defined by the iteration \eqref{eq5} converges weakly to a solution of \eqref{eq7}.
 \end{theorem}

\begin{proof}$P^{-1}M$ is maximally monotone in $\mathcal{Z}_P$~\cite[Proposition 20.24]{bau17}, so its resolvent $(P^{-1} M+\mathrm{Id})^{-1}$ is firmly nonexpansive in $\mathcal{Z}_P$. By virtue of the Krasnosel'ski\u{\i}--Mann theorem, 
$(z^{(i)})_{i\in \mathbb{N}}$ converges weakly in $\mathcal{Z}_P$ to some element $z^\star\in\mathcal{Z}$ with $0\in P^{-1} Mz^\star$, so that $0\in Mz^\star$.
\end{proof}

\subsection{The case where $\bl{C}$ is affine}

Let us continue with our analysis of the forward-backward iteration to solve \eqref{eq1}. In this section, we suppose that, in addition to being $\xi$-cocoercive, $\bl{C}$ is affine; that is, 
\begin{equation}
\bl{C}:z\in\mathcal{Z}\mapsto \bl{Q}z + c
\end{equation}
 for some self-adjoint, positive, nonzero, bounded linear operator $\bl{Q}$ on $\mathcal{Z}$ and some element $c\in\mathcal{Z}$. We have $\xi=1/\|\bl{Q}\|$. 
Now, we can write \eqref{eq4} as
\begin{equation}
    0\in (M+\bl{C}) z^{(i+\frac{1}{2})}  + P (z^{(i+\frac{1}{2})} - z^{(i)}),\label{eq8}
\end{equation}
where
\begin{equation}
P={\textstyle\frac{1}{\gamma}}\mathrm{Id}-\bl{Q}.
\end{equation}
Therefore, the forward-backward iteration \eqref{eq3} can be interpreted as a preconditioned proximal point iteration \eqref{eq5}, applied to find a zero of $M+\bl{C}$! Since $P$ must be strongly positive, we must have $\gamma \in (0,\xi)$, so that the admissible range for $\gamma$ is halved. But in return we have the larger range $(0,2)$ for relaxation. Hence, we have the following new convergence result:

 \begin{theorem}[forward-backward algorithm \eqref{eq3}, affine case]\label{thm24}
 Let $z^{(0)}\in\mathcal{Z}$, let $\gamma \in (0, \xi)$, and let $(\rho^{(i)})_{i\in\mathbb{N}}$ be a sequence in $[0,2]$ such that $\sum_{i\in \mathbb{N}} \rho^{(i)} (2-\rho^{(i)})=+\infty$.
Then the sequence $(z^{(i)})_{i\in \mathbb{N}}$ defined by the iteration \eqref{eq3} converges weakly to a solution of \eqref{eq1}.\end{theorem}

\begin{proof}In view of \eqref{eq6} and \eqref{eq8}, this is \cref{thm23} applied to the problem \eqref{eq7} with $M+\bl{C}$ as the maximally monotone operator.
\end{proof}\smallskip

Now, let us improve the range $(0, \xi)$ for $\gamma$ in  \cref{thm24} to $(0, \xi]$. Let $\gamma\in(0,\xi]$. We first recall that for every self-adjoint, positive, bounded linear operator $P$ on $\mathcal{Z}$, there exists a unique self-adjoint, positive, bounded linear operator on $\mathcal{Z}$, denoted by $\sqrt{P}$ and called the positive square root of $P$, such that $\sqrt{P}\sqrt{P}=P$. So, let $A$ be a bounded linear operator from $\mathcal{Z}$ to some real Hilbert space $\mathcal{A}$, such that $A^*A=\bl{Q}$.  $A=\sqrt{\bl{Q}}$ is a valid choice, but in some cases, $\bl{Q}=A^*A$ for some $A$ in the first place, such as in least-squares problems; see \eqref{eq12}.  In any case, we do not need to exhibit $A$; the fact that it exists is sufficient here. Furthermore, let $B$ be a bounded linear operator from some real Hilbert space $\mathcal{B}$ to $\mathcal{A}$, such that $AA^*+BB^*=(1/\gamma)\mathrm{Id}$. Since $\|AA^*\|=\|A^*A\|=\|\bl{Q}\|=1/\xi\leq 1/\gamma$,  $(1/\gamma)\mathrm{Id}-AA^*$ is positive and we can choose $B=\sqrt{(1/\gamma)\mathrm{Id}-AA^*}$. Again, we do not have to exhibit $B$; the fact that it exists is sufficient. Then we can rewrite the problem \eqref{eq1} as finding a pair $w=(z,b)\in\mathcal{W}=\mathcal{Z}\times\mathcal{B}$ solution to the monotone inclusion
\begin{equation}\label{eqab}
 \left( \begin{array}{c}
   0\\0
   \end{array}\right)\in
    \underbrace{\left( \begin{array}{c}
     A^*(Az+Bb) \\
      B^*(Az+Bb)
   \end{array}\right)}_{Sw}+ \underbrace{\left( \begin{array}{c}
    Mz +c\\
  (\mathcal{B}\mbox{ if }b= 0,\ \emptyset\mbox{ otherwise})
   \end{array}\right)}_{Nw}.
\end{equation}
The operators $S$ and $N$ defined in \eqref{eqab} are maximally monotone in $\mathcal{W}$. To solve a monotone inclusion involving two monotone operators, the Douglas--Rachford splitting algorithm~\cite{lio79,eck92,com04,sva11} is a natural choice. So, let us consider a general monotone inclusion in some arbitrary real Hilbert space $\mathcal{W}$, with nonempty solution set:
\begin{equation}\label{eqab2}
0\in Sw + Nw,
\end{equation}
for some maximally monotone operators $S$ and $N$ on $\mathcal{W}$. Let $w^{(0)}\in\mathcal{W}$ and $y^{(0)}\in\mathcal{W}$, let $\tau >0$ be a real parameter, and let $(\rho^{(i)})_{i\in\mathbb{N}}$ be a sequence of relaxation parameters.
The Douglas--Rachford iteration can be written in several ways, one of which is as follows:
\ceq{align}{
&\,\mbox{\textbf{Douglas--Rachford iteration for \eqref{eqab2}}: for }i=0,1,\ldots,\notag\\[-1mm]
    &\left\lfloor
    \begin{array}{l}
    y^{(i+\frac{1}{2})}= J_{S^{-1}/\tau} (y^{(i)}+w^{(i)}/\tau)\\
    w^{(i+\frac{1}{2})}= J_{\tau N} \big(w^{(i)}-\tau (2y^{(i+\frac{1}{2})}-y^{(i)})\big)\\
     y^{(i+1)}=y^{(i)}+\rho^{(i)} (y^{(i+\frac{1}{2})}-y^{(i)})\\
     w^{(i+1)}=w^{(i)}+\rho^{(i)} (w^{(i+\frac{1}{2})}-w^{(i)}).
      \end{array}\right.\label{eqdrm}}
We have the following convergence theorem~\cite[Theorem 26.11]{bau17}:

\begin{lemma}[Douglas--Rachford algorithm \eqref{eqdrm}]\label{lem25}
Let $w^{(0)}\in\mathcal{W}$ and $y^{(0)}\in\mathcal{W}$, let $\tau >0$, and let $(\rho^{(i)})_{i\in\mathbb{N}}$ be a sequence in $[0,2]$ such that $\sum_{i\in \mathbb{N}} \rho^{(i)} (2-\rho^{(i)})=+\infty$. 
Then the sequence $(w^{(i)})_{i\in \mathbb{N}}$ defined by the iteration \eqref{eqdrm} converges weakly to a solution of \eqref{eqab2}.
\end{lemma}

Thus, let us apply the Douglas--Rachford algorithm to \eqref{eqab} in the lifted space $\mathcal{W}=\mathcal{Z}\times\mathcal{B}$, with $\tau=\gamma$. For this, we need to compute the resolvents; we have \cite[equation~15]{oco18}:
\begin{align}
&J_{S^{-1}/\gamma}:(z,b)\in\mathcal{W}\mapsto (A^*a,B^*a)\ \mbox{with}\ a={\textstyle\frac{\gamma}{2}}(Az+Bb)\\
&J_{\gamma N}:(z,b)\in\mathcal{W}\mapsto \big(J_{\gamma M}(z-\gamma c),0\big)
\end{align}
Hence, assuming that $w^{(0)}=(z^{(0)},0)$ for some $z^{(0)}\in\mathcal{Z}$ and subsequently keeping only the first part of the variable $w$ as the variable $z$, the Douglas--Rachford iteration in this setting is:
\begin{align}
&\,\mbox{\textbf{Douglas--Rachford iteration for \eqref{eqab}}: for }i=0,1,\ldots,\notag\\[-1mm]
    &\left\lfloor
    \begin{array}{l}
    a^{(i+\frac{1}{2})}= \frac{\gamma}{2}\big(A(y_z^{(i)}+z^{(i)}/\gamma)+By_b^{(i)}\big)\\
    y_z^{(i+\frac{1}{2})}= A^* a^{(i+\frac{1}{2})}\\
    y_b^{(i+\frac{1}{2})}=B^* a^{(i+\frac{1}{2})}\\
    z^{(i+\frac{1}{2})}=J_{\gamma M}\big(z^{(i)}-\gamma (2y_z^{(i+\frac{1}{2})}-y_z^{(i)})-\gamma c\big)\\
     y_z^{(i+1)}=y_z^{(i)}+\rho^{(i)} (y_z^{(i+\frac{1}{2})}-y_z^{(i)})\\
      y_b^{(i+1)}=y_b^{(i)}+\rho^{(i)} (y_b^{(i+\frac{1}{2})}-y_b^{(i)})\\
      z^{(i+1)}=z^{(i)}+\rho^{(i)} (z^{(i+\frac{1}{2})}-z^{(i)}).
      \end{array}\right.\label{eqdrm2}
\end{align}
Suppose that $(y_z^{(0)},y_b^{(0)})=(A^*a^{(0)},B^*a^{(0)})$ for some $a^{(0)}\in \mathcal{A}$ and set  $a^{(i+1)}=a^{(i)}+\rho^{(i)} (a^{(i+\frac{1}{2})}-a^{(i)})$ for every $i\in\mathbb{N}$. Since for every $i\in\mathbb{N}$, $y_z^{(i)}=A^*a^{(i)}$ and $y_b^{(i)}=B^*a^{(i)}$, then $Ay_z^{(i)}+B y_b^{(i)}=(AA^*+BB^*)a^{(i)}=a^{(i)}/\gamma$. Therefore, we have
\begin{align}
z^{(i+\frac{1}{2})}&=J_{\gamma M}\big(z^{(i)}-\gamma (2y_z^{(i+\frac{1}{2})}-y_z^{(i)})-\gamma c\big)\\
&=J_{\gamma M}\big(z^{(i)}-\gamma^2A^*(Ay_z^{(i)}+Az^{(i)}/\gamma+By_b^{(i)})+\gamma y_z^{(i)}-\gamma c\big)\\
&=J_{\gamma M}\big(z^{(i)}-\gamma A^*Az^{(i)} -\gamma A^*a^{(i)}+\gamma y_z^{(i)}-\gamma c\big)\\
&=J_{\gamma M}\big(z^{(i)}-\gamma (\bl{Q}z^{(i)}+c)\big).
\end{align}
So, we can remove all variables but $z$ and we recover the forward-backward iteration \eqref{eq3}! It is remarkable that the forward-backward algorithm  can be viewed as an instance of the Douglas--Rachford algorithm, since they are both fundamental algorithms to find a zero of a sum of two monotone operators. This is the case only because the cocoercive operator $\bl{C}$ is supposed here to be affine, though.

Hence, as an application of \cref{lem25}, we have the following convergence result, which extends  \cref{thm24} with the possibility of setting $\gamma=\xi$:

\begin{theorem}[forward-backward algorithm \eqref{eq3}, affine case]\label{thm26}
Let $z^{(0)}\in\mathcal{Z}$, let $\gamma \in (0, \xi]$, and let $(\rho^{(i)})_{i\in\mathbb{N}}$ be a sequence in $[0,2]$ such that $\sum_{i\in \mathbb{N}} \rho^{(i)} (2-\rho^{(i)})=+\infty$. 
Then the sequence $(z^{(i)})_{i\in \mathbb{N}}$ defined by the iteration \eqref{eq3} converges \br{weakly} to a solution of \eqref{eq1}.\end{theorem}

Thus, if $\gamma=\xi$, we can set $\rho^{(i)}=1.4$ according to \cref{lemma21}; with \cref{thm26}, we can do better and set $\rho^{(i)}=1.9$. \br{So, which values of the parameters should be used in practice? As discussed in \cref{secadd}, 
our recommendation for a given practical problem is to try the following three settings: 

\noindent$\bullet\;$ $\gamma=\xi$ and $\rho^{(i)}=1$, 

\noindent$\bullet\;$ $\gamma=1.9\xi$ and $\rho^{(i)}=1$,

\noindent$\bullet\;$ $\gamma=\xi$ and $\rho^{(i)}=1.9$.

\noindent We illustrate these three choices within an image processing experiment in~\cref{figu0}.
}

\subsection{Applications to convex optimization}

Let $\mathcal{X}$ be a real Hilbert space. 
Let $\ff\in \Gamma_0(\mathcal{X})$ and let $\h:\mathcal{X}\rightarrow \mathbb{R}$ be a convex and 
differentiable function whose gradient  $\bl{\nabla h}$ is  $\beta$-Lipschitz continuous for some real $\beta> 0$. We consider the convex optimization problem:
\begin{equation}
    \minimize_{x\in\mathcal{X}}\, \dg{f(x)} + \bl{h(x)}.\label{eq10}
\end{equation}
From Fermat's rule, 
the problem \eqref{eq10} is equivalent to \eqref{eq1} with $M=\dg{\partial f}$, which is maximally monotone, and $\bl{C}=\bl{\nabla h}$, which is $\xi$-cocoercive, with $\xi=1/\beta$ 
(and, as mentioned in \cref{secadd}, a zero of $\dg{\partial f}+\bl{\nabla h}$ is supposed to exist; equivalently here,  a minimizer of  $\dg{f} + \bl{h}$ is supposed to exist). So, it is natural to use the forward-backward iteration to solve \eqref{eq10}. 
Let $\gamma>0$, let $x^{(0)}\in\mathcal{X}$ and let $(\rho^{(i)})_{i\in\mathbb{N}}$ be a sequence of relaxation parameters. The 
forward-backward iteration, also called the proximal gradient algorithm, is:
\ceq{align}{
&\,\mbox{\textbf{Forward-backward iteration for \eqref{eq10}}: for }i=0,1,\ldots,\notag\\[-1mm]
    &\left\lfloor
    \begin{array}{l}
    x^{(i+\frac{1}{2})}= \mathrm{prox}_{\gamma \ff}\big(x^{(i)}-\gamma \bl{\nabla h} (x^{(i)})\big)\\
     x^{(i+1)}=x^{(i)}+\rho^{(i)} (x^{(i+\frac{1}{2})}-x^{(i)}).
      \end{array}\right.
    \label{eq9}}%
As a direct consequence of \cref{lemma21}, we have:

\begin{theorem}[forward-backward algorithm \eqref{eq9}]\label{thm27}
Let $x^{(0)}\in\mathcal{X}$ and let $\gamma \in (0, 2/\beta)$. Set $\delta=2-\gamma\beta/2$. Let $(\rho^{(i)})_{i\in\mathbb{N}}$ be a sequence in $[0,\delta]$ such that $\sum_{i\in \mathbb{N}} \rho^{(i)} (\delta-\rho^{(i)})=+\infty$. 
Then the sequence $(x^{(i)})_{i\in \mathbb{N}}$ defined by the iteration \eqref{eq9} converges weakly to a solution of \eqref{eq10}.\end{theorem}

Now, let us focus on the case where $\h$ is quadratic:
\begin{equation}
\h:x\mapsto {\textstyle\frac{1}{2}}\langle x,\bl{Q}x\rangle + \langle x,c\rangle +t \label{eq11}
\end{equation}
for some self-adjoint, positive, nonzero, bounded linear operator $\bl{Q}$ on $\mathcal{X}$, some $c\in\mathcal{X}$, and some $t\in\mathbb{R}$. We omit the constant $t$ in the notations of quadratic functions in what follows, since it does not play any role, so we can assume that $t=0$.

A very common example is a least-squares penalty, used in particular to solve inverse problems~\cite{cha16}; 
that is,
\begin{equation}
\h:x\mapsto {\textstyle\frac{1}{2}}\|Ax-y\|^2,\label{eq12}
\end{equation}
for some bounded linear operator $A$ from $\mathcal{X}$ to some real Hilbert space $\mathcal{Y}$ and some $y\in\mathcal{Y}$.
Clearly, \eqref{eq12} is an instance of \eqref{eq11}, with $\bl{Q}=A^*A$ and $c=-A^*y$. 

We have, for every $x\in\mathcal{X}$,
\begin{equation}
    \bl{\nabla h}(x)=\bl{Q}x+c, 
\end{equation}
so $\beta=\|\bl{Q}\|$. Suppose that $\gamma \in (0, 1/\beta)$ and set $P={\textstyle\frac{1}{\gamma}}\mathrm{Id}-\bl{Q}$, which is strongly positive. Then  the update step in \eqref{eq9} can be written as
\begin{equation}
    x^{(i+\frac{1}{2})}=\argmin_{x\in\mathcal{X}} \big( \dg{f(x)}+ \bl{h(x)}+{\textstyle\frac{1}{2}}\|x-x^{(i)}\|_P^2\big),
\end{equation}
where we introduce the norm $\|\cdot\|_P:x\mapsto \sqrt{\langle x,Px\rangle}$. So, $x^{(i+\frac{1}{2})}$ can be viewed as the result of applying the proximity operator of $\ff+\h$ with the  preconditioned norm $\|\cdot\|_P$.  Thus, convergence follows from  \cref{thm24}. Furthermore, as a direct consequence of \cref{thm26}, we have:

\begin{theorem}[forward-backward algorithm \eqref{eq9}, quadratic case]\label{thm28}
Suppose that $\h$ is quadratic. Let $x^{(0)}\in\mathcal{X}$, let $\gamma \in (0, 1/\beta]$, and let  $(\rho^{(i)})_{i\in\mathbb{N}}$ be a sequence in $[0,2]$ such that $\sum_{i\in \mathbb{N}} \rho^{(i)} (2-\rho^{(i)})=+\infty$. 
Then the sequence $(x^{(i)})_{i\in \mathbb{N}}$ defined by the iteration \eqref{eq9} converges weakly to a solution of \eqref{eq10}.\end{theorem}

\section{The Loris--Verhoeven Algorithm}\label{sec3}

Let $\mathcal{X}$ and $\mathcal{U}$ be two real Hilbert spaces. Let $\g\in \Gamma_0(\mathcal{U})$ and let $\h:\mathcal{X}\rightarrow \mathbb{R}$ be a convex and 
differentiable function with  $\beta$-Lipschitz continuous gradient for some real $\beta> 0$. Let $L:\mathcal{X}\rightarrow \mathcal{U}$ be a bounded linear operator.

Often, the template problem \eqref{eq10} of minimizing the sum of two functions is too simple and we would like, instead, to
\begin{equation}
    \minimize_{x\in\mathcal{X}}\, \dr{g(Lx)} + \bl{h(x)}.\label{eq13}
\end{equation}
We assume that there is no simple way to compute the proximity operator of $\g\circ L$. The corresponding monotone inclusion, which we  will actually solve, is
\begin{equation}
    0 \in L^*\dr{\partial g} (Lx) + \bl{\nabla h}(x).
\end{equation}
To bypass the annoying operator $L$, we introduce an auxiliary variable $u\in\dr{\partial g} (Lx)$, which shall be called the \emph{dual} variable, so that the problem now consists in finding $x\in\mathcal{X}$ and $u\in\mathcal{U}$ such that
\begin{equation}\label{eqlvinc0}
  \left\{\begin{array}{l}
  u\in\dr{\partial g} (Lx)\\
  0 = L^*u + \bl{\nabla h}(x)
  \end{array}\right. .
\end{equation}
The interest in increasing the dimension of the problem is that we obtain a system of two monotone inclusions, which are decoupled: $\dr{\partial g}$ and $\bl{\nabla h}$ appear separately in the two inclusions. 
So, equivalently, the problem is to find 
 a pair of objects $z=(x,u)$   in 
$\mathcal{Z}=\mathcal{X}\times\mathcal{U}$
such that
\begin{equation}\label{variational_prim_dual_prob}
   \left( \begin{array}{c}
   0\\0
   \end{array}\right)\in
    \underbrace{\left( \begin{array}{c}
     L^*u \\
   -Lx + (\dr{\partial g})^{-1}u
   \end{array}\right)}_{Mz}+ \underbrace{\left( \begin{array}{c}
     \bl{\nabla h}(x) \\
  0
   \end{array}\right)}_{\bl{C}z} .
\end{equation}
The operator $M:\mathcal{Z}\rightarrow 2^\mathcal{Z},(x,u)\mapsto (L^*u, -Lx + (\dr{\partial g})^{-1}u)$ is maximally monotone~\cite[Proposition 26.32 (iii)]{bau17} and $\bl{C}:\mathcal{Z}\rightarrow\mathcal{Z},(x,u)\mapsto (\bl{\nabla h}(x),0)$ is $\xi$-cocoercive, with $\xi=1/\beta$. Thus, it is natural to think of using the forward-backward iteration to solve the problem \eqref{variational_prim_dual_prob}. However, to make the resolvent of $M$ computable with the proximity operator of $\g$, we need  
preconditioning. The solution consists in the following iteration, which we first write in implicit form:
\begin{align}\label{var_inclusion}
	&\left( \begin{array}{c}
	0  \\ 
	0  \\
	\end{array} \right) 
	\in 
	\underbrace{\left( \begin{array}{c}
	L^* u^{(i+\frac{1}{2})} \\ 
	-L x^{(i+\frac{1}{2})}+(\dr{\partial g})^{-1} u^{(i+\frac{1}{2})} \\
	\end{array} \right)}_{Mz^{(i+\frac{1}{2})}}+
	\underbrace{\left( \begin{array}{c}
	\bl{\nabla h}(x^{(i)})  \\ 
	0  \\
	\end{array} \right)}_{\bl{C}z^{(i)}}\\
	&\qquad\qquad	\qquad+
	\underbrace{\left( \begin{array}{cc} 
	\frac{1}{\tau}\mathrm{Id} & 0 \\ 
	0 & \frac{1}{\sigma}\mathrm{Id}-\tau LL^* \\
	\end{array} \right)}_{P}
	\underbrace{\left( \begin{array}{c}
	x^{(i+\frac{1}{2})}-x^{(i)}  \\ 
	u^{(i+\frac{1}{2})}-u^{(i)}\\
	\end{array} \right)}_{z^{(i+\frac{1}{2})}-z^{(i)}},\notag
\end{align}
where $\tau>0$ and $\sigma>0$ are two real parameters, $z^{(i)}=(x^{(i)},u^{(i)})$, and $z^{(i+\frac{1}{2})}=(x^{(i+\frac{1}{2})},u^{(i+\frac{1}{2})})$. It is not straightforward to see that this yields an explicit iteration. The key is to remark that we have
\begin{equation}
    	x^{(i+\frac{1}{2})}=x^{(i)}-\tau 	\bl{\nabla h}(x^{(i)}) - \tau L^* u^{(i+\frac{1}{2})},\label{eq15}
\end{equation}
so that we can update the primal variable $x^{(i+\frac{1}{2})}$, once the dual variable $u^{(i+\frac{1}{2})}$ is available. So, the first step of the algorithm is to construct $u^{(i+\frac{1}{2})}$. It depends on $Lx^{(i+\frac{1}{2})}$, which is not yet available, but using \eqref{eq15}, we can express it using $x^{(i)}$ and $LL^*u^{(i+\frac{1}{2})}$. This last term is canceled in the preconditioner $P$ to make the update explicit. That is,
\begin{align}
    		&0\in -L x^{(i+\frac{1}{2})}+(\dr{\partial g})^{-1} u^{(i+\frac{1}{2})}+({\textstyle\frac{1}{\sigma}}\mathrm{Id}-\tau LL^*)(u^{(i+\frac{1}{2})}-u^{(i)})\\
\Leftrightarrow{}&0\in -L x^{(i)} \!+\! \tau 	L\bl{\nabla h}(x^{(i)})\!+\!\tau LL^* u^{(i+\frac{1}{2})}
\!+\!(\dr{\partial g})^{-1} u^{(i+\frac{1}{2})}\!+\!({\textstyle\frac{1}{\sigma}}\mathrm{Id}\!-\!\tau LL^*)(u^{(i+\frac{1}{2})}\!-\!u^{(i)})\\
\Leftrightarrow{}&0\in -L x^{(i)} + \tau 	L\bl{\nabla h}(x^{(i)})+\tau LL^* u^{(i)}
+(\dr{\partial g})^{-1} u^{(i+\frac{1}{2})}+{\textstyle\frac{1}{\sigma}}
(u^{(i+\frac{1}{2})}-u^{(i)})\\
\Leftrightarrow{}&\big(\sigma(\dr{\partial g})^{-1}+\mathrm{Id}\big) u^{(i+\frac{1}{2})}
\ni \sigma L x^{(i)} - \tau\sigma 	L\bl{\nabla h}(x^{(i)})-\tau \sigma LL^* u^{(i)}+u^{(i)}\\
\Leftrightarrow{}&u^{(i+\frac{1}{2})}
= \mathrm{prox}_{\sigma \dr{g^*}}\Big(\sigma L \big(x^{(i)} - \tau\bl{\nabla h}(x^{(i)})\big)+u^{(i)}-\tau \sigma LL^* u^{(i)}\Big).
\end{align}
We note that \eqref{eqlvinc0} is equivalent to 
\begin{equation}
  \left\{\begin{array}{l}
  Lx\in (\dr{\partial g})^{-1}u\\
  x \in (\bl{\nabla h})^{-1}(-L^*u) 
  \end{array}\right. ,
\end{equation}
which implies that $0\in (\dr{\partial g})^{-1}u-L(\bl{\nabla h})^{-1}(-L^*u)$; this is the first-order characterization  of the convex optimization problem
\begin{equation}
    \minimize_{u\in\mathcal{U}}\, \dr{g^*(u)} + \bl{h^*(-L^*u)},\label{eq21}
\end{equation}
which is called the \emph{dual} problem to the \emph{primal} problem \eqref{eq13}. So, if a pair $(x,u)\in\mathcal{X}\times\mathcal{U}$ is a solution to \eqref{variational_prim_dual_prob}, then $x$ is a solution to \eqref{eq13} and $u$ is a solution to \eqref{eq21}.

Thus, let $\tau>0$ and $\sigma>0$, let $x^{(0)}\in\mathcal{X}$ and 
$u^{(0)}\in\mathcal{U}$, and let $(\rho^{(i)})_{i\in\mathbb{N}}$ be a sequence of relaxation parameters. The primal-dual forward-backward iteration, which we call the Loris--Verhoeven iteration is:
\ceq{align}{\label{algorithm_1}
&\,\mbox{\textbf{Loris--Verhoeven iteration for \eqref{eq13} and \eqref{eq21}}: for }i=0,1,\ldots,\notag\\[-1mm]
    &\left\lfloor
    \begin{array}{l}
    u^{(i+\frac{1}{2})}
= \mathrm{prox}_{\sigma \dr{g^*}}\Big(u^{(i)}+\sigma L \big(x^{(i)} - \tau\bl{\nabla h}(x^{(i)})-\tau L^* u^{(i)}\big)\Big)\\
x^{(i+1)}=x^{(i)}-\rho^{(i)}\tau \big(
\bl{\nabla h}(x^{(i)}) + L^* u^{(i+\frac{1}{2})}
\big)\\
u^{(i+1)}=u^{(i)}+\rho^{(i)} (u^{(i+\frac{1}{2})}-u^{(i)}).
\end{array}\right.}%
This algorithm was first proposed by Loris and Verhoeven in the case where $\h$ is a least-squares term~\cite{lor11}. 
It was then rediscovered several times and named \emph{Primal-Dual Fixed-Point algorithm based on the Proximity Operator}  (PDFP2O)~\cite{che13} or \emph{Proximal Alternating Predictor-Corrector} (PAPC) algorithm~\cite{dro15}. The above interpretation of the algorithm as a primal-dual forward-backward iteration has been presented in~\cite{com14}.

We note that the Loris--Verhoeven iteration can be written in such a way that there is only one call to $\bl{\nabla h}$ and $L^*$ per iteration; see the form \eqref{algorithm_1p}, for instance.

As an application of \cref{thm22}, we obtain the following convergence result:

\begin{theorem}[Loris--Verhoeven algorithm \eqref{algorithm_1}]\label{thm31}Let $x^{(0)}\in\mathcal{X}$ and 
$u^{(0)}\in\mathcal{U}$. Let $\tau \in (0, 2/\beta)$ and $\sigma>0$ be such that $\sigma\tau\|L\|^2<1$. Set $\delta=2-\tau\beta/2$. Let $(\rho^{(i)})_{i\in\mathbb{N}}$ be a sequence in $[0,\delta]$ such that $\sum_{i\in \mathbb{N}} \rho^{(i)} (\delta-\rho^{(i)})=+\infty$.
Then the sequences $(x^{(i)})_{i\in \mathbb{N}}$ and $(u^{(i)})_{i\in \mathbb{N}}$ defined by the iteration \eqref{algorithm_1} converge weakly to a solution of \eqref{eq13} and a solution of \eqref{eq21}, respectively.\end{theorem}

\begin{proof}In view of \eqref{var_inclusion} and \eqref{eq4p}, this is \cref{thm22} applied to the problem \eqref{variational_prim_dual_prob}. For this, $P$ must be strongly positive, which is the case if and only if $\sigma\tau\|L\|^2<1$.
Moreover, $P^{-1}\bl{C}$ is $1/(\tau\beta)$-cocoercive in $\mathcal{Z}_P$.
\end{proof}\smallskip

The following result makes it possible to have $\sigma\tau\|L\|^2=1$; it is a consequence of the analysis by O'Connor and Vandenberghe~\cite{oco18} of the PD3O algorithm~\cite{yan18} discussed in \cref{sec7}, 
of which the Loris--Verhoeven algorithm is a particular case. See also Theorems 3.4 and 3.5 in \cite{che13} for the same result, but without relaxation.

\begin{theorem}[Loris--Verhoeven algorithm \eqref{algorithm_1}]\label{thm32}
Suppose that $\mathcal{X}$ and $\mathcal{U}$ are of finite dimension. 
Let $x^{(0)}\in\mathcal{X}$ and 
$u^{(0)}\in\mathcal{U}$.
Let $\tau \in (0, 2/\beta)$ and  $\sigma>0$ be such that $\sigma\tau\|L\|^2\leq 1$. Set $\delta=2-\tau\beta/2$. Let $(\rho^{(i)})_{i\in\mathbb{N}}$ be a sequence in $[0,\delta]$ such that $\sum_{i\in \mathbb{N}} \rho^{(i)} (\delta-\rho^{(i)})=+\infty$.
Then the sequences $(x^{(i)})_{i\in \mathbb{N}}$ and $(u^{(i)})_{i\in \mathbb{N}}$ defined by the iteration \eqref{algorithm_1} converge to a solution of \eqref{eq13} and a solution of \eqref{eq21}, respectively.\end{theorem}

Note that the finite dimension assumption  is not necessary in the proof technique \cite{oco18}, so it could be removed.

Thus, in practice, one can keep $\tau$ as the only parameter to tune and set $\sigma=1/(\tau\|L\|^2)$. For primal-dual algorithms whose iteration satisfies an inclusion like \eqref{var_inclusion}, we want $P$ to be as close to $0$ as possible, broadly speaking. We also want the proximal parameters, here $\sigma$ and $\tau$, to be as large as possible. But these objectives are antagonistic: if $\tau$ is large, $\sigma$ must be small. Nevertheless, for a given value of $\tau$, setting $\sigma$ to the largest possible value $1/(\tau\|L\|^2)$ is the best we can do.\medskip
	
If $\mathcal{X}=\mathcal{U}$, $L=\mathrm{Id}$, and we set $\sigma=1/\tau$, the Loris--Verhoeven iteration becomes:
\begin{align}\label{LV_LId}
    &\left\lfloor
    \begin{array}{l}
    u^{(i+\frac{1}{2})}
= \mathrm{prox}_{\dr{g^*}/\tau}\big(x^{(i)}/\tau - \bl{\nabla h}(x^{(i)})\big)\\
x^{(i+\frac{1}{2})}=x^{(i)}-\tau 
\bl{\nabla h}(x^{(i)}) -\tau u^{(i+\frac{1}{2})}\\
\hphantom{x^{(i+1)}}=\mathrm{prox}_{\tau \g}\big(x^{(i)} - \tau\bl{\nabla h}(x^{(i)})\big)\\
x^{(i+1)}=x^{(i)}+\rho^{(i)} (x^{(i+\frac{1}{2})}-x^{(i)})\\
u^{(i+1)}=u^{(i)}+\rho^{(i)} (u^{(i+\frac{1}{2})}-u^{(i)}).
\end{array}\right.
\end{align}
So, we can discard the dual variable and recover the forward-backward iteration \eqref{eq9}. It is interesting that in this primal algorithm, there is an implicit dual variable $u^{(i+\frac{1}{2})}=-\bl{\nabla h}(x^{(i)})+(x^{(i)}-x^{(i+\frac{1}{2})})/\tau$, which converges to a solution of the dual problem; that is, to a minimizer of $\dr{g^*}(u)+\bl{h^*}(-u)$.\medskip

Again, let us focus on the case where $\h$ is quadratic; that is,
$\h:x\mapsto {\textstyle\frac{1}{2}}\langle x,\bl{Q}x\rangle + \langle x,c\rangle$  
for some self-adjoint, positive, nonzero, bounded linear operator $\bl{Q}$ on $\mathcal{X}$,  and some $c\in\mathcal{X}$. We have $\beta=\|\bl{Q}\|$. We can rewrite the primal-dual inclusion \eqref{var_inclusion}, which characterizes the Loris--Verhoeven iteration  \eqref{algorithm_1}, as
\begin{equation}\label{var_inclusion2}
	\left(\!\begin{array}{c}
	0  \\ 
	0  \\
	\end{array}\!\right) 
	\!\in\! 
	\left(\!\begin{array}{c}
\bl{\nabla h}(x^{(i+\frac{1}{2})})+	L^* u^{(i+\frac{1}{2})} \\ 
	-L x^{(i+\frac{1}{2})}+(\dr{\partial g})^{-1} u^{(i+\frac{1}{2})} \\
	\end{array}\!\right)
	\!+\!
	\left(\!\begin{array}{cc} 
	\frac{1}{\tau}\mathrm{Id}\!-\!\bl{Q} & 0 \\ 
	0 & \frac{1}{\sigma}\mathrm{Id}\!-\!\tau LL^* \\
	\end{array}\!\right)\!
	\left(\!\begin{array}{c}
	x^{(i+\frac{1}{2})}-x^{(i)}  \\ 
	u^{(i+\frac{1}{2})}-u^{(i)}\\
	\end{array}\!\right)\!.
\end{equation}
We recognize the structure of a primal-dual proximal point algorithm, so we can apply  \cref{thm23}. We can do better and apply \cref{thm26} to obtain convergence under the conditions $\tau \in (0, 1/\beta]$ and $\sigma\tau\|L\|^2<1$.
Indeed, we want to solve $0\in P^{-1}Mz + P^{-1}\bl{C}z$ in $\mathcal{Z}_P$, with $P$ defined in \eqref{var_inclusion}. We have $P^{-1}\bl{C}:(x,u)\mapsto (\tau\bl{Q}x+\tau c,0)$, which is affine and $1/(\tau\beta)$-cocoercive in $\mathcal{Z}_P$. Hence, \cref{thm26} can be applied with $\gamma=1$.

We can do even better and apply \cref{thm52} below, to obtain the condition $\sigma\tau\|L\|^2\leq 1$ instead of $\sigma\tau\|L\|^2< 1$, since according to \eqref{eqlvgcp1}--\eqref{eqlvgcp2}, the Loris--Verhoeven algorithm with $\h$ quadratic can be viewed as a primal-dual Douglas--Rachford algorithm. Hence, we have:

\begin{theorem}[Loris--Verhoeven algorithm \eqref{algorithm_1}, quadratic case]\label{thm33}Suppose that $\h$ is quadratic. Let $x^{(0)}\in\mathcal{X}$ and 
$u^{(0)}\in\mathcal{U}$. Let $\tau \in (0, 1/\beta]$ and $\sigma>0$ be such that $\sigma\tau\|L\|^2\leq 1$. Let $(\rho^{(i)})_{i\in\mathbb{N}}$ be a sequence in $[0,2]$ such that $\sum_{i\in \mathbb{N}} \rho^{(i)} (2-\rho^{(i)})=+\infty$. 
Then the sequences $(x^{(i)})_{i\in \mathbb{N}}$ and $(u^{(i)})_{i\in \mathbb{N}}$ defined by the iteration \eqref{algorithm_1} converge weakly to a solution of \eqref{eq13} and a solution of \eqref{eq21}, respectively.\end{theorem}

\begin{figure}[t]
\centering
\begin{tabular}{ccc}
\includegraphics[scale=0.53]{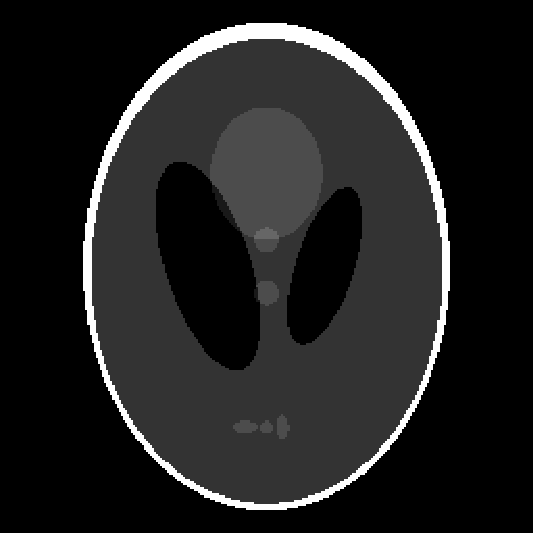}&\includegraphics[scale=0.53]{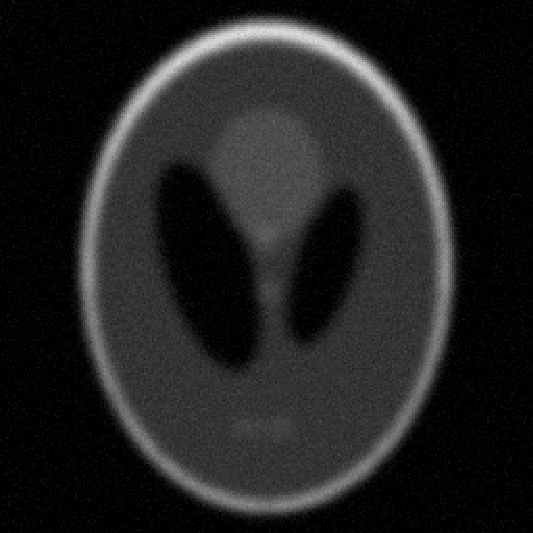}&\includegraphics[scale=0.53]{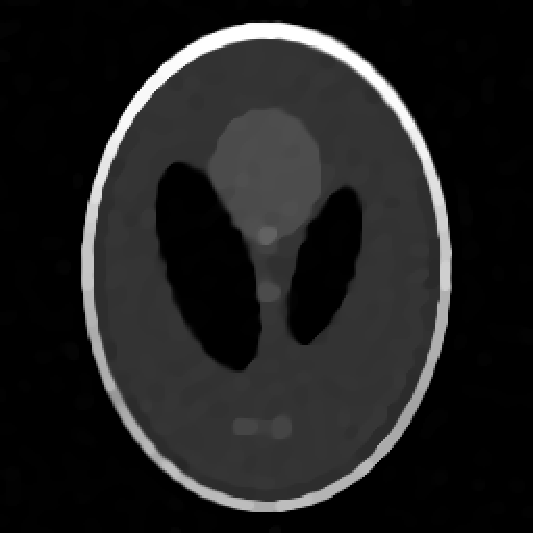}\\
(a)&(b)&(c)
\end{tabular}
\begin{tabular}{cc}
\hspace{-4.5mm}\includegraphics[scale=0.76]{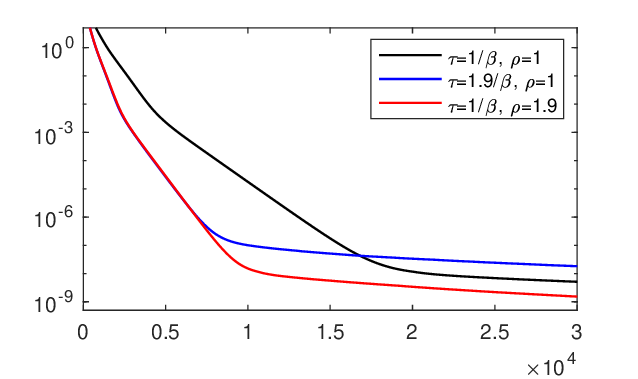}\hspace{-4mm}&\hspace{-4mm}\includegraphics[scale=0.76]{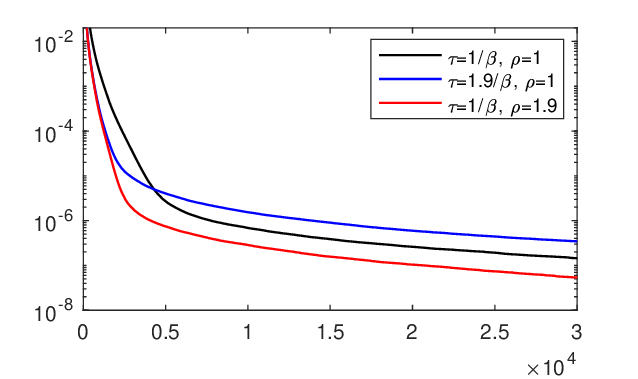}\hspace{-7mm}{\,}\\
$\!\!\!\!\!$(d) $\|x^{(i)}-x^\star\|^2$ w.r.t. $i\geq 1$&\ \ (e) $\Psi(x^{(i)})-\Psi(x^\star)$ w.r.t. $i\geq 1$
\end{tabular}
\caption{\br{Image deblurring experiment: the goal is to estimate the reference image, which is the classical Shepp--Logan phantom  (a), given the degraded image corrupted by blur and noise (b). The restored image (c) is the solution of a regularized least-squares problem of the form  \eqref{eq13}. The Loris--Verhoeven algorithm \eqref{algoCP} is used with $\sigma=1/(8\tau)$, according to \cref{thm32} and \cref{thm33}, since $\|L\|^2\leq 8$. The squared distance to the solution $x^\star$ and the decay of the objective function $\Psi$ along the iterations are shown in (d) and (e), respectively, for different values of $\tau$ and $\rho^{(i)}$ constant and equal to $1$ (no relaxation) or $1.9$ (overrelaxation). The benefit of overrelaxation (red curve), as allowed in \cref{thm33},  is clearly visible.}}
\label{figu0}
\end{figure}

In \cref{figu0}, we illustrate how overrelaxation can accelerate the Loris--Verhoeven algorithm when $\h$ is quadratic. We consider the problem of image deblurring, which consists in restoring an image $y$ corrupted by Gaussian blur and Gaussian noise, shown in \cref{figu0}~(b). We formulate deblurring as an optimization problem of the form \eqref{eq13}; that is, we aim at minimizing  $\Psi = \h +  \g\circ L$, where $\h:x\mapsto {\textstyle\frac{1}{2}}\|Ax-y\|^2$ is a least-squares data-fidelity term, with $A$ the convolution with a lowpass filter ($\beta=\|A\|=1$), 
and $\g\circ L$ is the ``isotropic'' total variation~\cite{rud92}: $L$ maps the image to the vector field made of all horizontal and vertical finite differences, and $ \dr{g}$ is $\lambda$ times the $l_{1,2}$ norm (see \cite{cha16,con17} for details), where $\lambda=0.002$ is the regularization parameter. So, we use the Loris--Verhoeven algorithm to solve the problem; a solution $x^\star$ (which seems to be unique and is estimated by $x^{(i)}$ with $i=10^6$) is shown in \cref{figu0}~(c). We can see from the plots of $\|x^{(i)}-x^\star\|^2$ and $\Psi(x^{(i)})$ that the choice $\tau=\frac{1}{\beta}$, $\sigma\tau\|L\|^2=1$, and $\rho^{(i)}=1.9$, as allowed by \cref{thm33}, speeds up the convergence in comparison with the standard choices allowed by \cref{thm32}.

\subsection{The Primal-Dual Fixed-Point (PDFP) algorithm}

Let us extend the problem \eqref{eq13} to
\begin{equation}
    \minimize_{x\in\mathcal{X}}\, \dg{f(x)}+\dr{g(Lx)}  + \bl{h(x)},\label{eq13e}
\end{equation}
using the same assumptions as before in this section, with an additional function  $\ff\in\Gamma_0(\mathcal{X})$. The dual problem is
\begin{equation}
    \minimize_{u\in\mathcal{U}}\, (\ff+\h)^*(-L^*u) + \dr{g^*(u)}.\label{eq21e}
\end{equation}
Can we modify the Loris--Verhoeven algorithm to handle this extended problem? We will see one way to do this in \cref{sec7}, with the PD3O algorithm. Meanwhile, let us extend the Loris--Verhoeven algorithm by keeping its primal-dual forward-backward structure. For this, we suppose  that for any $\gamma>0$, the proximity operator of $\gamma \ff$ is affine; that is, there exist some self-adjoint, positive, bounded linear operator $\dg{R_{\gamma}}$ on $\mathcal{X}$, with $\|\dg{R_{\gamma}}\|\leq 1$, and some element $a_\gamma\in\mathcal{X}$, such that
\begin{equation}
   \mathrm{prox}_{\gamma \ff}:x\in\mathcal{X}\mapsto \dg{R_{\gamma}}x+ a_{\gamma}.
\end{equation}
We mention two cases of practical interest, where the proximity operator is affine. First, this is the case if $\ff$ is quadratic: $\ff:x\in\mathcal{X}\mapsto \frac{1}{2}\langle x,\dg{Q'}x\rangle+\langle x,c'\rangle+t'$ for some self-adjoint, positive, bounded linear operator $\dg{Q'}$,  some $c'\in\mathcal{X}$,  and some $t'\in\mathbb{R}$, with  $\dg{R_{\gamma}}=(\gamma \dg{Q'}+\mathrm{Id})^{-1}$ and $a_{\gamma}=-(\gamma \dg{Q'}+\mathrm{Id})^{-1}c'$. Second, this is the case if $\ff$  is the indicator function of an affine subspace: $\ff:x\mapsto (0$ if $x\in\mathcal{A}$, $+\infty$ otherwise$)$, where $\mathcal{A}$ is a closed affine subspace of $\mathcal{X}$; for instance, $\mathcal{A}=\{x\in\mathcal{X} : Ax=y\}$ for some bounded linear operator $A$ and some element $y$ in the range of $A$; then $\mathrm{prox}_{\gamma \ff}$ is the projector onto $\mathcal{A}$, which does not depend on $\gamma$.

Let $\tau>0$ and $\sigma>0$. The extended Loris--Verhoeven iteration, written in implicit form, is
\begin{align}\label{var_inclusione}
	\left(\!\begin{array}{c}
	0  \\ 
	0  \\
	\end{array}\!\right) 
	&\in 
	\underbrace{\left( \begin{array}{c}
	\dg{\partial f}(x^{(i+\frac{1}{2})})+ L^* u^{(i+\frac{1}{2})} \\ 
	-L x^{(i+\frac{1}{2})}+(\dr{\partial g})^{-1} u^{(i+\frac{1}{2})} \\
	\end{array} \right)}_{Mz^{(i+\frac{1}{2})}}\\
	&\quad\ +
	\underbrace{\left( \begin{array}{c}
	\bl{\nabla h}(x^{(i)})  \\ 
	0  \\
	\end{array} \right)}_{\bl{C}z^{(i)}}
	+
	\underbrace{\left( \begin{array}{cc} 
	\frac{1}{\tau}\mathrm{Id} & 0 \\ 
	0 & \frac{1}{\sigma}\mathrm{Id}-\tau L\dg{R_\tau}L^* \\
	\end{array} \right)}_{P}
	\underbrace{\left( \begin{array}{c}
	x^{(i+\frac{1}{2})}-x^{(i)}  \\ 
	u^{(i+\frac{1}{2})}-u^{(i)}\\
	\end{array} \right)}_{z^{(i+\frac{1}{2})}-z^{(i)}}.\notag
\end{align}
This corresponds to the iteration:
\ceq{align}{\label{algorithm_1e}
&\,\mbox{\textbf{PDFP iteration for \eqref{eq13e} and \eqref{eq21e}}: for }i=0,1,\ldots,\notag\\[-1mm]
    &\left\lfloor
    \begin{array}{l}
     \tilde{x}^{(i)}
= \mathrm{prox}_{\tau \ff}\big(x^{(i)}-\tau\bl{\nabla h}(x^{(i)})-\tau L^* u^{(i)} \big)\\
    u^{(i+\frac{1}{2})}
= \mathrm{prox}_{\sigma \dr{g^*}}\big(u^{(i)}+\sigma L \tilde{x}^{(i)} \big)\\
    x^{(i+\frac{1}{2})}
= \mathrm{prox}_{\tau \ff}\big(x^{(i)}-\tau\bl{\nabla h}(x^{(i)})-\tau L^* u^{(i+\frac{1}{2})} \big)\\
x^{(i+1)}=x^{(i)}+\rho^{(i)} (x^{(i+\frac{1}{2})}-x^{(i)})\\
u^{(i+1)}=u^{(i)}+\rho^{(i)} (u^{(i+\frac{1}{2})}-u^{(i)}).
\end{array}\right.}%
This algorithm has been analyzed by Chen et al.~\cite{che16} and called the Primal-Dual Fixed-Point (PDFP) algorithm, so we give it the same  name.

Since the algorithm has the structure of a primal-dual forward-backward algorithm, we can, again, invoke \cref{thm22}, and we obtain:

\begin{theorem}[PDFP algorithm \eqref{algorithm_1e}, affine $\mathrm{prox}_{\ff}$ case]\label{thm34}Suppose that $\mathrm{prox}_{\ff}$ is affine. Let $x^{(0)}\in\mathcal{X}$ and 
$u^{(0)}\in\mathcal{U}$. Let $\tau \in (0, 2/\beta)$ and $\sigma>0$ be such that $\sigma\tau\|L\dg{R_\tau}L^*\|<1$. Set $\delta=2-\tau\beta/2$. Let $(\rho^{(i)})_{i\in\mathbb{N}}$ be a sequence in $[0,\delta]$ such that $\sum_{i\in \mathbb{N}} \rho^{(i)} (\delta-\rho^{(i)})=+\infty$. 
Then the sequences $(x^{(i)})_{i\in \mathbb{N}}$ and $(u^{(i)})_{i\in \mathbb{N}}$ defined by the iteration \eqref{algorithm_1e} converge weakly to a solution of \eqref{eq13e} and a solution of \eqref{eq21e}, respectively.\end{theorem}

 We note that $\|L\dg{R_\tau}L^*\|\leq \|L\|^2$, so if $\sigma\tau\|L\|^2<1$, then $\sigma\tau\|L\dg{R_\tau}L^*\|<1$.

Chen et al. have proved convergence of the PDFP algorithm with any function $\ff$, not only in the affine proximity operator case, as follows~\cite[Theorem 3.1]{che16}:

\begin{theorem}[PDFP algorithm \eqref{algorithm_1e}, general case]\label{lem35}
Suppose that $\mathcal{X}$ and $\mathcal{U}$ are of finite dimension. Let $x^{(0)}\in\mathcal{X}$ and 
$u^{(0)}\in\mathcal{U}$. Let $\tau \in (0, 2/\beta)$ and $\sigma>0$ be such that $\sigma\tau\|L\|^2<1$. Set $\rho^{(i)}=1$ for every $i\in \mathbb{N}$. Then the sequences $(x^{(i)})_{i\in \mathbb{N}}$ and $(u^{(i)})_{i\in \mathbb{N}}$ defined by the iteration \eqref{algorithm_1e} converge  to a solution of \eqref{eq13e} and a solution of \eqref{eq21e}, respectively.\end{theorem}

It remains an open question whether the PDFP algorithm can be relaxed like in \cref{thm34} for an arbitrary $\ff$.

We can check that when $\h=\frac{\beta}{2}\|\cdot\|^2$, it is equivalent to applying the Loris--Verhoeven algorithm to minimize $\g\circ L +\bl{h}$ or to applying the PDFP algorithm to minimize $\ff+\g\circ L$, with $\ff=\frac{\beta}{2}\|\cdot\|^2$. Indeed, in the latter case,  let $\tau'>0$ and $\sigma>0$ be the parameters of the PDFP algorithm. Set $\tau =\frac{\tau'}{1+\tau'\beta}$.  We have $\mathrm{prox}_{\tau' \ff}:x\in\mathcal{X}\mapsto \frac{1}{1+\tau'\beta}x$, which is linear. Then the PDFP algorithm \eqref{algorithm_1e} becomes:
\begin{align}\label{algorithm_1f}
&\,\mbox{\textbf{PDFP iteration}: for }i=0,1,\ldots,\notag\\[-1mm]
    &\left\lfloor
    \begin{array}{l}
    u^{(i+\frac{1}{2})}
= \mathrm{prox}_{\sigma \dr{g^*}}\Big(u^{(i)}+\frac{\sigma}{1+\tau'\beta} L \big(x^{(i)}-\tau' L^* u^{(i)} \big) \Big)\\
\hphantom{u^{(i+\frac{1}{2})}}=\mathrm{prox}_{\sigma \dr{g^*}}\Big(u^{(i)}+\sigma L \big(x^{(i)}-\tau\bl{\nabla h}(x^{(i)})-\tau L^* u^{(i)} \big) \Big)\\
    x^{(i+\frac{1}{2})}
=\frac{1}{1+\tau'\beta}\big(x^{(i)}-\tau' L^* u^{(i+\frac{1}{2})} \big)\\
\hphantom{x^{(i+\frac{1}{2})}}=x^{(i)}-\tau\bl{\nabla h}(x^{(i)})-\tau L^* u^{(i+\frac{1}{2})}\\
x^{(i+1)}=x^{(i)}+\rho^{(i)} (x^{(i+\frac{1}{2})}-x^{(i)})\\
u^{(i+1)}=u^{(i)}+\rho^{(i)} (u^{(i+\frac{1}{2})}-u^{(i)}),
\end{array}\right.\end{align}
and we recover the Loris--Verhoeven algorithm  \eqref{algorithm_1} applied to minimize $\g\circ L +\bl{h}$, with parameters $\tau$ and $\sigma$. The conditions for convergence from the PDFP interpretation are $\tau<1/\beta$, $\sigma\tau\|L\|^2<1$, and $\delta=2$, so that we do not gain anything in comparison with \cref{thm33}.

\section{The Chambolle--Pock and the Douglas--Rachford Algorithms}\label{sec4}

We now consider a problem similar to \eqref{eq13}, but this time we want to make use of the proximity operators of the two functions.
So, let $\mathcal{X}$ and $\mathcal{U}$ be two real Hilbert spaces. Let $\ff\in \Gamma_0(\mathcal{X})$ and $\g\in \Gamma_0(\mathcal{U})$. Let $L:\mathcal{X}\rightarrow \mathcal{U}$ be a nonzero bounded linear operator.
We want to
\begin{equation}
    \minimize_{x\in\mathcal{X}}\, \dg{f(x)} + \dr{g(Lx)}.\label{eq22}
\end{equation}
The corresponding monotone inclusion is
\begin{equation}
    0 \in \dg{\partial f}(x) + L^*\dr{\partial g} (Lx).
\end{equation}
Again, to bypass the annoying operator $L$, we introduce an auxiliary variable $u\in\dr{\partial g} (Lx)$, which will be called the dual  variable, so that the problem now consists in finding $x\in\mathcal{X}$ and $u\in\mathcal{U}$ such that
\begin{equation}\label{variational_prim_dual_prob2}
  \left\{\begin{array}{l}
  u\in\dr{\partial g} (Lx)\\
  0 \in L^*u + \dg{\partial f}(x)
  \end{array}\right. .
\end{equation}
Let us define the dual convex optimization problem associated to the primal problem \eqref{eq22}:
\begin{equation}
    \minimize_{u\in\mathcal{U}}\, \dg{f^*(-L^*u)} + \dr{g^*(u)}.\label{eq23}
\end{equation}
If a pair $(x,u)\in\mathcal{X}\times\mathcal{U}$ is a solution to \eqref{variational_prim_dual_prob2}, then $x$ is a solution to \eqref{eq22} and $u$ is a solution to \eqref{eq23}.

To solve the primal and dual problems \eqref{eq22} and \eqref{eq23} jointly, Chambolle and Pock~\cite{cha11a} proposed the following algorithms (without relaxation); see also Esser et al.~\cite{ess10} and Zhang et al.~\cite{zha10}:
\ceq{align}{\label{algoCP}
&\,\mbox{\textbf{Chambolle--Pock iteration, form I, for \eqref{eq22} and \eqref{eq23}}: for }i=0,1,\ldots,\notag\\[-1mm]
    &\left\lfloor
    \begin{array}{l}
    x^{(i+\frac{1}{2})}
= \mathrm{prox}_{\tau \ff}\big(x^{(i)}-\tau L^* u^{(i)} \big)\\
u^{(i+\frac{1}{2})}
= \mathrm{prox}_{\sigma \dr{g^*}}\big(u^{(i)}+\sigma L(2x^{(i+\frac{1}{2})}-x^{(i)}) \big)\\
x^{(i+1)}=x^{(i)}+\rho^{(i)} (x^{(i+\frac{1}{2})}-x^{(i)})\\
u^{(i+1)}=u^{(i)}+\rho^{(i)} (u^{(i+\frac{1}{2})}-u^{(i)}).
\end{array}\right.}

\ceq{align}{\label{algoCP2}
&\,\mbox{\textbf{Chambolle--Pock iteration, form II, for \eqref{eq22} and \eqref{eq23}}:}\notag\\[-1mm]
&\,\mbox{for }i=0,1,\ldots,\notag\\[-1mm]
    &\left\lfloor
    \begin{array}{l}
u^{(i+\frac{1}{2})} = \mathrm{prox}_{\sigma \dr{g^*}}\big(u^{(i)}+\sigma L x^{(i)} \big)\\
    x^{(i+\frac{1}{2})}
= \mathrm{prox}_{\tau \ff}\big(x^{(i)}-\tau L^* (2u^{(i+\frac{1}{2})}-u^{(i)}) \big)\\
u^{(i+1)}=u^{(i)}+\rho^{(i)} (u^{(i+\frac{1}{2})}-u^{(i)})\\
x^{(i+1)}=x^{(i)}+\rho^{(i)} (x^{(i+\frac{1}{2})}-x^{(i)}).
\end{array}\right.}\smallskip

We note that the Chambolle--Pock algorithm is self-dual: if we apply the Chambolle--Pock iteration form I to the problem \eqref{eq23} to minimize $\dg{\tilde{f}}+\dr{\tilde{g}}\circ \tilde{L}$ with $\dg{\tilde{f}}=\dr{g^*}$, $\dr{\tilde{g}}=\dg{f^*}$, $\tilde{L}=-L^*$, and the roles of $x$ and $u$ are switched, as are the roles of $\tau$ and $\sigma$, we obtain exactly the Chambolle--Pock iteration form II.

Chambolle and Pock proved the convergence in the finite-dimensional case, assuming that 
$\tau\sigma\|L\|^2<1$ and $\rho^{(i)}= 1$~\cite{cha11a}. 
The convergence was proved in a different way by He and Yuan, with a constant relaxation parameter $\rho^{(i)}= \rho\in (0,2)$ and the other hypotheses the same~\cite{he10}; indeed, they observed that the Chambolle--Pock algorithm is a primal-dual proximal point algorithm to find a primal-dual pair  $z=(x,u)$   in 
$\mathcal{Z}=\mathcal{X}\times\mathcal{U}$, a solution to the monotone inclusion 
\begin{equation}\label{variational_prim_dual_probcp}
   \left( \begin{array}{c}
   0\\0
   \end{array}\right)\in
    \underbrace{\left( \begin{array}{c}
     \dg{\partial f}(x) + L^*u \\
   -Lx + (\dr{\partial g})^{-1}u
   \end{array}\right)}_{Mz}.
\end{equation}
The operator $M:\mathcal{Z}\rightarrow 2^\mathcal{Z},(x,u)\mapsto (\dg{\partial f}(x)+L^*u, -Lx + (\dr{\partial g})^{-1}u)$ is maximally monotone~\cite[Proposition 26.32 (iii)]{bau17}. 
Then one can observe that the Chambolle--Pock iteration form I satisfies
\begin{equation}\label{var_inclusioncp}
	\left( \begin{array}{c}
	0  \\ 
	0  \\
	\end{array} \right) 
	\in 
	\underbrace{\left( \begin{array}{c}
	\dg{\partial f} (x^{(i+\frac{1}{2})})+L^* u^{(i+\frac{1}{2})} \\ 
	-L x^{(i+\frac{1}{2})}+(\dr{\partial g})^{-1} u^{(i+\frac{1}{2})} \\
	\end{array} \right)}_{Mz^{(i+\frac{1}{2})}}
	+
	\underbrace{\left( \begin{array}{cc} 
	\frac{1}{\tau}\mathrm{Id} & -L^* \\ 
	-L & \frac{1}{\sigma}\mathrm{Id} \\
	\end{array} \right)}_{P}
	\underbrace{\left( \begin{array}{c}
	x^{(i+\frac{1}{2})}-x^{(i)}  \\ 
	u^{(i+\frac{1}{2})}-u^{(i)}\\
	\end{array} \right)}_{z^{(i+\frac{1}{2})}-z^{(i)}}.
\end{equation}
The Chambolle--Pock iteration form II satisfies the same primal-dual inclusion, but with $-L$ replaced by $L$ in $P$.

Thus, the Chambolle--Pock iteration is a preconditioned primal-dual proximal point algorithm and, as a consequence of  \cref{thm23}, we have~\cite[Theorem 3.2]{con13}:

\begin{theorem}[Chambolle--Pock algorithm \eqref{algoCP} or \eqref{algoCP2}]\label{thm41}
Let $x^{(0)}\in\mathcal{X}$ and 
$u^{(0)}\in\mathcal{U}$. Let $\tau >0$ and $\sigma>0$  be such that $\sigma\tau\|L\|^2<1$. Let $(\rho^{(i)})_{i\in\mathbb{N}}$ be a sequence in $[0,2]$ such that $\sum_{i\in \mathbb{N}} \rho^{(i)} (2-\rho^{(i)})=+\infty$. Then the sequences $(x^{(i)})_{i\in \mathbb{N}}$ and $(u^{(i)})_{i\in \mathbb{N}}$, defined either by the iteration \eqref{algoCP} or by the iteration \eqref{algoCP2}, converge weakly to a solution of \eqref{eq22} and a solution of \eqref{eq23}, respectively.\end{theorem}

In addition, the first author proved that in the finite-dimensional setting, one can set $\sigma\tau\|L\|^2=1$~\cite[Theorem 3.3]{con13}, see also O'Connor and Vandenberghe~\cite{oco18} for another proof and \cite{bri19} for a discussion in a generalized setting. If we apply \cref{thm52} below with $K=\mathrm{Id}$ and $\eta=1/\tau$, we obtain:

\begin{theorem}[Chambolle--Pock algorithm \eqref{algoCP} or \eqref{algoCP2}]\label{thm42}
Let $x^{(0)}\in\mathcal{X}$ and 
$u^{(0)}\in\mathcal{U}$. Let $\tau >0$ and $\sigma>0$  be such that $\sigma\tau\|L\|^2\leq 1$. Let $(\rho^{(i)})_{i\in\mathbb{N}}$ be a sequence in $[0,2]$ such that $\sum_{i\in \mathbb{N}} \rho^{(i)} (2-\rho^{(i)})=+\infty$.  Then the sequences $(x^{(i)})_{i\in \mathbb{N}}$ and $(u^{(i)})_{i\in \mathbb{N}}$, defined either by the iteration \eqref{algoCP} or by the iteration \eqref{algoCP2}, converge weakly to a solution of \eqref{eq22} and a solution of \eqref{eq23}, respectively.\end{theorem}

The difference between \cref{thm41} and \cref{thm42}
 is that $\sigma\tau\|L\|^2=1$ is allowed in the latter. This is a significant improvement: in practice, one can set $\sigma=1/(\tau\|L\|^2)$ in the algorithms and have only one parameter remaining, namely $\tau$, to tune.\medskip

\begin{figure}[t]
\centering
\begin{tabular}{ccc}
\includegraphics[scale=0.343]{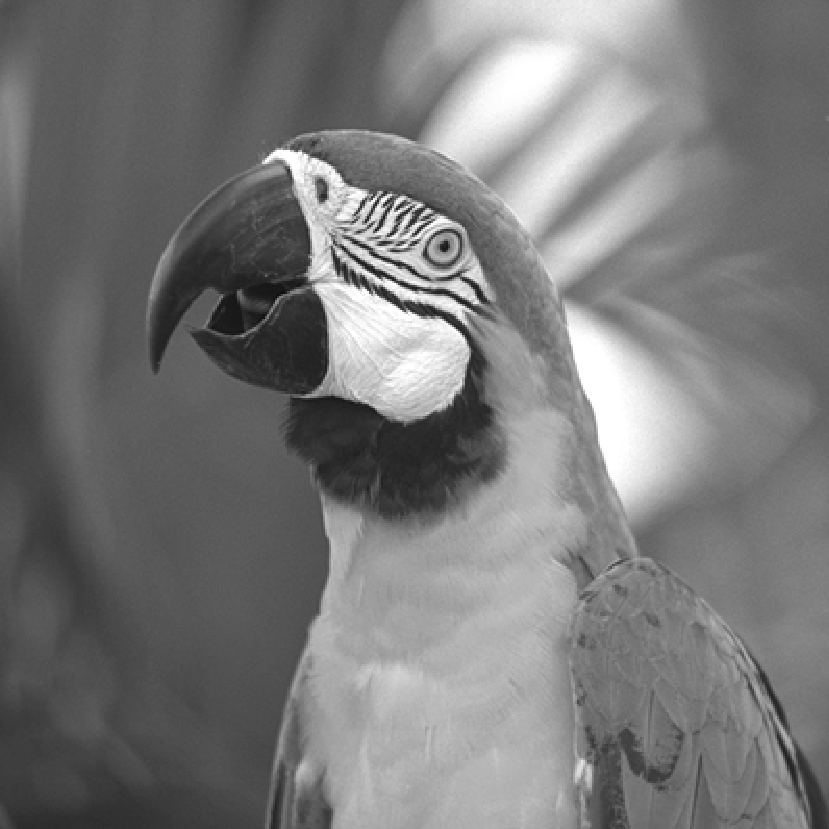}&\includegraphics[scale=0.343]{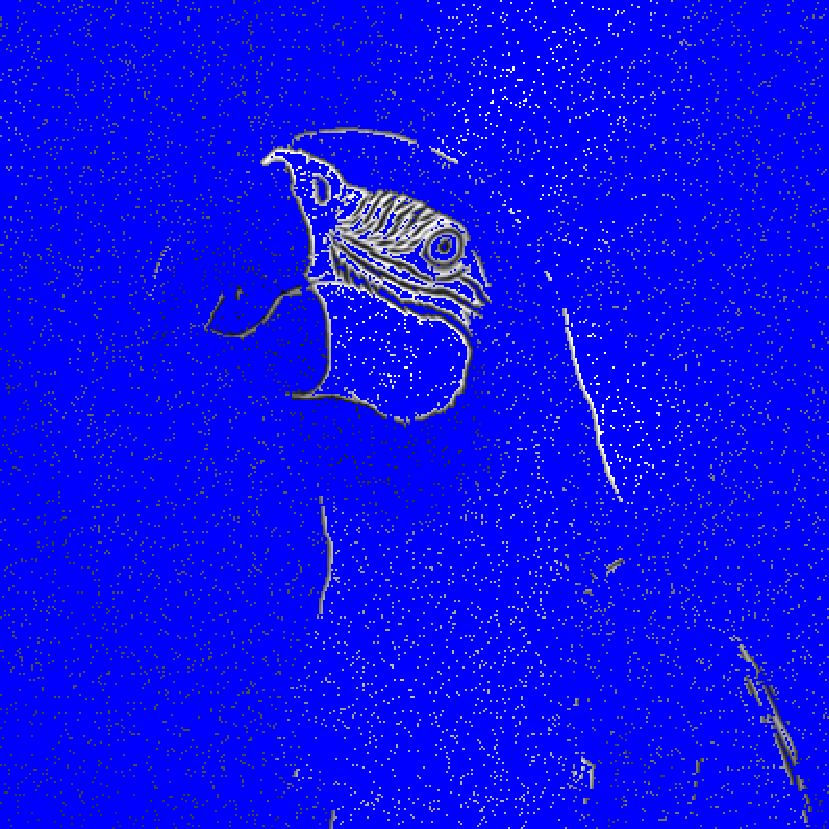}&\includegraphics[scale=0.343]{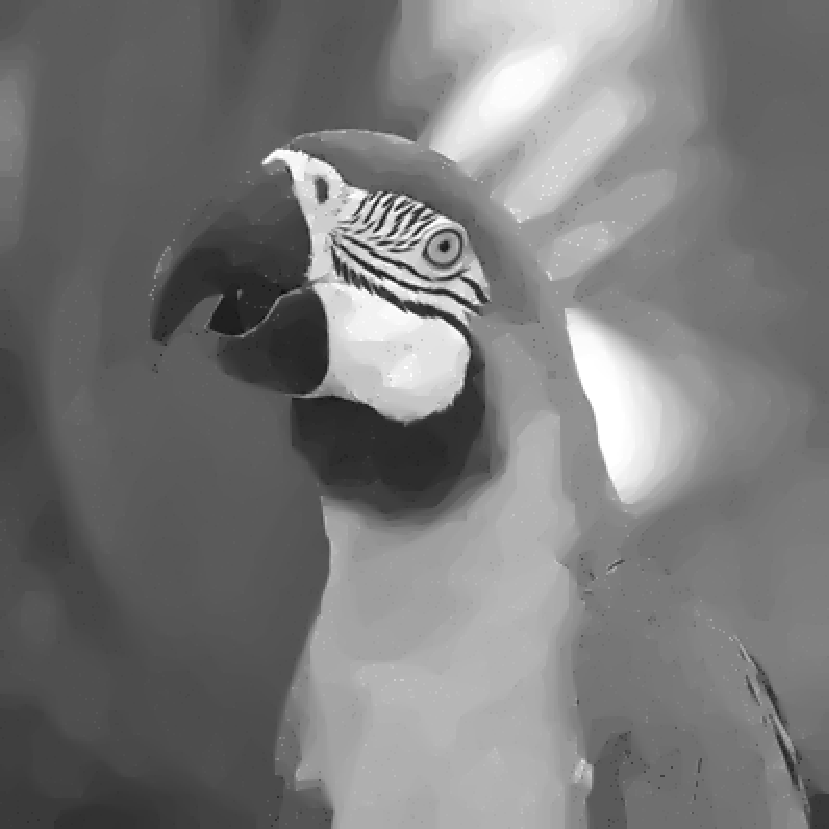}\\
(a)&(b)&(c)
\end{tabular}
\begin{tabular}{cc}
\hspace{-5mm}\includegraphics[scale=0.75]{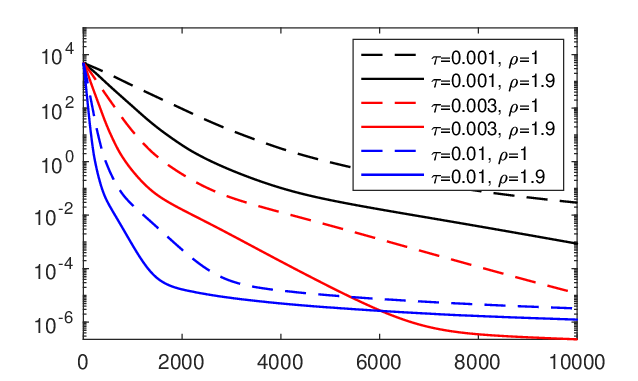}\hspace{-4mm}&\hspace{-4mm}\includegraphics[scale=0.75]{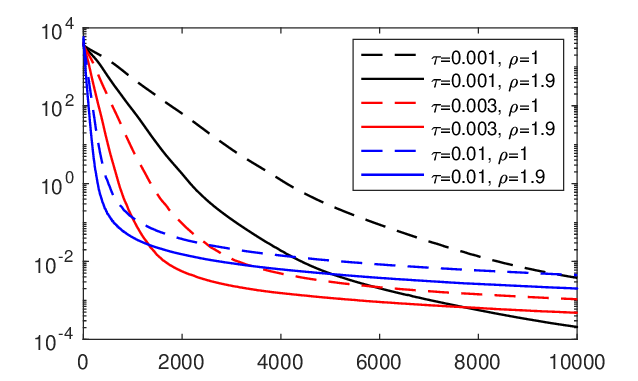}\hspace{-5mm}{\,}\\
\ (d) $\|x^{(i)}-x^\star\|^2$ w.r.t. $i\geq 1$&\ \ (e) $\Psi(x^{(i)})-\Psi(x^\star)$ w.r.t. $i\geq 1$
\end{tabular}
\caption{\br{Image inpainting experiment: the goal is to estimate the reference image (a) given only a subset (8\%) of its pixel values (b), with the missing values shown in blue. The reconstructed image (c) is the solution of the convex optimization problem of the form  \eqref{eq22} of minimizing the total variation, while keeping the available pixel values. The Chambolle--Pock algorithm \eqref{algoCP} is used with $\sigma=1/(8\tau)$, according to \cref{thm42}, since $\|L\|^2\leq 8$. The squared distance to the solution $x^\star$ and the decay of the objective function $\Psi$ along the iterations are shown in (d) and (e), respectively, for three different values of $\tau$ and $\rho^{(i)}$ constant and equal to $1$ (no relaxation) or $1.9$ (overrelaxation). The benefit of overrelaxation is clearly visible.}}
\label{figu1}
\end{figure}

\br{In \cref{figu1}, we illustrate how overrelaxation can accelerate the Chambolle--Pock algorithm. We consider the problem of image inpainting, which consists in estimating missing pixel values of an image. We formulate inpainting as an optimization problem of the form \eqref{eq22}; that is, we aim at minimizing  $\Psi = \dg{f} + \g\circ L$, where $\dg{f}$ is the indicator function of the affine subspace of images having the prescribed pixel values, and $\g\circ L$ is the ``isotropic'' total variation~\cite{rud92,cha16,con17}. 
So, we use the Chambolle--Pock algorithm to solve the problem; a solution $x^\star$ (which seems to be unique and is estimated by $x^{(i)}$ with $i=10^6$) is shown in \cref{figu1}~(c). We can see from the plots of $\|x^{(i)}-x^\star\|^2$ and $\Psi(x^{(i)})-\Psi(x^\star)=\g(Lx^{(i)})-\g(Lx^\star)$ (since  $\ff(x^{(i)})=0$ for every $i$) that the best value of $\tau$ depends on the desired accuracy or number of iterations: $\tau=0.01$ is better after $10^3$ iterations, whereas a smaller value is better after $10^4$ iterations. In any case, $\rho^{(i)}=1.9$ significantly speeds up the convergence, in comparison with $\rho^{(i)}=1$.}

\subsection{The Proximal Method of Multipliers}\label{secpmm}

Let $\mathcal{X}$ and $\mathcal{U}$ be two real Hilbert spaces. Let $\g\in \Gamma_0(\mathcal{U})$, $c\in\mathcal{X}$, and let $L:\mathcal{X}\rightarrow \mathcal{U}$ be a bounded linear operator. We consider the problem
\begin{equation}
    \minimize_{x\in\mathcal{X}}\, \dr{g(Lx)} + \langle x,c\rangle\label{eqpalm1}
\end{equation}
and the dual problem
\begin{equation}
    \minimize_{u\in\mathcal{U}}\, \dr{g^*(u)}\quad\mbox{s.t.}\quad L^*u+c=0.\label{eqpalm2}
\end{equation}
We can view the term $\langle x,c\rangle$ as a smooth function $\h$ with constant gradient $\bl{\nabla h}=c$ and apply the Loris--Verhoeven algorithm. Thus, let $\tau>0$ and $\sigma>0$, let $x^{(0)}\in\mathcal{X}$ and 
$u^{(0)}\in\mathcal{U}$, and let $(\rho^{(i)})_{i\in\mathbb{N}}$ be a sequence of relaxation parameters. The iteration is (we change the name of the variable $x$ to $s$):
\begin{align}
\label{algorithm_pmm1}
&\,\mbox{\textbf{Loris--Verhoeven iteration for \eqref{eqpalm1} and \eqref{eqpalm2}}: for }i=0,1,\ldots,\notag\\[-1mm]
    &\left\lfloor
    \begin{array}{l}
    u^{(i+\frac{1}{2})}
= \mathrm{prox}_{\sigma \dr{g^*}}\big(u^{(i)}+\sigma L (s^{(i)} -\tau L^* u^{(i)}-\tau c)\big)\\
s^{(i+1)}=s^{(i)}-\rho^{(i)}\tau (
L^* u^{(i+\frac{1}{2})}+c)\\
u^{(i+1)}=u^{(i)}+\rho^{(i)} (u^{(i+\frac{1}{2})}-u^{(i)}).
\end{array}\right.
\end{align}
Each iteration of the algorithm satisfies the inclusion
\begin{equation}\label{var_inclusionpmm}
	\left(\!\begin{array}{c}
	0  \\ 
	0  \\
	\end{array}\!\right) 
	\in 
	\left( \begin{array}{c}
	c+L^* u^{(i+\frac{1}{2})} \\ 
	-L s^{(i+\frac{1}{2})}+(\dr{\partial g})^{-1} u^{(i+\frac{1}{2})} \\
	\end{array} \right)
	+
	\left( \begin{array}{cc} 
	\frac{1}{\tau}\mathrm{Id} & 0 \\ 
	0 & \frac{1}{\sigma}\mathrm{Id}-\tau LL^* \\
	\end{array} \right)
	\left( \begin{array}{c}
	s^{(i+\frac{1}{2})}-s^{(i)}  \\ 
	u^{(i+\frac{1}{2})}-u^{(i)}\\
	\end{array} \right)\!.
\end{equation}
Since the algorithm has the structure of a proximal point algorithm, we can apply  \cref{thm23} and obtain convergence under the condition $\sigma\tau\|L\|^2<1$.

An alternative  is to see the term $\langle x,c\rangle$ as a function $\ff$ with proximity operator $\mathrm{prox}_{\tau \ff}:x\mapsto x-\tau c$. Thus, we can solve \eqref{eqpalm1} and \eqref{eqpalm2} by applying the Chambolle--Pock algorithm (we change the name of its variable $x^{(i+\frac{1}{2})}$  to  $s^{(i)}$):
\begin{align}\label{algopmm2}
&\,\mbox{\textbf{Chambolle--Pock iteration, form I, for \eqref{eqpalm1} and \eqref{eqpalm2}}:}\notag\\[-1mm]
&\,\mbox{for }i=0,1,\ldots,\notag\\[-1mm]
    &\left\lfloor
    \begin{array}{l}
    s^{(i)}
= x^{(i)}-\tau L^* u^{(i)} -\tau c\\
u^{(i+\frac{1}{2})}
= \mathrm{prox}_{\sigma \dr{g^*}}\big(u^{(i)}+\sigma L(2s^{(i)}-x^{(i)}) \big)\\
\hphantom{u^{(i+\frac{1}{2})}}{}=\mathrm{prox}_{\sigma \dr{g^*}}\big(u^{(i)}+\sigma L (s^{(i)}-\tau L^* u^{(i)} -\tau c) \big)\\
x^{(i+1)}=x^{(i)}+\rho^{(i)} (s^{(i)}-x^{(i)})\\
\hphantom{x^{(i+1)}}{}=x^{(i)}-\rho^{(i)} \tau(L^* u^{(i)} + c)\\
u^{(i+1)}=u^{(i)}+\rho^{(i)} (u^{(i+\frac{1}{2})}-u^{(i)}).
\end{array}\right.\end{align}
We can discard the variable $x$ and the iteration becomes:
\begin{align}\label{algopmm22}
    &\left\lfloor
    \begin{array}{l}
  u^{(i+\frac{1}{2})}
=\mathrm{prox}_{\sigma \dr{g^*}}\big(u^{(i)}+\sigma L (s^{(i)}-\tau L^* u^{(i)} -\tau c) \big)\\
s^{(i+1)}=s^{(i)}-\rho^{(i)}\tau \big(
L^* u^{(i+\frac{1}{2})}+c\big)\\
u^{(i+1)}=u^{(i)}+\rho^{(i)} (u^{(i+\frac{1}{2})}-u^{(i)}),
\end{array}\right.\end{align}
and we exactly recover the iteration \eqref{algorithm_pmm1}. 
Therefore, the Loris--Verhoeven algorithm and the Chambolle--Pock algorithm are the same for the specific problems \eqref{eqpalm1} and \eqref{eqpalm2}; we call this algorithm the Proximal Method of Multipliers~\cite{roc762}.

Furthermore, the form \eqref{algopmm2} provides us with a compact form  of the algorithm:
\ceq{align}{\label{algorithm_pmm12}
&\,\mbox{\textbf{Proximal Method of Multipliers for \eqref{eqpalm1} and \eqref{eqpalm2}}: for }i=0,1,\ldots,\notag\\[-1mm]
    &\left\lfloor
    \begin{array}{l}
    a^{(i)}=L^* u^{(i)}+c\\
    u^{(i+1)}=u^{(i)}+\rho^{(i)} \Big(\mathrm{prox}_{\sigma \dr{g^*}}\big(u^{(i)}+\sigma L (x^{(i)}-2\tau a^{(i)} ) \big)-u^{(i)}\Big)\\
x^{(i+1)}=x^{(i)}-\rho^{(i)} \tau a^{(i)},
\end{array}\right.}
and since the two updates are independent, they can be performed in parallel.

\cref{thm42} yields:

\begin{theorem}[Proximal Method of Multipliers \eqref{algorithm_pmm12}]\label{thm43}
Let $x^{(0)}\in\mathcal{X}$ and 
$u^{(0)}\in\mathcal{U}$. Let $\tau >0$ and $\sigma>0$ be such that $\sigma\tau\|L\|^2\leq 1$. Let $(\rho^{(i)})_{i\in\mathbb{N}}$ be a sequence in $[0,2]$ such that $\sum_{i\in \mathbb{N}} \rho^{(i)} (2-\rho^{(i)})=+\infty$. 
Then the sequences $(x^{(i)})_{i\in \mathbb{N}}$ and $(u^{(i)})_{i\in \mathbb{N}}$ defined by the iteration \eqref{algorithm_pmm12} converge weakly to a solution of \eqref{eqpalm1} and a solution of \eqref{eqpalm2}, respectively.\end{theorem}

We can also apply the form II of the Chambolle--Pock algorithm, with the same convergence properties as in \cref{thm43}:
\ceq{align}{\label{algopmm3}
&\,\mbox{\textbf{Chambolle--Pock iteration, form II, for \eqref{eqpalm1} and \eqref{eqpalm2}}:}\notag\\[-1mm]
&\,\mbox{for }i=0,1,\ldots,\notag\\[-1mm]
    &\left\lfloor
    \begin{array}{l}
u^{(i+\frac{1}{2})} = \mathrm{prox}_{\sigma \dr{g^*}}\big(u^{(i)}+\sigma L x^{(i)} \big)\\
    x^{(i+1)}=x^{(i)}-\rho^{(i)} \tau\big(L^* (2u^{(i+\frac{1}{2})}-u^{(i)})+c\big)\\
u^{(i+1)}=u^{(i)}+\rho^{(i)} (u^{(i+\frac{1}{2})}-u^{(i)}).
\end{array}\right.}

\subsection{The Douglas--Rachford algorithm}\label{secdr}

Let us consider the particular case of the problem \eqref{eq22} when 
$\mathcal{X}=\mathcal{U}$ and $L=\mathrm{Id}$. The problem becomes
\begin{equation}
    \minimize_{x\in\mathcal{X}}\, \dg{f(x)} + \dr{g(x)}.\label{eqpdr}
\end{equation}
By setting $\sigma=1/\tau$ in the Chambolle--Pock algorithm, we recover the Douglas--Rachford algorithm \cite{lio79,eck92,com04,sva11} as a particular case:
\begin{align}\label{algorithm_dr}
&\,\mbox{\textbf{Douglas--Rachford iteration for \eqref{eqpdr}}: for }i=0,1,\ldots,\notag\\[-1mm]
    &\left\lfloor
    \begin{array}{l}
   x^{(i+\frac{1}{2})}=\mathrm{prox}_{\tau \ff }\big(x^{(i)}-\tau  u^{(i)}
\big)\\
    u^{(i+\frac{1}{2})}
= \mathrm{prox}_{ \dr{g^*}/\tau}\big(u^{(i)}+ (2x^{(i+\frac{1}{2})} - x^{(i)})/\tau\big)\\
x^{(i+1)}=x^{(i)}+\rho^{(i)} (x^{(i+\frac{1}{2})}-x^{(i)})\\
u^{(i+1)}=u^{(i)}+\rho^{(i)} (u^{(i+\frac{1}{2})}-u^{(i)}).
\end{array}\right.
\end{align}

Note that the Douglas--Rachford algorithm is not a primal-dual proximal point algorithm since the operator $P$ in \eqref{var_inclusioncp} is no longer strongly positive. So, its weak convergence, which is not strong in general~\cite{bui20}, is difficult to establish; it was only shown in 2011 by Svaiter~\cite{sva11}. 
Let us have a quick look at the structure of the Douglas--Rachford algorithm. First, we define the auxiliary variable 
\begin{equation}
s^{(i)}=x^{(i)}-\tau u^{(i)}.
\end{equation}
The Douglas--Rachford iteration depends only on this concatenated variable, not on the full pair $(x^{(i)},u^{(i)})$, and we can rewrite it as:
\begin{align}\label{algorithm_2}
&\,\mbox{\textbf{Douglas--Rachford iteration for \eqref{eqpdr}}: for }i=0,1,\ldots,\notag\\[-1mm]
    &\left\lfloor
    \begin{array}{l}
   x^{(i+\frac{1}{2})}=\mathrm{prox}_{\tau \ff }(s^{(i)})\\
    u^{(i+\frac{1}{2})}
= \mathrm{prox}_{\dr{g^*}/\tau}\big((2x^{(i+\frac{1}{2})} - s^{(i)})/\tau\big)\\
s^{(i+1)}=s^{(i)}+\rho^{(i)} (x^{(i+\frac{1}{2})}-\tau u^{(i+\frac{1}{2})}-s^{(i)}),
\end{array}\right.
\end{align}
or, equivalently,
\ceq{align}{\label{algorithm_22}
&\,\mbox{\textbf{Douglas--Rachford iteration for \eqref{eqpdr}}: for }i=0,1,\ldots,\notag\\[-1mm]
    &\left\lfloor
    \begin{array}{l}
   x^{(i+\frac{1}{2})}=\mathrm{prox}_{\tau \ff }(s^{(i)})\\
s^{(i+1)}=s^{(i)}+\rho^{(i)}\big( \mathrm{prox}_{ \tau \g}(2x^{(i+\frac{1}{2})} - s^{(i)})-x^{(i+\frac{1}{2})}\big).
\end{array}\right.}
Or, keeping only the variable $s$:
\begin{equation}
\textstyle\left\lfloor
    \begin{array}{l}
    s^{(i+\frac{1}{2})}=\big(\frac{1}{2}(2\mathrm{prox}_{\tau \g }-\mathrm{Id}) (2\mathrm{prox}_{\tau \ff }-\mathrm{Id})+\frac{1}{2}\mathrm{Id}\big)s^{(i)}\\
s^{(i+1)}=s^{(i)}+\rho^{(i)}(s^{(i+\frac{1}{2})}-s^{(i)}).
\end{array}\right.
\end{equation}
Since the operator mapping $s^{(i)}$ to $s^{(i+\frac{1}{2})}$ is firmly nonexpansive, $(s^{(i)})_{i\in \mathbb{N}}$ converges weakly, after the Krasnosel'ski\u{\i}--Mann theorem, for every $\tau>0$ and sequence $(\rho^{(i)})_{i\in\mathbb{N}}$ in $[0,2]$ such that $\sum_{i\in \mathbb{N}} \rho^{(i)} (2-\rho^{(i)})=+\infty$. But in infinite dimension, 
one cannot deduce 
that $(x^{(i+\frac{1}{2})})_{i\in \mathbb{N}}$ converges. In finite dimension, however, the proximity operator is continuous, so that one can indeed conclude that $(x^{(i+\frac{1}{2})})_{i\in \mathbb{N}}$ converges to some minimizer of $\ff+\g$.

As an application of Corollary 28.3 in \cite{bau17}, which proves weak convergence using fine properties of the operators, we have:

\begin{theorem}[Douglas--Rachford algorithm \eqref{algorithm_22}]\label{thm44}Let $s^{(0)}\in\mathcal{X}$. Let $\tau >0$ and let $(\rho^{(i)})_{i\in\mathbb{N}}$ be a sequence in $[0,2]$ such that $\sum_{i\in \mathbb{N}} \rho^{(i)} (2-\rho^{(i)})=+\infty$.  Then the sequence $(x^{(i+\frac{1}{2})})_{i\in \mathbb{N}}$ defined by the iteration \eqref{algorithm_22} converges weakly to a solution of \eqref{eqpdr}.\end{theorem}

Note that it would be misleading to state \cref{thm44} by invoking \cref{thm42}, since \cref{thm42} is itself based on the convergence result of the Douglas--Rachford algorithm~\cite[Corollary 28.3]{bau17}. In other words, weak convergence of the Douglas--Rachford algorithm is the fundamental result from which we deduce weak convergence of the Chambolle--Pock algorithm, which generalizes it.\medskip

\br{As mentioned in \cite{com10}, the limit case of  the Douglas--Rachford algorithm with  $\rho^{(i)}=2$ is the Peaceman--Rachford algorithm~\cite{lio79,eck92,com04}. But  in that case, the operator $T=(2\mathrm{prox}_{\tau \g }-\mathrm{Id}) (2\mathrm{prox}_{\tau \ff }-\mathrm{Id})
$ mapping $s^{(i)}$ to $s^{(i+1)}$ is merely nonexpansive and not averaged. For instance, with  $\ff=0$ and $\g:x\mapsto(0\ $if$\ x=0,\ +\infty\ $otherwise$)$, so that the minimizer of $\ff+\g$ is $0$, we have $T=-\mathrm{Id}$ and the iteration $x^{(i+\frac{1}{2})} =s^{(i)}$, $s^{(i+1)}=-s^{(i)}$, which cycles and does not converge if $s^{(0)}\neq 0$. This shows the tightness of the Krasnosel'ski\u{\i}--Mann theorem. The Peaceman--Rachford algorithm converges under additional assumptions,  for instance if $\ff$ is strictly convex and real-valued \cite{com09}.
}\medskip

Let us show that, like the Chambolle--Pock algorithm, the Douglas--Rachford algorithm is self-dual \cite[Lemma~3.6]{eck89}\cite{yan16}: it is equivalent to apply it to the primal problem \eqref{eqpdr} or to the dual problem: 
\begin{equation}
    \minimize_{u\in\mathcal{X}}\, \dg{f^*(-u)} + \dr{g^*(u)}.\label{eqdrd}
\end{equation}
Using the Moreau identity \eqref{eqmoreau} and starting from the form \eqref{algorithm_dr}, we can write the Douglas--Rachford iteration  as:
\begin{align}\label{algorithm_280}
    &\left\lfloor
    \begin{array}{l}
    v^{(i+\frac{1}{2})}=\mathrm{prox}_{\dg{f^*(-\cdot)}/\tau }\big(u^{(i)} -\frac{1}{\tau} x^{(i)} \big)\\
   x^{(i+\frac{1}{2})}=x^{(i)}- \tau (u^{(i)} - v^{(i+\frac{1}{2})})\\
    u^{(i+\frac{1}{2})}
= \mathrm{prox}_{ \dr{g^*}/\tau}\big(  v^{(i+\frac{1}{2})} +\frac{1}{\tau}x^{(i+\frac{1}{2})} \big)\\
x^{(i+1)}=x^{(i)}+\rho^{(i)} (x^{(i+\frac{1}{2})}-x^{(i)})\\
u^{(i+1)}=u^{(i)}+\rho^{(i)} (u^{(i+\frac{1}{2})}-u^{(i)}).
\end{array}\right.\end{align}
If we introduce the variable $\tilde{s}^{(i)}=\frac{-1}{\tau}s^{(i)}=u^{(i)}-\frac{1}{\tau} x^{(i)} $, we can rewrite the iteration as:
\begin{align}\label{algorithm_281}
    &\left\lfloor
    \begin{array}{l}
    v^{(i+\frac{1}{2})}=\mathrm{prox}_{\dg{f^*(-\cdot)}/\tau }(\tilde{s}^{(i)} )\\
   x^{(i+\frac{1}{2})}=x^{(i)}- \tau (u^{(i)} - v^{(i+\frac{1}{2})})\\
    u^{(i+\frac{1}{2})}
= \mathrm{prox}_{ \dr{g^*}/\tau}\big(2v^{(i+\frac{1}{2})} -\tilde{s}^{(i)} \big)\\
x^{(i+1)}=x^{(i)}+\rho^{(i)} (x^{(i+\frac{1}{2})}-x^{(i)})\\
u^{(i+1)}=u^{(i)}+\rho^{(i)} (u^{(i+\frac{1}{2})}-u^{(i)})\\
\tilde{s}^{(i+1)}=u^{(i+1)}-\frac{1}{\tau} x^{(i+1)} \\
\hphantom{\tilde{s}^{(i+1)}}=\tilde{s}^{(i)}+\rho^{(i)} \big(\mathrm{prox}_{ \dr{g^*}/\tau}(2v^{(i+\frac{1}{2})} -\tilde{s}^{(i)}) -v^{(i+\frac{1}{2})}\big),
\end{array}\right.\end{align}
which can be simplified as:
\begin{align}\label{algorithm_282}
&\,\mbox{\textbf{Douglas--Rachford iteration for \eqref{eqdrd}}: for }i=0,1,\ldots,\notag\\[-1mm]
    &\left\lfloor
    \begin{array}{l}
    v^{(i+\frac{1}{2})}=\mathrm{prox}_{\dg{f^*(-\cdot)}/\tau }(\tilde{s}^{(i)} )\\
    \tilde{s}^{(i+1)}=\tilde{s}^{(i)}+\rho^{(i)} \big(\mathrm{prox}_{ \dr{g^*}/\tau}(2v^{(i+\frac{1}{2})} -\tilde{s}^{(i)}) -v^{(i+\frac{1}{2})}\big).
\end{array}\right.\end{align}
Thus, we recognize the form \eqref{algorithm_22} of the Douglas--Rachford algorithm applied to the dual problem \eqref{eqdrd}, with parameter $1/\tau$. All in all, there are only two different instances of the Douglas--Rachford algorithm: the one given here, and the one obtained by switching the roles of $\ff$ and $\g$ (which we could obtain as a particular case of the Chambolle--Pock algorithm, form II).

Douglas--Rachford splitting has been shown to be the only way, in some sense, to minimize a sum of two functions, by calling their proximity operators~\cite{ryu20b}.

\subsubsection{A slightly modified version of the Douglas--Rachford algorithm}

Let us consider a slightly modified version of the Douglas--Rachford algorithm, which we will use in the next section. 
In addition to $\ff$ and $\g$ as before, let $c\in\mathcal{X}$. 
We consider the problem:
\begin{equation}
    \minimize_{x\in\mathcal{X}}\, \dg{f(x)} + \dr{g(x)} + \langle x,c\rangle.\label{eq27}
\end{equation}
The dual problem is 
\begin{equation}
    \minimize_{u\in\mathcal{X}}\, \dg{f^*(-u-c)} + \dr{g^*(u)}.\label{eq28}
\end{equation}
The linear term $\langle x,c\rangle$ does not add any difficulty: one can apply the Douglas--Rachford algorithm  to minimize $\dg{\tilde{f}}+\g$, where $\dg{\tilde{f}}=\ff+\langle \cdot,c\rangle$, using the fact that $\mathrm{prox}_{\tau \dg{\tilde{f}}}(x)=\mathrm{prox}_{\tau \ff}(x-\tau c)$. So, the modified Douglas--Rachford algorithm is:
\begin{align}\label{algodrc}
&\,\mbox{\textbf{Douglas--Rachford iteration for \eqref{eq27}}: for }i=0,1,\ldots,\notag\\[-1mm]
    &\left\lfloor
    \begin{array}{l}
   x^{(i+\frac{1}{2})}=\mathrm{prox}_{\tau \ff }(s^{(i)}-\tau c)\\
s^{(i+1)}=s^{(i)}+\rho^{(i)}\big( \mathrm{prox}_{ \tau \g}(2x^{(i+\frac{1}{2})} - s^{(i)})-x^{(i+\frac{1}{2})}\big).
\end{array}\right.\end{align}
Using the Moreau identity \eqref{eqmoreau} and starting from the form \eqref{algorithm_dr}, we can write the algorithm as
\ceq{align}{\label{algorithm_28}
&\,\mbox{\textbf{Douglas--Rachford iteration for \eqref{eq27}}: for }i=0,1,\ldots,\notag\\[-1mm]
    &\left\lfloor
    \begin{array}{l}
    v^{(i+\frac{1}{2})}=\mathrm{prox}_{\dg{f^*}/\tau }\big(\frac{1}{\tau} x^{(i)} - u^{(i)} - c\big)\\
   x^{(i+\frac{1}{2})}=x^{(i)}- \tau (u^{(i)} + v^{(i+\frac{1}{2})}+c)\\
    u^{(i+\frac{1}{2})}
= \mathrm{prox}_{ \dr{g^*}/\tau}\big(\frac{1}{\tau} x^{(i+\frac{1}{2})} - v^{(i+\frac{1}{2})} - c\big)\\
x^{(i+1)}=x^{(i)}+\rho^{(i)} (x^{(i+\frac{1}{2})}-x^{(i)})\\
u^{(i+1)}=u^{(i)}+\rho^{(i)} (u^{(i+\frac{1}{2})}-u^{(i)})
\end{array}\right.}
(note that the variable $v$ here is not the same as the one in \eqref{algorithm_282} or \eqref{algorithm_admm}; it is its opposite).

As an application of Corollary 28.3 in \cite{bau17}, we have:

\begin{theorem}[Douglas--Rachford algorithm \eqref{algorithm_28}]\label{thm45}
Let $x^{(0)}\in\mathcal{X}$ and 
$u^{(0)}\in\mathcal{X}$. Let $\tau >0$. Let $(\rho^{(i)})_{i\in\mathbb{N}}$ be a sequence in $[0,2]$ such that $\sum_{i\in \mathbb{N}} \rho^{(i)} (2-\rho^{(i)})=+\infty$.  Then the sequences $(x^{(i)})_{i\in \mathbb{N}}$ and $(u^{(i)})_{i\in \mathbb{N}}$ defined by the iteration \eqref{algorithm_28} converge weakly to a solution of \eqref{eq27} and a solution of \eqref{eq28}, respectively. Moreover, 
$(u^{(i+\frac{1}{2})}+v^{(i+\frac{1}{2})}+c)_{i\in \mathbb{N}}$ converges strongly to $0$.
\end{theorem}

\subsection{The Alternating Direction Method of Multipliers (ADMM)}

The Alternating Direction Method of Multipliers (ADMM) goes back to Glowinski and Marocco \cite{glo75}, and Gabay and Mercier \cite{gab76}. This algorithm has been studied extensively; see e.g.~\cite{eck89,eck92,cai10,boy11,eck12,eck15,glo14,glo162}. 
The ADMM was rediscovered by Osher et al.\ and called the \emph{Split Bregman Algorithm}~\cite{gol09,set11}; 
this method has received significant attention in image processing~\cite{che132,wu10, gol10, cai09,cai10}. 
The equivalence between the ADMM and the Split Bregman Algorithm is now well established~\cite{set09,ess10,set11}.
Also, the ADMM has been popularized in image processing by a series of papers of Figueiredo, Bioucas-Dias and others;  e.g.~\cite{fig10,afo10,afo11}. 

It is well known that the ADMM is equivalent to the Douglas--Rachford algorithm~\cite{glo89,gab83,eck92,yan16}, but the possibility of relaxation is often ignored. So, let us show this equivalence again: starting from the form \eqref{algorithm_dr} of the Douglas--Rachford algorithm and using the Moreau identity \eqref{eqmoreau}, we can rewrite it  as
\begin{align}\label{algorithm_admm}
    &\left\lfloor
    \begin{array}{l}
   x^{(i+\frac{1}{2})}=\mathrm{prox}_{\tau \ff }\big(x^{(i)}-\tau  u^{(i)}\big)\\
 v^{(i+\frac{1}{2})}=u^{(i)}+\frac{1}{\tau}(x^{(i+\frac{1}{2})}-x^{(i)})\\
w^{(i+1)} = \mathrm{prox}_{ \tau\dr{g}}\big(\tau v^{(i+\frac{1}{2})} +x^{(i+\frac{1}{2})}\big)\\
    u^{(i+\frac{1}{2})}= v^{(i+\frac{1}{2})}+\frac{1}{\tau}(x^{(i+\frac{1}{2})}-w^{(i+1)})\\
x^{(i+1)}=x^{(i)}+\rho^{(i)} (x^{(i+\frac{1}{2})}-x^{(i)})\\
u^{(i+1)}=u^{(i)}+\rho^{(i)} (u^{(i+\frac{1}{2})}-u^{(i)}).
\end{array}\right.
\end{align}
Let us introduce the variable $v^{(i+1)}=v^{(i+\frac{1}{2})}+\frac{\rho^{(i)}-1}{\tau}(x^{(i+\frac{1}{2})}-w^{(i+1)})$ for every $i\in\mathbb{N}$. Let  $w^{(0)}\in\mathcal{X}$ and  $v^{(0)}\in\mathcal{X}$ be such that $x^{(0)}-\tau  u^{(0)}=w^{(0)}-\tau  v^{(0)}$. Then for every $i\in\mathbb{N}$, $x^{(i+1)}-\tau  u^{(i+1)}=(1-\rho^{(i)})(x^{(i)}-\tau  u^{(i)})+\rho^{(i)}(x^{(i+\frac{1}{2})}-\tau u^{(i+\frac{1}{2})})=
(1-\rho^{(i)})(x^{(i+\frac{1}{2})}-\tau  v^{(i+\frac{1}{2})})+\rho^{(i)}(w^{(i+1)}-\tau v^{(i+\frac{1}{2})})=w^{(i+1)}-\tau v^{(i+1)}$. Moreover, we introduce the scaled variables $\tilde{v}^{(i)}=\tau v^{(i)}$ and $\tilde{v}^{(i+\frac{1}{2})}=\tau v^{(i+\frac{1}{2})}$.
Hence, we can remove the variables $x$ and $u$ and rewrite the iteration as
\ceq{align}{\label{algorithm_admm2}
&\,\mbox{\textbf{ADMM iteration for \eqref{eqpdr} and \eqref{eqdrd}}: for }i=0,1,\ldots,\notag\\[-1mm]
     &\left\lfloor
    \begin{array}{l}
   x^{(i+\frac{1}{2})}=\mathrm{prox}_{\tau \ff }\big(w^{(i)}-\tilde{v}^{(i)}\big)\\
 \tilde{v}^{(i+\frac{1}{2})}=\tilde{v}^{(i)}+x^{(i+\frac{1}{2})}-w^{(i)}\\
w^{(i+1)} = \mathrm{prox}_{ \tau\dr{g}}\big(x^{(i+\frac{1}{2})}+\tilde{v}^{(i+\frac{1}{2})}\big)\\
  \tilde{v}^{(i+1)}=\tilde{v}^{(i+\frac{1}{2})}+(\rho^{(i)}-1)(x^{(i+\frac{1}{2})}-w^{(i+1)}).
\end{array}\right.}
As an application of Corollary 28.3 in \cite{bau17}, we have:

\begin{theorem}[ADMM \eqref{algorithm_admm2}]\label{thm46}
Let $w^{(0)}\in\mathcal{X}$ and $\tilde{v}^{(0)}\in\mathcal{X}$. Let $\tau >0$ and let $(\rho^{(i)})_{i\in\mathbb{N}}$ be a sequence in $[0,2]$ such that $\sum_{i\in \mathbb{N}} \rho^{(i)} (2-\rho^{(i)})=+\infty$.  Then the sequences $(x^{(i+\frac{1}{2})})_{i\in \mathbb{N}}$ and $(w^{(i)})_{i\in \mathbb{N}}$ defined by the iteration \eqref{algorithm_admm2} both converge weakly to some element $x^\star\in\mathcal{X}$ solution to \eqref{eqpdr}. Moreover, the sequences $(\tilde{v}^{(i+\frac{1}{2})}/\tau)_{i\in \mathbb{N}}$ and $(\tilde{v}^{(i)}/\tau)_{i\in \mathbb{N}}$  defined by the iteration \eqref{algorithm_admm2} both converge weakly to some element $u^\star\in\mathcal{X}$ solution to \eqref{eqdrd}.  In addition, $(x^{(i+\frac{1}{2})}-w^{(i+1)})_{i\in \mathbb{N}}$ converges strongly to $0$.\end{theorem}

Thus, in the Douglas--Rachford algorithm, or equivalently the ADMM, there are two primal variables $x$ and $w$ and two dual variables $u$ and $v$, as can be seen in \eqref{algorithm_admm}. However, depending on how the algorithm is written, only some of these variables appear. Also, two monotone inclusions are satisfied in parallel at every iteration:
\begin{equation}\label{var_inclusioncpa1}
	\left( \begin{array}{c}
	0  \\ 
	0  \\
	\end{array} \right) 
	\in 
	{\left( \begin{array}{c}
	\dg{\partial f} (x^{(i+\frac{1}{2})})+ u^{(i+\frac{1}{2})} \\ 
	- x^{(i+\frac{1}{2})}+(\dr{\partial g})^{-1} u^{(i+\frac{1}{2})} \\
	\end{array} \right)}
	+
	{\left( \begin{array}{cc} 
	\frac{1}{\tau}\mathrm{Id} & -\mathrm{Id} \\ 
	-\mathrm{Id} & \tau\mathrm{Id} \\
	\end{array} \right)}
	{\left( \begin{array}{c}
	x^{(i+\frac{1}{2})}-x^{(i)}  \\ 
	u^{(i+\frac{1}{2})}-u^{(i)}\\
	\end{array} \right)},
\end{equation}
\begin{equation}\label{var_inclusioncpa2}
	\left( \begin{array}{c}
	0  \\ 
	0  \\
	\end{array} \right) 
	\in 
	{\left( \begin{array}{c}
	-(\dg{\partial f})^{-1} (-v^{(i+\frac{1}{2})})+ w^{(i+1)} \\ 
	- v^{(i+\frac{1}{2})}+\dr{\partial g} (w^{(i+1)}) \\
	\end{array} \right)}
	+
	{\left( \begin{array}{cc} 
	\tau\mathrm{Id} & -\mathrm{Id} \\ 
	-\mathrm{Id} & \frac{1}{\tau}\mathrm{Id} \\
	\end{array} \right)}
	{\left( \begin{array}{c}
	v^{(i+\frac{1}{2})}-v^{(i)}  \\ 
	w^{(i+1)}-w^{(i)}\\
	\end{array} \right)}.
\end{equation}
Note that in the ADMM form \eqref{algorithm_admm2}, the final extrapolation step from $v^{(i+\frac{1}{2})}$ to $v^{(i+1)}$ accounts for the absence of relaxation on the variable $w$, in such a way that $w^{(i+1)}-\tau v^{(i+1)}$ takes the appropriate value. That is, $w^{(i+1)}-\tau v^{(i+1)}=w^{(i)}-\tau v^{(i)}+\rho^{(i)}\big((w^{(i+1)}-\tau v^{(i+1/2)})-(w^{(i)}-\tau v^{(i)})\big)$, as if the two variables $v$ and $w$ had been relaxed as usual.

We remark that $\tau$ should not be chosen to be as large as possible in the Douglas--Rachford algorithm or the ADMM: we see in the primal-dual inclusions \eqref{var_inclusioncpa1} and \eqref{var_inclusioncpa2} that the antagonistic values $\tau$ and $1/\tau$ control the primal and dual updates, so there is a tradeoff to achieve in the choice of $\tau$. We insist on the fact that the Douglas--Rachford algorithm and the ADMM are primal-dual algorithms: it is equivalent to apply them on the primal or on the dual problem. There are only two different algorithms to minimize $\ff+\g$: the one given here, with $\mathrm{prox}_{\ff}$ or $\mathrm{prox}_{\dg{f^*}}$ applied first,  and the one obtained by exchanging $\ff$ and $\g$.

Finally, one can find the following relaxed version of the ADMM in the literature:
\begin{align}\label{algorithm_admm3}
&\,\mbox{\textbf{ADMM iteration for \eqref{eqpdr} and \eqref{eqdrd}}: for }i=0,1,\ldots,\notag\\[-1mm]
     &\left\lfloor
    \begin{array}{l}
   x^{(i+\frac{1}{2})}=\mathrm{prox}_{\tau \ff }\big(w^{(i)}-\tilde{v}^{(i)}\big)\\
 \tilde{v}^{(i+1)}=\tilde{v}^{(i)}+\rho\big(x^{(i+\frac{1}{2})}-w^{(i)}\big)\\
w^{(i+1)} = \mathrm{prox}_{ \tau\dr{g}}\big(\tilde{v}^{(i+1)} +x^{(i+\frac{1}{2})}\big),
\end{array}\right.\end{align}
which is proved to converge for $\rho\in \big(0,(\sqrt{5}+1)/2\big)$~\cite{glo84}. This is different from the relaxation considered here. Both types of relaxation have been unified in an extended setting in \cite{he16}.

\subsection{Infimal postcompositions}

Suppose that we want to solve an optimization problem that involves a term $\dr{g(Lx)}$. Throughout the paper, we are studying splitting algorithms, which call the proximity operator of $\g$, $L$, and $L^*$ separately.
Here we look at an alternative, which consists in \emph{a change of variable}: we introduce the variable $r=Lx$ and we express the optimization problem with respect to $r$ instead of $x$. Then an algorithm to solve the reformulated problem will construct a sequence $(r^{(i)})_{i\in\mathbb{N}}$ converging to $r^\star=Lx^\star$, where $x^\star$ is a solution to the initial problem, that will be recovered as a byproduct. For instance, suppose that we want to minimize $\dg{f(x)}+\dr{g(Lx)}$ over $x\in\mathcal{X}$, where $\ff$, $\g$, $L:\mathcal{X}\rightarrow \mathcal{U}$ are as before. 
This is the same as
\begin{align}
&\minimize_{(x,r)\in\mathcal{X}\times\mathcal{U}}\, \dg{f(x)}+\dr{g(r)}\ \ \mbox{s.t.}\ r=Lx\\
\equiv\ &\br{\minimize_{r\in\mathcal{U}}\, \min_{x\in\mathcal{X}}\,\{ \dg{f(x)}:Lx=r\} +\dr{g(r)}}.
\end{align}
Thus, we can rewrite the problem with respect to $r$ only, by ``hiding'' $x$ as a sub-variable, as follows:
\begin{equation}
\minimize\, \dg{\tilde{f}(r)}+\dr{g(r)}, 
\end{equation}
where $\dg{\tilde{f}}:r\in\mathcal{U}\mapsto\inf_{x\in\mathcal{X}} \big(\dg{f}(x) + \imath_{\{r\}}(Lx)\big)$. 
$\dg{\tilde{f}}\in\Gamma_0(\mathcal{U})$ under some mild mathematical safeguards; for instance, we will require that the infimum is attained in the definition of $\dg{\tilde{f}}$.
Clearly, minimizing $\dg{\tilde{f}(r)}+\dr{g(r)}$ will force $r$ to be in the range of $L$, otherwise $\dg{\tilde{f}(r)}=+\infty$. That is, a solution will be $r^\star=Lx^\star$ for some $x^\star\in \argmin_{\{x\in\mathcal{X} : Lx=r^\star\}}\dg{f(x)}$, which is the actual element we are interested in. 
To apply a proximal splitting algorithm to this problem, we do not need to evaluate $\dg{\tilde{f}}$, but we must be able to apply its proximity operator 
\begin{equation}
\mathrm{prox}_\dg{\tilde{f}}:r\in\mathcal{U}\mapsto Lx',\  \mbox{where} \  x'\in\argmin_{x\in\mathcal{X}} \big(\dg{f(x)}+{\textstyle\frac{1}{2}}\|Lx-r\|^2\big).
\end{equation}
That is, $\mathrm{prox}_\dg{\tilde{f}}=L(\dg{\partial f}+L^*L)^{-1}L^*$, which is indeed single-valued: even if $x'$ in its definition is not unique, $Lx'$ is unique. We have $\dg{\tilde{f}^*}:u\in\mathcal{U}\mapsto \sup_{x\in\mathcal{X}} \big(\langle u,Lx\rangle-\dg{f(x)}\big)=\dg{f^*(L^* u)}$ \cite[Proposition 13.24]{bau17}. Thus, the dual problem to $\min_x \big(\dg{f(x)}+\dr{g(Lx)}\big)=\min_r \big(\dg{\tilde{f}(r)}+\dr{g(r)}\big)$ is $\min_u \big(\dg{\tilde{f}^*(-u)}+\dr{g^*(u)}\big)=\min_u \big(\dg{f^*(-L^* u)}+\dr{g^*(u)}\big)$, which is consistent with our formulations of dual problems so far. The function $\tilde{f}$ is called the \emph{infimal postcomposition} of $f$ by $L$, which is in what follows denoted by $\dg{L\triangleright f}$ \cite[Definition 12.34]{bau17}. The property to remember for it is:
\begin{equation}
(\dg{L\triangleright f})^* = \ff^* \circ L^*.
\end{equation}
Thus,  as we have seen, 
behind this technical notion, there is simply a change of variables, which is helpful if 
$\mathrm{prox}_\dg{L\triangleright f}$, or by the Moreau identity $\mathrm{prox}_{\ff^* \circ L^*}$, is easy to compute.\medskip

Let us apply this notion of infimal postcomposition. Let $\mathcal{X}$, $\mathcal{Y}$, $\mathcal{U}$ be real Hilbert spaces. Let $\ff\in \Gamma_0(\mathcal{X})$, $\g\in \Gamma_0(\mathcal{Y})$, $c\in\mathcal{U}$. Let $K:\mathcal{Y}\rightarrow \mathcal{U}$ and $L:\mathcal{X}\rightarrow \mathcal{U}$ be nonzero bounded linear operators. 
We consider the problem:
\begin{equation}
    \minimize_{(x,y)\in\mathcal{X}\times\mathcal{Y}}\, \dg{f(x)} + \dr{g(y)}\quad\mbox{s.t.}\quad Lx+Ky=c.\label{eqadmmip1}
\end{equation}
The dual problem is
\begin{equation}
    \minimize_{u\in\mathcal{U}}\, \dg{f^*(-L^*u)} + \dr{g^*(-K^*u)}+\langle u,c\rangle.\label{eqadmmid1}
\end{equation}
The corresponding primal-dual monotone inclusion is: find $(x,y,u)\in\mathcal{X}\times\mathcal{Y}\times\mathcal{U}$ such that
\begin{equation}\label{eqincip}
\left(\begin{array}{l}
0\\
0\\
  0
  \end{array}\right) \in
\left(\begin{array}{l}
 \dg{\partial f}(x) +L^*u\\
  \dr{\partial g}(y) +K^*u\\
  -Lx-Ky + c
  \end{array}\right).
\end{equation}

For instance, the minimization of $\dg{f(x)}+\dr{g(Lx)}$ is a particular case of this problem, with $c=0$, $\mathcal{U}=\mathcal{Y}$, $K=\mathrm{-Id}$. We want to solve \eqref{eqadmmip1} by applying two infimal postcompositions, so that the problem amounts to finding an element $r^\star=Lx^\star=-Ky^\star+c$, where $(x^\star,y^\star)$ is a solution to  \eqref{eqadmmip1}. That is, we rewrite \eqref{eqadmmip1} as 
\begin{equation}
    \minimize_{r\in\mathcal{U}}\, \dg{L\triangleright f(r)} + \dr{(-K)\triangleright g(r-c),}\label{eqadmmipa}
\end{equation}
where
\begin{align}
   \dg{L\triangleright f}:&r\in\mathcal{U}\mapsto\inf_{x\in\mathcal{X}} \big(\dg{f}(x) + \imath_{\{0\}}(Lx-r)\big)\label{eqipf1}\\
    \dr{(-K)\triangleright g}:&r\in\mathcal{U}\mapsto\inf_{y\in\mathcal{X}} \big(\dr{g}(y) + \imath_{\{0\}}(Ky+r)\big)\label{eqipf2}.
\end{align}
For the problem to be well posed, we suppose that the following assumptions hold:\medskip

\noindent (i) The solution set of \eqref{eqincip} is nonempty.

\noindent (ii) The infimal postcompositions are exact~\cite[Definition 12.34]{bau17}; that is, for every $r\in\mathcal{U}$, the infimum in 
\eqref{eqipf1} and \eqref{eqipf2} is attained, so it is a minimum (its value can be $+\infty$). In other words, for every $r\in\mathrm{ran}\, L$, the set of minimizers of $\ff$ over the set $L^{-1}\{r\}$ is nonempty. Likewise, for every $r\in\mathrm{ran}\,K$, the set of minimizers of $\g$ over $K^{-1}\{-r\}$ is nonempty.

\noindent (iii) $\dg{L\triangleright f}$ is lower semicontinuous. A sufficient condition for that is $0 \in \mathrm{sri} (\mathrm{ran}\,L^* - \mathrm{dom}\, \dg{f^*})$~\cite[Corollary 25.44]{bau17}. Likewise, $\dr{(-K)\triangleright g}$ is lower semicontinuous.\medskip

Note that $\dg{L\triangleright f}$ and $\dr{(-K)\triangleright g}$ are convex and proper~\cite[Proposition 12.36]{bau17}. So (iii) implies that $\dg{L\triangleright f}\in\Gamma_0(\mathcal{U})$ and $\dr{(-K)\triangleright g}\in\Gamma_0(\mathcal{U})$. Moreover, $\dg{(L\triangleright f)^*}=\dg{f^*}\circ L^*$ and  $\dr{((-K)\triangleright g)^*}=\dr{g^*}\circ -K^*$~\cite[Proposition 13.24]{bau17}.

Then, let us apply the ADMM \eqref{algorithm_admm2} to the problem \eqref{eqadmmipa}:
\ceq{align}{\label{algorithm_admm2inf}
&\,\mbox{\textbf{ADMM iteration for \eqref{eqadmmip1} and \eqref{eqadmmid1}}: for }i=0,1,\ldots,\notag\\[-1mm]
     &\left\lfloor
    \begin{array}{l}
   x^{(i+\frac{1}{2})}\in\argmin_{x\in\mathcal{X}} \big(\tau\dg{f(x)}+\frac{1}{2}\|Lx+Ky^{(i)}-c+\tilde{v}^{(i)}\|^2\big)\\
 \tilde{v}^{(i+\frac{1}{2})}=\tilde{v}^{(i)}+Lx^{(i+\frac{1}{2})}+Ky^{(i)}-c\\
y^{(i+1)}  \in\argmin_{y\in\mathcal{Y}} \big(\tau\g(y)+\frac{1}{2}\|Lx^{(i+\frac{1}{2})}+Ky-c+\tilde{v}^{(i+\frac{1}{2})} \|^2\big)\\
  \tilde{v}^{(i+1)}=\tilde{v}^{(i+\frac{1}{2})}+(\rho^{(i)}-1)\big(Lx^{(i+\frac{1}{2})}+Ky^{(i+1)}-c\big)
\end{array}\right.}
($x^{(i+\frac{1}{2})}$ and $w^{(i)}$ in \eqref{algorithm_admm2} become $Lx^{(i+\frac{1}{2})}$ and $-Ky^{(i)}+c$ in \eqref{eqadmmipa}, respectively).

Since the ADMM and the Douglas--Rachford algorithm are equivalent, we can write the iteration like in \eqref{algorithm_22} instead:
\ceq{align}{\label{algorithm_22ip}
&\,\mbox{\textbf{Douglas--Rachford iteration for \eqref{eqadmmip1}}: for }i=0,1,\ldots,\notag\\[-1mm]
    &\left\lfloor
    \begin{array}{l}
    x^{(i+\frac{1}{2})}\in\argmin_{x\in\mathcal{X}}\big( \tau\ff(x)+\frac{1}{2}\|Lx-s^{(i)}\|^2\big)\\
   y^{(i+1)}  \in\argmin_{y\in\mathcal{Y}}\big( \tau\g(y)+\frac{1}{2}\|Ky-c+2Lx^{(i+\frac{1}{2})} - s^{(i)}\|^2\big)\\
s^{(i+1)}=s^{(i)}-\rho^{(i)}\big(  Lx^{(i+\frac{1}{2})}+Ky^{(i+1)}-c\big).
\end{array}\right.}
In the case $\rho^{(i)}>1$, the Douglas--Rachford form \eqref{algorithm_22ip} is to be preferred over the ADMM form \eqref{algorithm_admm2inf}, since it involves fewer operations.

As an application of \cref{thm46}, we have:

\begin{theorem}[ADMM \eqref{algorithm_admm2inf} or Douglas--Rachford algorithm \eqref{algorithm_22ip}]
Let $y^{(0)}\in\mathcal{Y}$, $\tilde{v}^{(0)}\in\mathcal{U}$, $s^{(0)}\in\mathcal{U}$ be such that $s^{(0)}=-Ky^{(0)}+c-\tilde{v}^{(0)}$. 
Let $\tau >0$ and let $(\rho^{(i)})_{i\in\mathbb{N}}$ be a sequence in $[0,2]$ such that $\sum_{i\in \mathbb{N}} \rho^{(i)} (2-\rho^{(i)})=+\infty$.  Then the sequences $(Lx^{(i+\frac{1}{2})})_{i\in \mathbb{N}}$ and $(-Ky^{(i)}+c)_{i\in \mathbb{N}}$ defined by the iteration \eqref{algorithm_admm2inf}, or equivalently by the iteration \eqref{algorithm_22ip}, both converge weakly to some element $r^\star\in\mathcal{X}$ solution to \eqref{eqadmmipa}. Moreover, the sequences $(\tilde{v}^{(i+\frac{1}{2})}/\tau)_{i\in \mathbb{N}}$ and $(\tilde{v}^{(i)}/\tau)_{i\in \mathbb{N}}$  defined by the iteration \eqref{algorithm_admm2inf} both converge weakly to some element $u^\star\in\mathcal{X}$ solution to \eqref{eqadmmid1}.  In addition, $(Lx^{(i+\frac{1}{2})}+Ky^{(i+1)}-c)_{i\in \mathbb{N}}$ converges strongly to $0$.

Furthermore, suppose that $\mathcal{X}$, $\mathcal{Y}$, $\mathcal{U}$ are of finite dimension and that, for every $r\in\mathcal{U}$, the minimizers of $\tau \ff+\frac{1}{2}\|L\cdot-r\|^2$ and $\tau \g+\frac{1}{2}\|K\cdot+r\|^2$ are unique. Then the sequence $(x^{(i+\frac{1}{2})},y^{(i)})_{i\in \mathbb{N}}$ defined by the iteration \eqref{algorithm_admm2inf}, or equivalently by the iteration \eqref{algorithm_22ip}, converges to a solution of \eqref{eqadmmip1}.\end{theorem}

\begin{proof}Let us prove the second part of the theorem. The operators $(\tau\dg{\partial f}+L^*L)^{-1}L^*$ and $(\tau\dr{\partial g}+K^*K)^{-1}K^*$ are supposed single-valued on $\mathcal{U}$. Therefore, $(\tau\dg{\partial f}+L^*L)^{-1}$ is single-valued on $\mathrm{ran}\,L^*$ and, since it is maximally monotone, it is continuous on $\mathrm{ran}\,L^*$. Likewise, $(\tau\dr{\partial g}+K^*K)^{-1}$ is continuous on $\mathrm{ran}\,K^*$. Hence, $(x^{(i+\frac{1}{2})})_{i\in \mathbb{N}}$ converges to $x^\star=(\tau\dg{\partial f}+L^*L)^{-1}L^*(r^\star-\tau u^\star)$ and $(y^{(i)})_{i\in \mathbb{N}}$ converges to $y^\star=(\tau\dr{\partial g}+K^*K)^{-1}K^*(-r^\star+c-\tau u^\star)$. Since $Lx^\star=r^\star$ and $-Ky^\star+c=r^\star$, we have $0\in \tau\dg{\partial f}(x^\star)+\tau L^* u^\star$ and $0\in \tau\dr{\partial g}(y^\star)+\tau K^* u^\star$. In addition, $Lx^\star+Ky^\star=c$. Therefore, $(x^\star,y^\star,u^\star)$ is a solution to \eqref{eqincip}, which implies that $(x^\star,y^\star)$ is a solution to \eqref{eqadmmip1}.
\end{proof}\smallskip

Here we have used infimal postcompositions in the Douglas--Rachford algorithm, but they could be used in any other splitting algorithm. One could study, on a case-by-case basis, whether the subproblems corresponding to the proximity operators of the infimal postcompositions are easy to solve, and if so, whether it is better to do so or to split the problem with separate calls to the linear operators.

\section{The Generalized Chambolle--Pock Algorithm}\label{sec5}

Let $\mathcal{X}$, $\mathcal{U}$, $\mathcal{V}$ be real Hilbert spaces. Let $\ff\in \Gamma_0(\mathcal{V})$, $\g\in \Gamma_0(\mathcal{U})$, $c\in\mathcal{X}$. Let $K:\mathcal{X}\rightarrow \mathcal{V}$ and $L:\mathcal{X}\rightarrow \mathcal{U}$ be nonzero bounded linear operators. 
We consider the following problem, generalizing \eqref{eq22}:
\begin{equation}
    \minimize_{x\in\mathcal{X}}\, \dg{f(Kx)} + \dr{g(Lx)} + \langle x,c\rangle.\label{eq29}
\end{equation}
The corresponding monotone inclusion, which we will actually solve, is
\begin{equation}
   0\in K^*\dg{\partial f} (Kx) + L^*\dr{\partial g} (Lx) +c.
\end{equation}
We introduce two dual variables $u\in  \mathcal{U}$ and $v\in  \mathcal{V}$, so that the problem is to find $(x,u,v)$ such that 
\begin{equation}
\left\{\begin{array}{l}
v\in \dg{\partial f} (Kx)\\
  u\in\dr{\partial g} (Lx)\\
  0 \in K^*v+L^*u + c
  \end{array}\right. .
\end{equation}
Accordingly, the dual problem is 
\begin{equation}
    \minimize_{(u,v)\in\mathcal{U}\times\mathcal{V}}\, \dg{f^*(v)} + \dr{g^*(u)}\quad\mbox{s.t.}\quad K^*v+L^*u+c=0.\label{eq292}
\end{equation}
Clearly, if $K=\mathrm{Id}$ and $L=\mathrm{Id}$, \eqref{eq29} and \eqref{eq292} revert to \eqref{eq27} and \eqref{eq28}, and we can use the Douglas--Rachford algorithm \eqref{algorithm_28} to solve these problems. On the other hand, if $K=\mathrm{Id}$ and $c=0$,  \eqref{eq29} and \eqref{eq292} revert to \eqref{eq22} and \eqref{eq23}, and we can use the Chambolle--Pock algorithm  \eqref{algoCP} or \eqref{algoCP2}.\\

Let $\tau>0$, $\sigma>0$, $\eta>0$, let $x^{(0)}\in\mathcal{X}$, $u^{(0)}\in\mathcal{U}$, $v^{(0)}\in\mathcal{V}$, and let $(\rho^{(i)})_{i\in\mathbb{N}}$ be a sequence of relaxation parameters. We consider the following algorithm:
\ceq{align}{\label{eqhybalg}
&\,\mbox{\textbf{Generalized Chambolle--Pock iteration for \eqref{eq29} and \eqref{eq292}}:}\notag\\
&\,\mbox{ for }i=0,1,\ldots,\notag\\[-1mm]
    &\left\lfloor
\begin{array}{l}
v^{(i+\frac{1}{2})}
= \mathrm{prox}_{\eta \dg{f^*}}\Big(v^{(i)}+ \eta K 
\big( x^{(i)}-\tau (L^* u^{(i)}+K^* v^{(i)}+c )\big)\Big)\\
x^{(i+\frac{1}{2})}
= x^{(i)} -\tau (L^* u^{(i)}+K^* v^{(i+\frac{1}{2})}+ c)\\
u^{(i+\frac{1}{2})}
= \mathrm{prox}_{\sigma \dr{g^*}}\big(u^{(i)}+ \sigma L 
(2x^{(i+\frac{1}{2})}-x^{(i)})\big)\\
\hphantom{u^{(i+\frac{1}{2})}}
= \mathrm{prox}_{\sigma \dr{g^*}}\Big(u^{(i)}+ \sigma L 
\big(x^{(i+\frac{1}{2})}-\tau(L^* u^{(i)}+K^* v^{(i+\frac{1}{2})}+c)\big) \Big)\\
x^{(i+1)}=x^{(i)}+\rho^{(i)} (x^{(i+\frac{1}{2})}-x^{(i)})\\
u^{(i+1)}=u^{(i)}+\rho^{(i)} (u^{(i+\frac{1}{2})}-u^{(i)})\\
v^{(i+1)}=v^{(i)}+\rho^{(i)} (v^{(i+\frac{1}{2})}-v^{(i)}).
\end{array}\right.}
The primal-dual inclusion satisfied at every iteration is
\begin{align}\label{algohybrid}
\left(\begin{array}{l}
0\\
0\\
0
\end{array}\right)&\in
\left(\begin{array}{l}
 c+L^*u^{(i+\frac{1}{2})}+K^*v^{(i+\frac{1}{2})}\\
-Lx^{(i+\frac{1}{2})}+ \dr{\partial  g^*}(u^{(i+\frac{1}{2})})\\
-Kx^{(i+\frac{1}{2})}+ \dg{\partial  f^*}(v^{(i+\frac{1}{2})})\\
\end{array}\right)\\
&\qquad+
\underbrace{\left(\begin{array}{ccc}
\frac{1}{\tau}\mathrm{Id}&-L^*&0\\
-L&\frac{1}{\sigma}\mathrm{Id}&0\\
0& 0&\frac{1}{\eta}\mathrm{Id}-\tau KK^*
\end{array}\right)}_P
\left(\begin{array}{c}
x^{(i+\frac{1}{2})}-x^{(i)}\\
u^{(i+\frac{1}{2})}-u^{(i)}\\
v^{(i+\frac{1}{2})}-v^{(i)}
\end{array}\right).\notag
\end{align}
$P$ is strongly positive if and only if
\begin{equation}
    \tau\sigma\|L\|^2<1\quad\mbox{and}\quad\tau\eta\|K\|^2<1.\label{eqcondh2}
\end{equation}
Therefore, since the algorithm is a preconditioned primal-dual proximal point algorithm, we can apply  \cref{thm23} and we obtain:

\begin{theorem}[Generalized Chambolle--Pock algorithm \eqref{eqhybalg}]\label{thm51}
Let $x^{(0)}\in\mathcal{X}$, $u^{(0)}\in\mathcal{U}$, $v^{(0)}\in\mathcal{V}$. 
Let $\tau >0$, $\sigma>0$, $\eta >0$ be such that $\tau\sigma\|L\|^2<1$ and $\tau\eta\|K\|^2<1$. Let $(\rho^{(i)})_{i\in\mathbb{N}}$ be a sequence in $[0,2]$ such that $\sum_{i\in \mathbb{N}} \rho^{(i)} (2-\rho^{(i)})=+\infty$.  Then the sequences $(x^{(i)})_{i\in \mathbb{N}}$ and $(u^{(i)},v^{(i)})_{i\in \mathbb{N}}$ defined by the iteration \eqref{eqhybalg} converge weakly to a solution of \eqref{eq29} and a solution of \eqref{eq292}, respectively.
\end{theorem}

The Generalized Chambolle--Pock algorithm \eqref{eqhybalg} has appeared in the literature under different names, such as  
 Alternating Proximal Gradient Method~\cite{ma16}, Generalized Alternating Direction Method of Multipliers~\cite{den15}, or Preconditioned ADMM~\cite{bre17,sun19}. The convergence results derived in this section generalize previously known results.

If $K=\mathrm{Id}$ and we set $\eta=1/\tau$, we recover the Chambolle--Pock algorithm form I, to minimize $\dg{\tilde{f}}+\g\circ L$ with $\dg{\tilde{f}}=\ff + \langle\cdot,c\rangle$. 
Indeed, in that case, the first two updates become:
\begin{align}
    &\left|
\begin{array}{l}
v^{(i+\frac{1}{2})}
= \mathrm{prox}_{ \dg{f^*}/\tau}\big(  
 x^{(i)}/\tau- L^* u^{(i)}-c \big)\\
x^{(i+\frac{1}{2})}
= x^{(i)} -\tau L^* u^{(i)}-\tau v^{(i+\frac{1}{2})}-\tau c\\
\hphantom{x^{(i+\frac{1}{2})}}=\mathrm{prox}_{ \tau \ff}\big(  
 x^{(i)}- \tau L^* u^{(i)}-\tau c \big)\\
 \hphantom{x^{(i+\frac{1}{2})}}=\mathrm{prox}_{ \tau \dg{\tilde{f}}}\big(  
 x^{(i)}- \tau L^* u^{(i)} \big).
\end{array}\right. 
\end{align}
Therefore, the Generalized Chambolle--Pock algorithm indeed generalizes the Chambolle--Pock algorithm to any linear operator $K$. But we remark that setting $\eta=1/\tau$ is not allowed in \cref{thm51}. So, to extend the range of parameters to $\tau\sigma\|L\|^2\leq 1$ and $\tau\eta\|K\|^2\leq1$, we now analyze the algorithm from another point of view: we show that it is a Douglas--Rachford algorithm \br{in an augmented space. This analysis is inspired by and generalizes the  analysis of the Chambolle--Pock algorithm by O'Connor and Vandenberghe~\cite{oco18}.}

Let $\tau>0$, $\sigma>0$, $\eta>0$ be such that $\tau\sigma\|L\|^2\leq 1$ and $\tau\eta\|K\|^2\leq1$.
Let  $A$ be a linear operator from $\mathcal{A}$ to $\mathcal{V}$, for some 
real Hilbert space $\mathcal{A}$, such that $KK^*+AA^*=(\tau\eta)^{-1}\mathrm{Id}$; for instance, $A=\sqrt{(\tau\eta)^{-1}\mathrm{Id}-KK^*}$ is a valid choice. 
We do not need to exhibit $A$; the fact that it exists is sufficient here. Similarly, let  $B$ be a linear operator from $\mathcal{B}$ to $\mathcal{U}$ for some 
real Hilbert space $\mathcal{B}$, such that $LL^*+BB^*=(\tau\sigma)^{-1}\mathrm{Id}$. We introduce the real Hilbert space $\mathcal{Z}=\mathcal{X}\times\mathcal{A}\times\mathcal{B}$ and the following functions of  $\Gamma_0(\mathcal{Z})$:
\begin{align}
&F:(x,a,b)\in\mathcal{Z}\mapsto \ff(Kx+Aa)+\imath_{\{0\}}(b),\\
&G:(x,a,b)\in\mathcal{Z}\mapsto \g(Lx+Bb)+\imath_{\{0\}}(a).
\end{align}
Then we can rewrite the problem \eqref{eq29} as:
\begin{equation}
    \minimize_{z\in\mathcal{Z}}\, F(z) + G(z) + \langle z,(c,0,0)\rangle.\label{eq293}
\end{equation}
We can now apply the Douglas--Rachford algorithm \eqref{algorithm_28} in the augmented space $\mathcal{Z}$. For this, we need to observe that the proximity operators of $F^*$ and $G^*$ are easy to compute: for any $\tau>0$, we have \cite[equation~15]{oco18}
\begin{align}
&\mathrm{prox}_{F^*/\tau}:(x,a,b)\in\mathcal{Z}\mapsto (K^*v,A^*v,b)\ \mbox{with}\ v=\mathrm{prox}_{\eta \dg{f^*}}(\tau\eta Kx+\tau\eta Aa)\\
&\mathrm{prox}_{G^*/\tau}:(x,a,b)\in\mathcal{Z}\mapsto (L^*u,a,B^*u)\ \mbox{with}\ u=\mathrm{prox}_{\sigma \dr{g^*}}(\tau\sigma Lx+\tau\sigma Bb).
\end{align}
After some substitutions, notably replacing $AA^*$ by $(\tau\eta)^{-1}\mathrm{Id}-K^*K$ and $BB^*$ by $(\tau\sigma)^{-1}\mathrm{Id}-L^*L$, we recover exactly the algorithm in \eqref{eqhybalg}.

Hence, as an application of \cref{thm45}, we obtain:

\begin{theorem}[Generalized Chambolle--Pock algorithm \eqref{eqhybalg}]\label{thm52}
Let $x^{(0)}\in\mathcal{X}$, $u^{(0)}\in\mathcal{U}$, $v^{(0)}\in\mathcal{V}$. 
Let $\tau >0$, $\sigma>0$, $\eta >0$ be such that $\tau\sigma\|L\|^2\leq 1$ and $\tau\eta\|K\|^2\leq 1$. Let $(\rho^{(i)})_{i\in\mathbb{N}}$ be a sequence in $[0,2]$ such that $\sum_{i\in \mathbb{N}} \rho^{(i)} (2-\rho^{(i)})=+\infty$.  Then the sequences $(x^{(i)})_{i\in \mathbb{N}}$ and $(u^{(i)},v^{(i)})_{i\in \mathbb{N}}$ defined by the iteration \eqref{eqhybalg} converge \br{weakly} to a solution of \eqref{eq29} and a solution of \eqref{eq292}, respectively. Moreover, 
$(L^*u^{(i+\frac{1}{2})}+K^*v^{(i+\frac{1}{2})}+c)_{i\in \mathbb{N}}$ converges strongly  to $0$.\end{theorem}

Thus, in practice, one can keep $\tau$ as the single parameter to tune and set $\sigma=1/(\tau\|L\|^2)$ and $\eta=1/(\tau\|K\|^2)$.

We can write the Generalized Chambolle--Pock algorithm in a different form with only one call to $K$, $K^*$, $L$, $L^*$ per iteration. 
For this, we introduce the scaled variable $\tilde{x}=x/\tau$ and auxiliary variables  $b=K^*v$ and $r=L^* u+K^* v+c$.
Set $b^{(0)}=K^*v^{(0)}$ and $l^{(0)}=L^* u^{(0)}+b^{(0)}+c$. 
\br{This yields the iteration:}
\begin{align}
    &\left\lfloor
\begin{array}{l}
v^{(i+\frac{1}{2})}
= \mathrm{prox}_{\eta \dg{f^*}}\big(v^{(i)}+ \eta \tau K 
( \tilde{x}^{(i)}- l^{(i)} )\big)\\
l^{(i+\frac{1}{2})}=l^{(i)}+K^* v^{(i+\frac{1}{2})}-b^{(i)}\\
u^{(i+\frac{1}{2})}
= \mathrm{prox}_{\sigma \dr{g^*}}\big(u^{(i)}+ \sigma \tau L 
(\tilde{x}^{(i)}-2l^{(i+\frac{1}{2})})\big)\\
\tilde{x}^{(i+1)}=\tilde{x}^{(i)}-\rho^{(i)} l^{(i+\frac{1}{2})}\\
u^{(i+1)}=u^{(i)}+\rho^{(i)} (u^{(i+\frac{1}{2})}-u^{(i)})\\
v^{(i+1)}=v^{(i)}+\rho^{(i)} (v^{(i+\frac{1}{2})}-v^{(i)})\\
b^{(i+1)}=b^{(i)}+\rho^{(i)} (l^{(i+\frac{1}{2})}-l^{(i)})\\
l^{(i+1)}=L^*u^{(i+1)}+b^{(i+1)}+c.
\end{array}\right.
\end{align}

Furthermore, let us show that we recover as a particular case the Loris--Verhoeven algorithm applied to the minimization of $\ff\circ K + \h$, when $\h:x\mapsto\frac{1}{2}\langle x,\bl{Q}x\rangle + \langle x,c\rangle$ is quadratic with $\bl{Q}=L^*L$ (again, given $\bl{Q}$, such an operator $L$ exists). Indeed, minimizing $\ff\circ K + \h$ is equivalent to minimizing $\ff\circ K + \g\circ L +  \langle \cdot,c\rangle$, with $\g=\dr{g^*}=\frac{1}{2}\|\cdot\|^2$. Let us apply the generalized Chambolle--Pock algorithm \eqref{eqhybalg} on this latter problem, with $\sigma=1$. Since $\mathrm{prox}_{\sigma \dr{g^*}}=\frac{1}{2}\mathrm{Id}$, we can write the update of $u$ as:
\begin{equation}
\textstyle
u^{(i+\frac{1}{2})}=\frac{1}{2}u^{(i)}+Lx^{(i+\frac{1}{2})}-\frac{1}{2}Lx^{(i)},\label{eqlvgcp1}
\end{equation}
so that if $u^{(0)}=Lx^{(0)}$, we have $u^{(i)}=Lx^{(i)}$, for every $i\in\mathbb{N}$. Hence, we can remove the variable $u$ and rewrite the iteration as:
\begin{align}\label{eqlvgcp2}
    &\left\lfloor
\begin{array}{l}
v^{(i+\frac{1}{2})}
= \mathrm{prox}_{\eta \dg{f^*}}\Big(v^{(i)}+ \eta K 
\big( x^{(i)}-\tau (L^* Lx^{(i)}+K^* v^{(i)}+c )\big)\Big)\\
x^{(i+\frac{1}{2})}
= x^{(i)} -\tau (L^* Lx^{(i)}+K^* v^{(i+\frac{1}{2})}+ c)\\
x^{(i+1)}=x^{(i)}+\rho^{(i)} (x^{(i+\frac{1}{2})}-x^{(i)})\\
v^{(i+1)}=v^{(i)}+\rho^{(i)} (v^{(i+\frac{1}{2})}-v^{(i)}),
\end{array}\right.\end{align}
which is exactly the Loris--Verhoeven iteration \eqref{algorithm_1}. Thus, the Loris--Verhoeven algorithm can be viewed as a primal-dual forward-backward algorithm, but also as a primal-dual Douglas--Rachford algorithm, when $\h$ is quadratic.

\section{The Condat--V\~u Algorithm}\label{sec6}

Let us consider the primal optimization problem:
\begin{equation}
    \minimize_{x\in\mathcal{X}}\, \dg{f(x)} + \dr{g(Lx)} + \bl{h(x)},\label{eq18}
\end{equation}
where  $\ff \in \Gamma_0(\mathcal{X})$, $\g\in \Gamma_0(\mathcal{U})$,  $\h:\mathcal{X}\rightarrow \mathbb{R}$ is a convex and 
differentiable function with  $\beta$-Lipschitz continuous gradient $\bl{\nabla h}$, for some real $\beta> 0$, and $L:\mathcal{X}\rightarrow \mathcal{U}$ is a bounded linear operator. The corresponding monotone inclusion is
\begin{equation}
    0 \in \dg{\partial f}(x)+L^*\dr{\partial g} (Lx) + \bl{\nabla h}(x).
    \label{eq19}
\end{equation}
Again, we introduce a dual variable $u$, so that we can rewrite the problem \eqref{eq19} with respect to 
 a pair of objects $z=(x,u)$   in 
$\mathcal{Z}=\mathcal{X}\times\mathcal{U}$:
\begin{equation}\label{eqvarcv}
   \left( \begin{array}{c}
   0\\0
   \end{array}\right)\in
    \underbrace{\left( \begin{array}{c}
   \dg{\partial f}(x)+ L^*u \\
   -Lx + (\dr{\partial g})^{-1}(u)
   \end{array}\right)}_{Mz}+ \underbrace{\left( \begin{array}{c}
     \bl{\nabla h}(x) \\
  0
   \end{array}\right)}_{\bl{C}z} .
\end{equation}
If $(x,u)\in\mathcal{X}\times\mathcal{U}$ is a solution to \eqref{eqvarcv}, then $x$ is a solution to \eqref{eq18} and $u\in \dr{\partial g}(Lx)$ is a solution to the dual problem
\begin{equation}
    \minimize_{u\in\mathcal{U}}\, (\ff+\h)^*(-L^*u) + \dr{g^*(u)}.\label{eq30}
\end{equation}

The operator $M:\mathcal{Z}\rightarrow 2^\mathcal{Z},(x,u)\mapsto (\dg{\partial f}(x)+L^*u, -Lx + (\dr{\partial g})^{-1}u)$ is maximally monotone~\cite[Proposition 26.32 (iii)]{bau17} and $\bl{C}:\mathcal{Z}\rightarrow\mathcal{Z},(x,u)\mapsto (\bl{\nabla h}(x),0)$ is $\xi$-cocoercive with $\xi=1/\beta$. Thus, it is again natural to think of the forward-backward iteration with preconditioning. The difference to the construction in \cref{sec3} is the presence of the nonlinear operator $\dg{\partial f}$, which prevents us from expressing $x^{(i+\frac{1}{2})}$ in terms of $x^{(i)}$ and $u^{(i+\frac{1}{2})}$. Instead, the iteration is made explicit by canceling the dependence of $x^{(i+\frac{1}{2})}$ on $u^{(i+\frac{1}{2})}$ in $P$. That is, the iteration, written in implicit form, is
\begin{align}\label{var_inclusioncv}
	\left( \begin{array}{c}
	0  \\ 
	0  \\
	\end{array} \right) 
	&\in 
	\underbrace{\left( \begin{array}{c}
	\dg{\partial f} (x^{(i+\frac{1}{2})})+L^* u^{(i+\frac{1}{2})} \\ 
	-L x^{(i+\frac{1}{2})}+(\dr{\partial g})^{-1} u^{(i+\frac{1}{2})} \\
	\end{array} \right)}_{Mz^{(i+\frac{1}{2})}}
	+
	\underbrace{\left( \begin{array}{c}
	\bl{\nabla h}(x^{(i)})  \\ 
	0  \\
	\end{array} \right)}_{\bl{C}z^{(i)}}\\
	&\qquad\quad+
	\underbrace{\left( \begin{array}{cc} 
	\frac{1}{\tau}\mathrm{Id} & -L^* \\ 
	-L & \frac{1}{\sigma}\mathrm{Id} \\
	\end{array} \right)}_{P}
	\underbrace{\left( \begin{array}{c}
	x^{(i+\frac{1}{2})}-x^{(i)}  \\ 
	u^{(i+\frac{1}{2})}-u^{(i)}\\
	\end{array} \right)}_{z^{(i+\frac{1}{2})}-z^{(i)}},\notag
\end{align}
where   $\tau>0$ and $\sigma>0$ are two real parameters, $z^{(i)}=(x^{(i)},u^{(i)})$ and $z^{(i+\frac{1}{2})}=(x^{(i+\frac{1}{2})},u^{(i+\frac{1}{2})})$.
Thus, 
the primal-dual forward-backward iteration is: 
\ceq{align}{\label{algocv}
&\,\mbox{\textbf{Condat--V\~u iteration, form I, for \eqref{eq18} and \eqref{eq30}}: for }i=0,1,\ldots,\notag\\[-1mm]
    &\left\lfloor
    \begin{array}{l}
x^{(i+\frac{1}{2})}= \mathrm{prox}_{\tau \ff}\big(x^{(i)}-\tau\bl{\nabla h}(x^{(i)})-\tau L^*u^{(i)} \big)\\
    u^{(i+\frac{1}{2})}=\mathrm{prox}_{\sigma \dr{g^*}}\big(u^{(i)}+\sigma L (2x^{(i+\frac{1}{2})}-x^{(i)})\big)\\
    x^{(i+1)}=x^{(i)}+\rho^{(i)}(x^{(i+\frac{1}{2})}-x^{(i)})\\
    u^{(i+1)}=u^{(i)}+\rho^{(i)} (u^{(i+\frac{1}{2})}-u^{(i)}).
\end{array}\right.}%
This algorithm was proposed independently by the first author~\cite{con13} and by V\~u~\cite{vu13}. 

An alternative is to update $u$ before $x$, instead of vice versa. This yields a different algorithm, characterized by the primal-dual inclusion
\begin{align}\label{var_inclusioncv2}
	\left( \begin{array}{c}
	0  \\ 
	0  \\
	\end{array} \right) 
	&\in 
	\underbrace{\left( \begin{array}{c}
	\dg{\partial f} (x^{(i+\frac{1}{2})})+L^* u^{(i+\frac{1}{2})} \\ 
	-L x^{(i+\frac{1}{2})}+(\dr{\partial g})^{-1} u^{(i+\frac{1}{2})} \\
	\end{array} \right)}_{Mz^{(i+\frac{1}{2})}}
	+
	\underbrace{\left( \begin{array}{c}
	\bl{\nabla h}(x^{(i)})  \\ 
	0  \\
	\end{array} \right)}_{\bl{C}z^{(i)}}\\
	&\qquad\quad+
	\underbrace{\left( \begin{array}{cc} 
	\frac{1}{\tau}\mathrm{Id} & L^* \\ 
	L & \frac{1}{\sigma}\mathrm{Id} \\
	\end{array} \right)}_{P}
	\underbrace{\left( \begin{array}{c}
	x^{(i+\frac{1}{2})}-x^{(i)}  \\ 
	u^{(i+\frac{1}{2})}-u^{(i)}\\
	\end{array} \right)}_{z^{(i+\frac{1}{2})}-z^{(i)}}.\notag
\end{align}
This corresponds to the primal-dual forward-backward iteration:
\ceq{align}{\label{algocv2}
&\,\mbox{\textbf{Condat--V\~u iteration, form II, for \eqref{eq18} and \eqref{eq30}}: for }i=0,1,\ldots,\notag\\[-1mm]
    &\left\lfloor
    \begin{array}{l}
    u^{(i+\frac{1}{2})}=\mathrm{prox}_{\sigma \dr{g^*}}\big(u^{(i)}+\sigma L x^{(i)}\big)\\
x^{(i+\frac{1}{2})}= \mathrm{prox}_{\tau \ff}\big(x^{(i)}-\tau\bl{\nabla h}(x^{(i)})-\tau L^*(2u^{(i+\frac{1}{2})}-u^{(i)}) \big)\\
     u^{(i+1)}=u^{(i)}+\rho^{(i)} (u^{(i+\frac{1}{2})}-u^{(i)})\\
    x^{(i+1)}=x^{(i)}+\rho^{(i)}(x^{(i+\frac{1}{2})}-x^{(i)}).
\end{array}\right.}

As an application of \cref{thm22}, we obtain the following convergence result~\cite[Theorem 3.1]{con13}:

\begin{theorem}[Condat--V\~u algorithm \eqref{algocv} or \eqref{algocv2}]\label{thm61}
Let $x^{(0)}\in\mathcal{X}$ and 
$u^{(0)}\in\mathcal{U}$. Let $\tau >0$ and $\sigma>0$  be such that $\tau\big(\sigma\|L\|^2+\frac{\beta}{2}\big)<1$. Set $\delta=2-\frac{\beta}{2}\left(\frac{1}{\tau}-\sigma\|L\|^{2}\right)^{-1}>1$.
Let $(\rho^{(i)})_{i\in\mathbb{N}}$ be a sequence in $[0,\delta]$ such that $\sum_{i\in \mathbb{N}} \rho^{(i)} (\delta-\rho^{(i)})=+\infty$. Then the sequences $(x^{(i)})_{i\in \mathbb{N}}$ and $(u^{(i)})_{i\in \mathbb{N}}$, defined either by the iteration \eqref{algocv} or by the iteration \eqref{algocv2}, converge weakly to a solution of \eqref{eq18} and a solution of \eqref{eq30}, respectively.\end{theorem}

\begin{proof}In view of \eqref{var_inclusioncv} and \eqref{var_inclusioncv2}, this is \cref{thm22}  applied to the problem \eqref{eqvarcv}. The condition on $\tau$ and $\sigma$ implies that $\sigma\tau\|L\|^2<1$, so that $P$ is strongly positive  
by virtue of the properties of the Schur complement. Let us establish the cocoercivity of $P^{-1}\bl{C}$ in $\mathcal{Z}_P$. 
In both cases  \eqref{var_inclusioncv} and \eqref{var_inclusioncv2}, we have, for every $z=(x,u)$ and $z'=(x',u')$ in $\mathcal{Z}$, 
\begin{align}
    \|P^{-1}\bl{C}z-P^{-1}\bl{C}z'\|^2_P&= \langle P^{-1}\bl{C}z-P^{-1}\bl{C}z',\bl{C}z-\bl{C}z' \rangle\label{eqcoer1}\\
    &=\left\langle{\textstyle\frac{1}{\sigma} \left(\frac{1}{\sigma\tau}\mathrm{Id}-L^*L\right)^{-1}}\big(\bl{\nabla h}(x)-\bl{\nabla h}(x')\big),\bl{\nabla h}(x)-\bl{\nabla h}(x')\right\rangle\\
    &\leq {\textstyle \left(\frac{1}{\tau}-\sigma\|L\|^2\right)^{-1}}\big\langle\bl{\nabla h}(x)-\bl{\nabla h}(x'),\bl{\nabla h}(x)-\bl{\nabla h}(x')\big\rangle\\
    &\leq {\textstyle\beta \left(\frac{1}{\tau}-\sigma\|L\|^2\right)^{-1}}\big\langle
    x-x',\bl{\nabla h}(x)-\bl{\nabla h}(x')\big\rangle\\
    &= {\textstyle\beta \left(\frac{1}{\tau}-\sigma\|L\|^2\right)^{-1}}\big\langle
    z-z',\bl{C}z-\bl{C}z'\big\rangle\\
    &= {\textstyle\beta \left(\frac{1}{\tau}-\sigma\|L\|^2\right)^{-1}}\big\langle
    z-z',P^{-1}\bl{C}z-P^{-1}\bl{C}z'\big\rangle_P.\label{eqcoer2}
\end{align}
Thus, $P^{-1}\bl{C}$ is $\chi$-cocoercive in $\mathcal{Z}_P$ with $\chi=\frac{1}{\beta}\big(\frac{1}{\tau}-\sigma\|L\|^2\big)$. Moreover, $\chi>1/2$ if and only if $\tau(\sigma\|L\|^2+\beta/2)<1$. Finally, $\delta = 2-1/(2\chi)$.
\end{proof}\smallskip

We can observe that if $\h=0$, the Condat--V\~u iteration reverts to the Chambolle--Pock iteration, so the former can be viewed as a generalization of the latter. 
Accordingly, if we set $\beta=0$ in \cref{thm61}, we recover \cref{thm41}.

The Condat--V\~u algorithm and the Loris--Verhoeven algorithms are both primal-dual forward-backward algorithms, but they are different. When $\ff=0$, larger values of $\tau$ and $\sigma$ are allowed in the latter than in the former; this may be beneficial to the convergence speed in practice. 

We can mention a different proximal splitting algorithm to solve \eqref{eq18} and \eqref{eq30}, proposed by Combettes and Pesquet~\cite{com12} before the Condat--V\~u algorithm. It is based on the forward-backward-forward splitting technique~\cite{tse00}\cite[section 26.6]{bau17}.\\

For the Condat--V\~u  algorithm too, let us focus on the case where $\h$ is quadratic; that is,
\begin{equation}
\h:x\mapsto {\textstyle\frac{1}{2}}\langle x,\bl{Q}x\rangle + \langle x,c\rangle 
\end{equation}
for some self-adjoint, positive, nonzero, bounded linear operator $\bl{Q}$ on $\mathcal{X}$ and some $c\in\mathcal{X}$. We have $\beta=\|\bl{Q}\|$. We can rewrite the primal-dual inclusion \eqref{var_inclusioncv}, which characterizes the Condat--V\~u iteration  \eqref{algocv}, as
\begin{equation}\label{var_inclusioncv3}
	\left(\! \begin{array}{c}
	0  \\ 
	0  \\
	\end{array}\!\right) \!
	\in\! 
	\underbrace{\left(\!\begin{array}{c}
(\dg{\partial f} +\bl{\nabla h})(x^{(i+\frac{1}{2})})+	L^* u^{(i+\frac{1}{2})} \\ 
	-L x^{(i+\frac{1}{2})}+(\dr{\partial g})^{-1} u^{(i+\frac{1}{2})} 
	\end{array}\!\right)}_{Mz^{(i+\frac{1}{2})}}
	+
	\underbrace{\left(\!\begin{array}{cc} 
	\frac{1}{\tau}\mathrm{Id}-\bl{Q} & -L^* \\ 
	-L & \frac{1}{\sigma}\mathrm{Id} 
	\end{array}\!\right)}_{P}
	\underbrace{\left(\!\begin{array}{c}
	x^{(i+\frac{1}{2})}-x^{(i)}  \\ 
	u^{(i+\frac{1}{2})}-u^{(i)}
	\end{array}\!\right)}_{z^{(i+\frac{1}{2})}-z^{(i)}}.
\end{equation}
Similarly, we can rewrite the primal-dual inclusion \eqref{var_inclusioncv2}, which characterizes the second form of the Condat--V\~u iteration  \eqref{algocv2}, as \eqref{var_inclusioncv3}, with $-L$ replaced by $L$ in $P$.

In both cases, using the properties of the Schur complement, $P$ is strongly positive if and only if 
\begin{equation}
    \tau \|\bl{Q}+\sigma L^*L\|<1
\end{equation}
(which implies that $\tau<1/\beta$).
A sufficient condition for this inequality to hold is $\tau(\sigma\|L\|^2+\beta)<1$. However, in some applications,  $\|\bl{Q}+\sigma L^*L\|$ may be smaller than $\sigma \|L\|^2+\beta$, so that larger stepsizes $\tau$ and $\sigma$ may be used when $\h$ is quadratic, leading to faster convergence.

Thus, when $\h$ is quadratic, the Condat--V\~u iteration can be viewed as a preconditioned primal-dual proximal point algorithm, just like the Chambolle--Pock iteration. Accordingly, we can apply  \cref{thm23} and obtain convergence under the condition $\tau \|\bl{Q}+\sigma L^*L\|<1$. 
Instead, let us apply the stronger \cref{thm26}
to allow $\tau \|\bl{Q}+\sigma L^*L\|=1$. For this, let us go back to the forward-backward analysis in 
\eqref{var_inclusioncv} and \eqref{var_inclusioncv2}. We suppose that $\tau\sigma\|L\|^2<1$. Then we can 
 strengthen the analysis in \eqref{eqcoer1}--\eqref{eqcoer2}: $P^{-1}\bl{C}$ is $\chi$-cocoercive in $\mathcal{Z}_P$ with $\chi=\|\left(\frac{1}{\tau}\mathrm{Id}-\sigma L^*L\right)^{-1}\bl{Q}\|^{-1}$. Then $\chi\geq 1$ 
if $\tau \|\bl{Q}+\sigma L^*L\|\leq 1$. Hence, we have:

\begin{theorem}[Condat--V\~u algorithm \eqref{algocv} or \eqref{algocv2}, quadratic case]\label{thm62}
Suppose that $\h$ is quadratic. Let $x^{(0)}\in\mathcal{X}$ and 
$u^{(0)}\in\mathcal{U}$. Let $\tau >0$ and $\sigma>0$  be such that $\tau \sigma \|L\|^2<1$ and $\tau \|\bl{Q}+\sigma L^*L\|\leq 1$.
Let $(\rho^{(i)})_{i\in\mathbb{N}}$ be a sequence in $[0,2]$ such that $\sum_{i\in \mathbb{N}} \rho^{(i)} (2-\rho^{(i)})=+\infty$. Then the sequences $(x^{(i)})_{i\in \mathbb{N}}$ and $(u^{(i)})_{i\in \mathbb{N}}$, defined  either by the iteration \eqref{algocv} or by the iteration \eqref{algocv2}, converge weakly to a solution of \eqref{eq18} and a solution of \eqref{eq30}, respectively.\end{theorem}

We note that $\tau (\beta+ \sigma \|L\|^2)\leq 1$ implies $\tau \sigma \|L\|^2<1$ and $\tau \|\bl{Q}+\sigma L^*L\|\leq 1$.

\subsection{An Extended Generalized Chambolle--Pock algorithm}\label{sec61}

Now that we have seen the different ways to construct primal-dual algorithms based on the forward-backward or proximal point algorithms, we can imagine further variations. For instance, we can extend the Generalized Chambolle--Pock algorithm in \cref{sec5} to deal with additional smooth terms on $u$ and $v$ in the dual problem. Let us look at one particular case of this extension: we add a smooth quadratic term on $v$. The problem is:
\begin{equation}
    \minimize_{(u,v)\in\mathcal{U}\times\mathcal{V}}\, \dg{f^*(v)} + \bl{\textstyle\frac{1}{2}\langle v,Q v\rangle} + \langle v,t\rangle + \dr{g^*(u)}\quad\mbox{s.t.}\quad K^*v+L^*u+c=0\label{eq2922}
\end{equation}
for some self-adjoint, positive, nonzero, bounded linear operator $\bl{Q}$ on $\mathcal{V}$ and some element $t\in\mathcal{V}$.  
The problem \eqref{eq2922} has applications, for instance, to solving inverse problems in imaging regularized with the total variation as defined by the first author~\cite{con17}. Again, we design a primal-dual  forward-backward algorithm, whose iteration satisfies 
\begin{align}\label{algogcp2}
\left(\!\begin{array}{l}
0\\
0\\
0
\end{array}\!\right)&\in
\underbrace{\left(\begin{array}{l}
 c+L^*u^{(i+\frac{1}{2})}+K^*v^{(i+\frac{1}{2})}\\
-Lx^{(i+\frac{1}{2})}+ \dr{\partial  g^*}(u^{(i+\frac{1}{2})})\\
-Kx^{(i+\frac{1}{2})}+ \dg{\partial  f^*}(v^{(i+\frac{1}{2})})+t\\
\end{array}\!\right)}_{Mz^{(i+\frac{1}{2})}}+
\underbrace{\left(\!\!\begin{array}{c}
0\\
0\\
\bl{Q}v^{(i)}
\end{array}\!\!\!\right)}_{\bl{C}z^{(i)}}\\
&\qquad\quad+
\underbrace{\left(\!\begin{array}{ccc}
\frac{1}{\tau}\mathrm{Id}&-L^*&0\\
-L&\frac{1}{\sigma}\mathrm{Id}&0\\
0& 0&\frac{1}{\eta}\mathrm{Id}-\tau KK^*
\end{array}\right)}_{P}
\underbrace{\left(\begin{array}{c}
x^{(i+\frac{1}{2})}-x^{(i)}\\
u^{(i+\frac{1}{2})}-u^{(i)}\\
v^{(i+\frac{1}{2})}-v^{(i)}
\end{array}\!\right)}_{z^{(i+\frac{1}{2})}-z^{(i)}}.\notag
\end{align}
Accordingly, the algorithm is:
\ceq{align}{\label{eqhybalg5}
&\,\mbox{\textbf{Primal-dual iteration for \eqref{eq2922}}: for }i=0,1,\ldots,\notag\\[-1mm]
    &\left\lfloor
\begin{array}{l}
v^{(i+\frac{1}{2})}
= \mathrm{prox}_{\eta \dg{f^*}}\Big(v^{(i)}+ \eta K 
\big( x^{(i)}-\tau (L^* u^{(i)}+K^* v^{(i)}+c )\big)\\
\qquad\qquad\qquad\qquad{}-\eta (\bl{Q}v^{(i)}+t)\Big)\\
x^{(i+\frac{1}{2})}
= x^{(i)} -\tau (L^* u^{(i)}+K^* v^{(i+\frac{1}{2})}+ c)\\
u^{(i+\frac{1}{2})}
= \mathrm{prox}_{\sigma \dr{g^*}}\big(u^{(i)}+ \sigma L 
(2x^{(i+\frac{1}{2})}-x^{(i)})\big)\\
x^{(i+1)}=x^{(i)}+\rho^{(i)} (x^{(i+\frac{1}{2})}-x^{(i)})\\
u^{(i+1)}=u^{(i)}+\rho^{(i)} (u^{(i+\frac{1}{2})}-u^{(i)})\\
v^{(i+1)}=v^{(i)}+\rho^{(i)} (v^{(i+\frac{1}{2})}-v^{(i)}).
\end{array}\right.}
The analysis is the same as for the Condat--V\~u algorithm: suppose that $\tau\sigma\|L\|^2<1$ and $\tau\eta\|K\|^2< 1$. Then $P$ is strongly positive and $P^{-1}\bl{C}$ is $\chi$-cocoercive in $\mathcal{Z}_P$ with $\chi=\|(\frac{1}{\eta}\mathrm{Id}-\tau KK^*)^{-1}\bl{Q}\|^{-1}$. Moreover, $\chi\geq 1$ if $\eta \|\tau KK^*+\bl{Q}\|\leq 1$. Hence, we can apply \cref{thm26} with $\gamma=1$ and we obtain:

\begin{theorem}[primal-dual algorithm \eqref{eqhybalg5}]\label{thm63}
Let $x^{(0)}\in\mathcal{X}$, $u^{(0)}\in\mathcal{U}$, $v^{(0)}\in\mathcal{V}$. 
Let $\tau >0$, $\sigma>0$, $\eta >0$ be such that $\tau\sigma\|L\|^2<1$, $\tau\eta\|K\|^2<1$, and $\eta\|\tau KK^*+\bl{Q}\|\leq 1$. Let $(\rho^{(i)})_{i\in\mathbb{N}}$ be a sequence in $[0,2]$ such that $\sum_{i\in \mathbb{N}} \rho^{(i)} (2-\rho^{(i)})=+\infty$.  Then  $(u^{(i)},v^{(i)})_{i\in \mathbb{N}}$ defined by the iteration \eqref{eqhybalg5} converges weakly to a solution of  \eqref{eq2922}.\end{theorem}

\section{The Primal-Dual Three-Operator Splitting (PD3O) and the Davis--Yin Algorithms}\label{sec7}

We consider the same problem as in the previous section, which involves three functions and a linear operator:
\begin{equation}
    \minimize_{x\in\mathcal{X}}\, \dg{f(x)} + \dr{g(Lx)} + \bl{h(x)}.\label{eq182}
\end{equation}
The dual problem is, again:
\begin{equation}
    \minimize_{u\in\mathcal{U}}\, (\ff+\h)^*(-L^*u) + \dr{g^*(u)}.\label{eq302}
\end{equation}

Let $\tau>0$ and $\sigma>0$, let $s^{(0)}\in\mathcal{X}$ and 
$u^{(0)}\in\mathcal{U}$, and let $(\rho^{(i)})_{i\in\mathbb{N}}$ be a sequence of relaxation parameters. 
The Primal-Dual Three-Operator Splitting (PD3O) algorithm proposed by Yan~\cite{yan18} is:
\ceq{align}{\label{eqpd3o1}
&\,\mbox{\textbf{PD3O iteration for \eqref{eq182}}: for }i=0,1,\ldots,\notag\\[-1mm]
    &\left\lfloor
    \begin{array}{l}
x^{(i+\frac{1}{2})}=\mathrm{prox}_{\tau \ff}(s^{(i)})\\
u^{(i+\frac{1}{2})}=\mathrm{prox}_{\sigma \dr{g^*}}\Big(u^{(i)}+\sigma L\big(
2x^{(i+\frac{1}{2})}-s^{(i)}-\tau \bl{\nabla h}(x^{(i+\frac{1}{2})})
-\tau L^*u^{(i)}
\big)\Big)\\
s^{(i+1)}=s^{(i)}+\rho^{(i)} \big(x^{(i+\frac{1}{2})} -s^{(i)}-\tau \bl{\nabla h}(x^{(i+\frac{1}{2})})-\tau L^*u^{(i+\frac{1}{2})}\big)\\
u^{(i+1)}=u^{(i)}+\rho^{(i)} (u^{(i+\frac{1}{2})}-u^{(i)})
\end{array}\right.}
(see \eqref{eqpd3o1p} for a more compact form with one call to $\bl{\nabla h}$ and $L^*$ per iteration).

If $\h$ is quadratic ($\bl{\nabla h}$ is affine) or if $\rho^{(i)}=1$, we have 
$s^{(i)}=
x^{(i)}  -\tau \bl{\nabla h}(x^{(i)})-\tau L^*u^{(i)}$ and
the PD3O iteration can be rewritten as
\begin{align}
    &\left\lfloor
    \begin{array}{l}
x^{(i+\frac{1}{2})}=\mathrm{prox}_{\tau \ff}\big(x^{(i)}-\tau \bl{\nabla h}(x^{(i)})-\tau L^*u^{(i)} \big)\\
u^{(i+\frac{1}{2})}=\mathrm{prox}_{\sigma \dr{g^*}}\Big(u^{(i)}
+\sigma L\big(
2x^{(i+\frac{1}{2})}-x^{(i)}+\tau\bl{\nabla h}(x^{(i)})-\tau \bl{\nabla h}(x^{(i+\frac{1}{2})})
\big)\Big)\\
x^{(i+1)}=x^{(i)}+\rho^{(i)} (x^{(i+\frac{1}{2})}-x^{(i)})\\
u^{(i+1)}=u^{(i)}+\rho^{(i)} (u^{(i+\frac{1}{2})}-u^{(i)}).
\end{array}\right. 
\end{align}
Thus, we recover the Chambolle--Pock algorithm when $\h=0$. When $\ff=0$, the PD3O algorithm reverts to the Loris--Verhoeven algorithm (with $s^{(i)}=x^{(i+\frac{1}{2})}$ in the former playing the role of $x^{(i)}$ in the latter).

Unlike the Condat--V\~u algorithm, the PD3O algorithm is not a primal-dual forward-backward algorithm; if $\sigma\tau\|L\|^2< 1$, it can be viewed as a primal-dual Davis--Yin algorithm~\cite{yan18,sal20} in the same way that the Loris--Verhoeven is a primal-dual forward-backward algorithm. Accordingly, we have:

\begin{theorem}[PD3O algorithm \eqref{eqpd3o1}]\label{thm71}
Let $s^{(0)}\in\mathcal{X}$ and 
$u^{(0)}\in\mathcal{U}$. Let $\tau \in (0, 2/\beta)$ and $\sigma>0$ be such that $\sigma\tau\|L\|^2<1$. Set $\delta=2-\tau\beta/2$. Let $(\rho^{(i)})_{i\in\mathbb{N}}$ be a sequence in $[0,\delta]$ such that $\sum_{i\in \mathbb{N}} \rho^{(i)} (\delta-\rho^{(i)})=+\infty$.
Then the sequences $(x^{(i+\frac{1}{2})})_{i\in \mathbb{N}}$ and $(u^{(i)})_{i\in \mathbb{N}}$ defined by the iteration \eqref{eqpd3o1} converge weakly to a solution of \eqref{eq182} and a solution of \eqref{eq302}, respectively.\end{theorem}

O'Connor and Vandenberghe showed that the PD3O algorithm can be obtained from the Davis--Yin algorithm, shown below, by a reformulation in an augmented space of the same type as the one shown for the Generalized Chambolle--Pock algorithm in \cref{sec5}. Their analysis makes it possible to have $\sigma\tau\|L\|^2=1$~\cite{oco18}. So, from the results on the Davis--Yin algorithm stated in \cref{thm73}, we can extend \cref{thm71} to $\sigma\tau\|L\|^2\leq 1$ as follows:

\begin{theorem}[PD3O algorithm \eqref{eqpd3o1}]\label{thm72}
Suppose that $\mathcal{X}$ and $\mathcal{U}$ are of finite dimension. 
Let $s^{(0)}\in\mathcal{X}$ and 
$u^{(0)}\in\mathcal{U}$.
Let $\tau \in (0, 2/\beta)$ and  $\sigma>0$ be such that $\sigma\tau\|L\|^2\leq 1$. Set $\delta=2-\tau\beta/2$. Let $(\rho^{(i)})_{i\in\mathbb{N}}$ be a sequence in $[0,\delta]$ such that $\sum_{i\in \mathbb{N}} \rho^{(i)} (\delta-\rho^{(i)})=+\infty$.
Then the sequences $(x^{(i+\frac{1}{2})})_{i\in \mathbb{N}}$ and $(u^{(i)})_{i\in \mathbb{N}}$ defined by the iteration \eqref{eqpd3o1} converge to a solution of \eqref{eq182} and a solution of \eqref{eq302}, respectively.\end{theorem}

Note that the finite dimension assumption is not necessary in the proof of this result \cite{oco18}, so it could be removed.

\br{Thus, the conditions for convergence of the PD3O algorithm are the same as for the Loris--Verhoeven algorithm, so that we get the additional function $\ff$ for free.}
\medskip

When $\mathcal{X}=\mathcal{U}$ and $L=\mathrm{Id}$, the problem \eqref{eq182} becomes:
\begin{equation}
    \minimize_{x\in\mathcal{X}}\, \dg{f(x)} + \dr{g(x)} + \bl{h(x)},\label{eq183}
\end{equation}
and if we set $\sigma=1/\tau$, the PD3O algorithm reverts to the Davis--Yin algorithm~\cite{dav17}:
\begin{align}
    &\left\lfloor
    \begin{array}{l}
x^{(i+\frac{1}{2})}=\mathrm{prox}_{\tau \ff}(s^{(i)})\\
u^{(i+\frac{1}{2})}=\mathrm{prox}_{ \dr{g^*}/\tau}\Big( \big(
2x^{(i+\frac{1}{2})}-s^{(i)}-\tau \bl{\nabla h}(x^{(i+\frac{1}{2})})\big)/\tau
\Big)\\
\hphantom{u^{(i+\frac{1}{2})}}=(\mathrm{Id}-\mathrm{prox}_{\tau \g}) \big(
2x^{(i+\frac{1}{2})}-s^{(i)}-\tau\bl{ \nabla h}(x^{(i+\frac{1}{2})})\big)/\tau\\
s^{(i+1)}=s^{(i)}+\rho^{(i)} \big(x^{(i+\frac{1}{2})}-s^{(i)} -\tau\bl{\nabla h}(x^{(i+\frac{1}{2})})-\tau u^{(i+\frac{1}{2})}\big)\\
\hphantom{s^{(i+1)}}=s^{(i)}+\rho^{(i)} \Big(\mathrm{prox}_{\tau \g}\big(
2x^{(i+\frac{1}{2})}-s^{(i)}-\tau \bl{\nabla h}(x^{(i+\frac{1}{2})})\big)
-x^{(i+\frac{1}{2})}\Big),
\end{array}\right.
\end{align}
which can be simplified as:
\ceq{align}{\label{eqdy}
&\,\mbox{\textbf{Davis--Yin iteration for \eqref{eq183}}: for }i=0,1,\ldots,\notag\\[-1mm]
    &\left\lfloor
    \begin{array}{l}
x^{(i+\frac{1}{2})}=\mathrm{prox}_{\tau \ff}(s^{(i)})\\
s^{(i+1)} =s^{(i)}+\rho^{(i)}\Big(
\mathrm{prox}_{\tau \g}\big(
2x^{(i+\frac{1}{2})}-s^{(i)}-\tau \bl{\nabla h}(x^{(i+\frac{1}{2})})\big)
-x^{(i+\frac{1}{2})}\Big).
\end{array}\right.}

If $\ff=0$, the Davis--Yin algorithm reverts to the forward-backward algorithm.  If $\h=0$, it reverts to  the Douglas--Rachford algorithm. 
The Davis--Yin algorithm generalizes the forward-Douglas--Rachford algorithm and the Generalized forward-backward algorithm~\cite{rag13,bri15,rag19}. 

Its main convergence result is:~\cite{dav17,ara21}:

\begin{theorem}[Davis--Yin algorithm \eqref{eqdy}]\label{thm73}
Let $s^{(0)}\in\mathcal{X}$. Let $\tau \in (0, 2/\beta)$. Set $\delta=2-\tau\beta/2$. Let $(\rho^{(i)})_{i\in\mathbb{N}}$ be a sequence in $[0,\delta]$ such that $\sum_{i\in \mathbb{N}} \rho^{(i)} (\delta-\rho^{(i)})=+\infty$. 
Then the sequence $(x^{(i+\frac{1}{2})})_{i\in \mathbb{N}}$ defined by the iteration \eqref{eqdy} converges weakly to a solution of \eqref{eq183}.\end{theorem}

\subsection{A new algorithm for the quadratic case}\label{secnaq}

Let us continue with the minimization of $\ff+\g \circ L + \h$ in the quadratic case where $\h(x)= {\textstyle\frac{1}{2}}\langle x,\bl{Q}x\rangle + \langle x,c\rangle$
for some self-adjoint, positive, nonzero, bounded linear operator $\bl{Q}$ on $\mathcal{X}$ and some $c\in\mathcal{X}$. Set $\beta=\|\bl{Q}\|$. There seems to be no way to interpret the PD3O algorithm in a different way, to enlarge the relaxation range to $(0,2)$. However, we can propose a new algorithm that has the widest range of parameters seen so far.

Let $\tau\in(0,2/\beta)$ and $\sigma>0$ be such that $\sigma\tau\|L\|^2<1$.
We introduce the following self-adjoint, strongly positive, bounded linear operator $P$ in $\mathcal{Z}$:
\begin{equation}
P=	\left( \begin{array}{cc} 
	\frac{1}{\tau}\mathrm{Id}-\frac{1}{2}\bl{Q} & 0 \\ 
	0 & \frac{1}{\sigma}\mathrm{Id}-\tau LL^* \\
	\end{array} \right).
\end{equation}
Note that unlike what we have done so far, only one half of $\bl{Q}$ appears in $P$. 
We define the real Hilbert space $\mathcal{Z}_P$ as $\mathcal{Z}$ endowed with the inner product $\langle\cdot,\cdot\rangle_P:(z,z')\mapsto \langle z,Pz '\rangle$.
We will solve the following monotone inclusion in $\mathcal{Z}_P$:
\begin{equation}\label{variational_prim_dual_prob4}
   \left( \begin{array}{c}
   0\\0
   \end{array}\right)\in
    \underbrace{P^{-1}\left( \begin{array}{c}
     \dg{\partial f}(x)+\frac{1}{2}\bl{Q}x+c \\
   0
   \end{array}\right)}_{M z}+
    \underbrace{P^{-1} \left(\begin{array}{c}
    \frac{1}{2}\bl{Q}x+L^*u \\
   -Lx + (\dr{\partial g})^{-1}u
   \end{array}\right)}_{N z}.
\end{equation}
Let $s^{(0)}\in\mathcal{Z}$. The operators $M$ and $N$ are maximally monotone in 	$\mathcal{Z}_P$, so we will use the Douglas--Rachford algorithm to solve \eqref{variational_prim_dual_prob4} (with stepsize $\gamma=1)$; it can be written as:
\begin{align}\label{algodrc_1}
&\,\mbox{\textbf{Douglas--Rachford iteration for \eqref{variational_prim_dual_prob4}}: for }i=0,1,\ldots,\notag\\[-1mm]
    &\left\lfloor
    \begin{array}{l}
   z^{(i+\frac{1}{2})}=J_{ M}(s^{(i)})\\
   w^{(i+\frac{1}{2})}=J_{ N}(2z^{(i+\frac{1}{2})} - s^{(i)})\\
s^{(i+1)}=s^{(i)}+\rho^{(i)}\big(w^{(i+\frac{1}{2})} -z^{(i+\frac{1}{2})}\big).
\end{array}\right.\end{align}
We have $J_M:(x,u)\mapsto \big(\mathrm{prox}_{\tau \ff}(x-\frac{\tau}{2}\bl{Q}x-\tau c),u\big)$, and the resolvent $J_N$ amounts to one iteration of the Loris--Verhoeven algorithm \eqref{var_inclusion2}, so it maps $(x,u)$ to $(x',u')$ with:
\begin{align}\label{algorithm_reslv}
 &\left\lfloor
    \begin{array}{l}
    u'
= \mathrm{prox}_{\sigma \dr{g^*}}\big(u+\sigma L (x - \frac{\tau}{2} \bl{Q}x-\tau L^* u)\big)\\
x'=x-\frac{\tau}{2} \bl{Q}x -\tau L^* u'.
\end{array}\right. 
\end{align}
Substituting these two expressions of the resolvents into \eqref{algodrc_1}, we obtain the new primal-dual Douglas--Rachford iteration:
\ceq{align}{\label{algodrc_2}
&\,\mbox{\textbf{Primal-dual Douglas--Rachford iteration, form I, for \eqref{eq182}}:}\notag\\ 
&\,\mbox{for }i=0,1,\ldots,\notag\\[-1mm]
    &\left\lfloor
    \begin{array}{l}
   x^{(i+\frac{1}{2})}=\mathrm{prox}_{\tau \ff}\big(s^{(i)}-\frac{\tau}{2}\bl{Q}s^{(i)}-\tau c\big)\\
   u^{(i+\frac{1}{2})}
= \mathrm{prox}_{\sigma \dr{g^*}}\big(u^{(i)}+\sigma L (2x^{(i+\frac{1}{2})}-s^{(i)} - \frac{\tau}{2} \bl{Q}(2x^{(i+\frac{1}{2})}-s^{(i)})\\
\qquad\qquad\qquad\qquad{}-\tau L^* u^{(i)})\big)\\
 s^{(i+\frac{1}{2})}=x^{(i+\frac{1}{2})}-\frac{\tau}{2} \bl{Q}(2x^{(i+\frac{1}{2})}-s^{(i)}) -\tau L^* u^{(i+\frac{1}{2})}\\
s^{(i+1)}=s^{(i)}+\rho^{(i)}( s^{(i+\frac{1}{2})}-s^{(i)})\\
u^{(i+1)}=u^{(i)}+\rho^{(i)}(u^{(i+\frac{1}{2})} -u^{(i)}).
\end{array}\right.}

If we switch the roles of $M$ and $N$ in \eqref{algodrc_1}, we obtain the iteration:
\ceq{align}{\label{algodrc_3}
&\,\mbox{\textbf{Primal-dual Douglas--Rachford iteration, form II, for \eqref{eq182}}:}\notag\\
&\,\mbox{for }i\!=\!0,1,\ldots,\notag\\[-1mm]
    &\left\lfloor
    \begin{array}{l}
  u^{(i+\frac{1}{2})}
= \mathrm{prox}_{\sigma \dr{g^*}}\big(u^{(i)}+\sigma L (s^{(i)} - \frac{\tau}{2} \bl{Q}s^{(i)}-\tau L^* u^{(i)})\big)\\
y^{(i+\frac{1}{2})}=s^{(i)}-\frac{\tau}{2} \bl{Q}s^{(i)} -\tau L^* u^{(i+\frac{1}{2})}\\
x^{(i+\frac{1}{2})}=\mathrm{prox}_{\tau \ff}\big(2y^{(i+\frac{1}{2})}-s^{(i)}-\frac{\tau}{2}\bl{Q}(2y^{(i+\frac{1}{2})}-s^{(i)})-\tau c\big)\\
s^{(i+1)}=s^{(i)}+\rho^{(i)}( x^{(i+\frac{1}{2})}-y^{(i+\frac{1}{2})})\\
u^{(i+1)}=u^{(i)}+\rho^{(i)}(u^{(i+\frac{1}{2})}-u^{(i)}).
\end{array}\right.}

From the convergence of the Douglas--Rachford algorithm \eqref{algodrc_1}~\cite[Theorem 26.11]{bau17}, we obtain:

\begin{theorem}[Primal-dual Douglas--Rachford algorithm \eqref{algodrc_2} or \eqref{algodrc_3}, quadratic case]\label{thm74}
Let $s^{(0)}\in\mathcal{X}$ and $u^{(0)}\in\mathcal{U}$.
Let $\tau \in (0, 2/\beta)$ and  $\sigma>0$ be such that $\sigma\tau\|L\|^2<1$.  Let $(\rho^{(i)})_{i\in\mathbb{N}}$ be a sequence in $[0,2]$ such that $\sum_{i\in \mathbb{N}} \rho^{(i)} (2-\rho^{(i)})=+\infty$.
Then the sequences $(x^{(i+\frac{1}{2})})_{i\in \mathbb{N}}$ and $(u^{(i)})_{i\in \mathbb{N}}$, defined either by the iteration \eqref{algodrc_2} or by the iteration \eqref{algodrc_3}, converge weakly to a solution of \eqref{eq182} and a solution of \eqref{eq302}, respectively.\end{theorem}

In comparison to the Condat--V\~u and the PD3O algorithms, with the new algorithm, $\tau$ is allowed to be in $(0,2/\beta)$, the condition on $\tau$ and $\sigma$ does not depend on $\bl{Q}$, and  the relaxation range is $(0,2)$. Like them,  the new algorithm reverts to the Chambolle--Pock algorithm when $\bl{Q}=0$. The price to pay is two evaluations of $\bl{Q}$ per iteration.

\section{Parallel Versions of the Algorithms}\label{sec8}

In this section, we generalize some of the algorithms seen so far, to make it possible to solve optimization problems involving an arbitrary number of functions $\dr{g_m}$, $m=1,\ldots,M$, instead of only one function $\g$. The derived algorithms are parallel, in the sense that the proximity operators of the functions $\dr{g_m}$ or $\dr{g_m^*}$ are independent of each other and can be applied in parallel. We present two techniques for this generalization, which consist in reformulating the problems in product spaces, an idea which has been known for a long time \cite{pie84}. 
When $M=1$, we recover the initial algorithms. However, we should keep in mind that the variables, which are updated in parallel with respect to antagonistic functions, are then essentially averaged to form the update of the variable $x$ of interest. This means that the larger $M$, the slower the convergence to a solution. For instance, it should be 
 faster to use the Douglas--Rachford algorithm to minimize $\ff +\g$ than to apply any parallel algorithm to minimize $\dr{g_1}+\dr{g_2}$. Similarly, the Generalized Chambolle--Pock algorithm in  \cref{sec5} to minimize $\dg{f\circ K}+\dr{g\circ L}$ should be faster than the parallel version of the Proximal Method of Multipliers in  \cref{secpmm} to minimize $\dr{g_1\circ L_1}+\dr{g_2\circ L_2}$.\medskip

\noindent\textbf{Technique 1}: Let $M\geq 1$ be an integer. Let $\mathcal{X}$ be a real Hilbert space. We consider the minimization of $\dg{f(x)}+\dr{\sum_{m=1}^M g_m(x)} + \bl{h(x)}$ over $x\in\mathcal{X}$ for some $\ff\in\Gamma_0(\mathcal{X})$, $\dr{g_m}\in \Gamma_0(\mathcal{X})$, and differentiable and convex function $\h$ with $\beta$-Lipschitz continuous gradient for some $\beta\geq 0$. 
Let $(\omega_m)_{m=1}^M$ be a sequence of positive weights whose sum is $1$. We introduce the real Hilbert space $\boldsymbol{\mathcal{X}}=\mathcal{X}\times\cdots\times\mathcal{X}$ ($M$ times), endowed with the inner product $\langle\cdot\,,\cdot\rangle_{\boldsymbol{\mathcal{X}}}:(\boldsymbol{x},\boldsymbol{x'})\mapsto \sum_{m=1}^M \omega_m \langle x_m,x'_m\rangle$. We introduce the function $\imath_=: \boldsymbol{x}=(x_m)_{m=1}^M\in \boldsymbol{\mathcal{X}}\mapsto (0$ if $x_1=\cdots= x_M,$ $+\infty$ otherwise). Then we can rewrite the function to minimize as $\dg{\boldsymbol{f}(\boldsymbol{x})}+\dr{\boldsymbol{g}(\boldsymbol{x})}+\bl{\boldsymbol{h}(\boldsymbol{x})}$ with $\dg{\boldsymbol{f}}:\boldsymbol{x}\in\boldsymbol{\mathcal{X}}\mapsto \ff(x_1)+\imath_=(\boldsymbol{x})$, $\dr{\boldsymbol{g}}:\boldsymbol{x}\in\boldsymbol{\mathcal{X}}\mapsto \sum_{m=1}^M \dr{g_m}(x_m)$, and $\bl{\boldsymbol{h}}:\boldsymbol{x}\in\boldsymbol{\mathcal{X}}\mapsto\sum_{m=1}^M \omega_m\bl{ h}(x_m)$. Thus, we have the gradient $\bl{\nabla \boldsymbol{h}}:\boldsymbol{x}\mapsto\big(\bl{\nabla h}(x_1),\ldots,\bl{\nabla h}(x_M)\big)$, which is $\beta$-Lipschitz continuous, and the following proximity operators: for any $\gamma>0$, $\mathrm{prox}_{\gamma\dg{\boldsymbol{f}}}:\boldsymbol{x}\mapsto (x',\ldots,x')$ with
$x'=\mathrm{prox}_{\gamma\ff}\big(\sum_{m=1}^M \omega_m x_m\big)$, and $\mathrm{prox}_{\gamma\dr{\boldsymbol{g}}}:\boldsymbol{x}\mapsto\big(\mathrm{prox}_{(\gamma/\omega_1)\dr{g_1}}(x_1),\ldots,\mathrm{prox}_{(\gamma/\omega_M)\dr{g_M}}(x_M)\big)$. Note that one must be careful when changing the metric; for instance, we have $\dr{\partial \boldsymbol{g}}:\boldsymbol{x}\mapsto\big(\frac{1}{\omega_1}\dr{\partial g_1}(x_1),\ldots,\frac{1}{\omega_M}\dr{\partial g_M}(x_M)\big)$, 
$\dr{\boldsymbol{g^*}}:\boldsymbol{x}\in\boldsymbol{\mathcal{X}}\mapsto \sum_{m=1}^M \dr{g_m^*}(\omega_m x_m)$, and 

\noindent
$\mathrm{prox}_{\gamma\dr{\boldsymbol{g^*}}}:\boldsymbol{x}\mapsto\big(\frac{1}{\omega_1}\mathrm{prox}_{\gamma\omega_1\dr{g_1^*}}(\omega_1 x_1),\ldots,\frac{1}{\omega_M}\mathrm{prox}_{\gamma\omega_M\dr{g_M^*}}(\omega_M x_M)\big)$.\bigskip

\noindent \textbf{Technique 2}: Let $M\geq 1$ be an integer. Let $\mathcal{X}$ and $\mathcal{U}_m$ be real Hilbert spaces. We consider an optimization problem that involves a term $\dr{\sum_{m=1}^M g_m(L_m x)}$, for some $\dr{g_m}\in \Gamma_0(\mathcal{U}_m)$ and bounded linear operators $L_m:\mathcal{X}\rightarrow \mathcal{U}_m$. Let $(\omega_m)_{m=1}^M$ be a sequence of positive weights. We introduce the real Hilbert space $\boldsymbol{\mathcal{U}}=\mathcal{U}_1\times\cdots\times\mathcal{U}_M$ endowed with the inner product $\langle\cdot\,,\cdot\rangle_{\boldsymbol{\mathcal{U}}}:(\boldsymbol{u},\boldsymbol{u'})\mapsto \sum_{m=1}^M \omega_m \langle u_m,u'_m\rangle$. Then $\dr{\sum_{m=1}^M g_m \circ L_m}=\dr{\boldsymbol{g} \circ \boldsymbol{L}}$ with $\dr{\boldsymbol{g}}:\boldsymbol{u}\in\boldsymbol{\mathcal{U}}\mapsto \sum_{m=1}^M \dr{g_m}(u_m)$ and $\boldsymbol{L}:x\in\mathcal{X}\mapsto \big(L_1 x,\ldots,L_M x\big)\in\boldsymbol{\mathcal{U}}$. Then we have, for any $\gamma>0$, $\mathrm{prox}_{\gamma\dr{\boldsymbol{g}}}:\boldsymbol{u}\mapsto\big(\mathrm{prox}_{(\gamma/\omega_1)\dr{g_1}}(u_1),  
\ldots,\mathrm{prox}_{(\gamma/\omega_M)\dr{g_M}}(u_M)\big)$. One must be again careful when changing the metric; for instance, we have $\boldsymbol{L^*}:u\in\boldsymbol{\mathcal{U}}\mapsto \sum_{m=1}^M \omega_m L^*_m u_m$, $\|\boldsymbol{L^*}\boldsymbol{L}\|=\|\sum_{m=1}^M \omega_m L^*_mL_m\|$, and
$\mathrm{prox}_{\gamma\dr{\boldsymbol{g^*}}}:\boldsymbol{u}\mapsto\big(\frac{1}{\omega_1}\mathrm{prox}_{\gamma\omega_1\dr{g_1^*}}(\omega_1 u_1),\ldots,\frac{1}{\omega_M}\mathrm{prox}_{\gamma\omega_M\dr{g_M^*}}(\omega_M u_M)\big)$.\medskip

With both techniques, the dual variables must be multiplied by $\omega_m$, so that they converge to solutions to the original dual problems, not the problems expressed with $\boldsymbol{\mathcal{U}}$ and its distorted metric.

\subsection{The Douglas--Rachford algorithm}\label{secdrp}

Let $M\geq 1$ be an integer. Let $\mathcal{X}$ be a real Hilbert space. Let $\ff\in\Gamma_0(\mathcal{X})$ and $\dr{g_m}\in \Gamma_0(\mathcal{X})$, for $m=1,\ldots,M$. We consider the convex optimization problem:
\begin{equation}
\minimize_{x\in\mathcal{X}} \dg{f(x)}+\dr{\sum_{m=1}^M g_m(x)}.\label{eqpdrp}
\end{equation}

Let $(\omega_m)_{m=1}^M$ be a sequence of positive weights whose sum is $1$. Let $\tau>0$ and 
let $(\rho^{(i)})_{i\in\mathbb{N}}$ be a sequence of nonnegative reals. Let $s_m^{(0)}\in\mathcal{X}$ for $m=1,\ldots,M$. By applying Technique 1 to the Douglas--Rachford algorithm \eqref{algorithm_22}, we obtain the iteration:
\ceq{align}{\label{algorithm_22p}
&\,\mbox{\textbf{Douglas--Rachford iteration, form I, for \eqref{eqpdrp}}: for }i=0,1,\ldots,\notag\\[-1mm]
    &\left\lfloor
    \begin{array}{l}
   x^{(i+\frac{1}{2})}=\mathrm{prox}_{\tau \ff }\big(\sum_{m=1}^M \omega_m s^{(i)}_m\big)\\
   \mbox{For $m=1,\ldots,M$,}\\
   \big\lfloor\;\,
s_m^{(i+1)}=s_m^{(i)}+\rho^{(i)}\big( \mathrm{prox}_{ \frac{\tau}{\omega_m} \dr{g_m}}(2x^{(i+\frac{1}{2})} - s_m^{(i)})-x^{(i+\frac{1}{2})}\big).
\end{array}\right.}

If we switch the roles of $\boldsymbol{f}$ and $\boldsymbol{g}$, we obtain the iteration:
\ceq{align}{\label{algorithm_22p2}
&\,\mbox{\textbf{Douglas--Rachford iteration, form II, for \eqref{eqpdrp}}: for }i=0,1,\ldots,\notag\\[-1mm]
    &\left\lfloor
    \begin{array}{l}
    \mbox{For $m=1,\ldots,M$,}\\
   \big\lfloor\;\,x^{(i+\frac{1}{2})}_m= \mathrm{prox}_{ \frac{\tau}{\omega_m} \dr{g_m}}(s_m^{(i)})\\
x^{(i+\frac{1}{2})}=\mathrm{prox}_{\tau \ff }\big(\sum_{m=1}^M \omega_m(2x^{(i+\frac{1}{2})}_m-s^{(i)}_m)\big)\\
\mbox{For $m=1,\ldots,M$,}\\
 \big\lfloor\;\,s_m^{(i+1)}=s_m^{(i)}+\rho^{(i)}(x^{(i+\frac{1}{2})}-x^{(i+\frac{1}{2})}_m)
\end{array}\right.}

As an application of \cref{thm44}, we have:

\begin{theorem}[Douglas--Rachford algorithm \eqref{algorithm_22p}]\label{thm81}Let $s_m^{(0)}\in\mathcal{X}$, for $m=1,\ldots,M$. Let $\tau >0$, let $\omega_1>0,\ldots,\omega_M>0$ be reals, whose sum is $1$, and let $(\rho^{(i)})_{i\in\mathbb{N}}$ be a sequence in $[0,2]$ such that $\sum_{i\in \mathbb{N}} \rho^{(i)} (2-\rho^{(i)})=+\infty$.  Then the sequence $(x^{(i+\frac{1}{2})})_{i\in \mathbb{N}}$, defined either by the iteration \eqref{algorithm_22p} or by the iteration \eqref{algorithm_22p2}, converges weakly to a solution of \eqref{eqpdrp}.\end{theorem}

\subsection{The Chambolle--Pock algorithm}

Let $M\geq 1$ be an integer. Let $\mathcal{X}$ and $\mathcal{U}_m$, $m=1,\ldots,M$, be real Hilbert spaces. Let $\ff\in\Gamma_0(\mathcal{X})$ and let $\dr{g_m}\in \Gamma_0(\mathcal{U}_m)$, $m=1,\ldots,M$. Let $L_m:\mathcal{X}\rightarrow \mathcal{U}_m$, $m=1,\ldots,M$, be bounded linear operators. We consider the convex optimization problem:
\begin{equation}
\minimize_{x\in\mathcal{X}} \dg{f(x)}+\dr{\sum_{m=1}^M g_m(L_m x)}.\label{eq22p}
\end{equation}
The dual problem is:
\begin{equation}
    \minimize_{(u_1,\ldots,u_M)\in\mathcal{U}_1\times\cdots\times\mathcal{U}_M}\,\dg{f^*({\textstyle -\sum_{m=1}^M L_m^*u_m})} +  \dr{\sum_{m=1}^M g^*_m(u_m)}.\label{eq23p}
\end{equation}

Let $(\omega_m)_{m=1}^M$ be a sequence of positive weights.
Let $\tau>0$, $\sigma>0$, and 
let $(\rho^{(i)})_{i\in\mathbb{N}}$ be a sequence of nonnegative reals. Let $x^{(0)}\in\mathcal{X}$ and $u_m^{(0)}\in\mathcal{U}_m$ for $m=1,\ldots,M$. Applying Technique 2 to the Chambolle--Pock algorithm  \eqref{algoCP} or  \eqref{algoCP2} and setting $\sigma_m=\sigma\omega_m$, we obtain the iterations:
\ceq{align}{\label{algoCPp}
&\,\mbox{\textbf{Chambolle--Pock iteration, form I, for \eqref{eq22p} and \eqref{eq23p}}: for }i=0,1,\ldots,\notag\\[-1mm]
    &\left\lfloor
    \begin{array}{l}
    x^{(i+\frac{1}{2})}
= \mathrm{prox}_{\tau \ff}\big(x^{(i)}-\tau \sum_{m=1}^M L^*_m u_m^{(i)} \big)\\
x^{(i+1)}=x^{(i)}+\rho^{(i)} (x^{(i+\frac{1}{2})}-x^{(i)})\\
 \mbox{For $m=1,\ldots,M$,}\\
\left\lfloor\begin{array}{l}
u_m^{(i+\frac{1}{2})}
= \mathrm{prox}_{\sigma_m \dr{g^*_m}}\big(u_m^{(i)}+\sigma_m L_m(2x^{(i+\frac{1}{2})}-x^{(i)}) \big)\\
u_m^{(i+1)}=u_m^{(i)}+\rho^{(i)} (u_m^{(i+\frac{1}{2})}-u_m^{(i)}).
\end{array}\right.
\end{array}\right.}

\ceq{align}{\label{algoCP2p}
&\,\mbox{\textbf{Chambolle--Pock iteration, form II, for \eqref{eq22p} and \eqref{eq23p}}:}\notag\\
&\,\mbox{ for }i=0,1,\ldots,\notag\\[-1mm]
    &\left\lfloor
    \begin{array}{l}
 \mbox{For $m=1,\ldots,M$,}\\
\left\lfloor\begin{array}{l}
u_m^{(i+\frac{1}{2})} = \mathrm{prox}_{\sigma_m \dr{g_m^*}}\big(u_m^{(i)}+\sigma_m L_m x^{(i)} \big)\\
u_m^{(i+1)}=u_m^{(i)}+\rho^{(i)} (u_m^{(i+\frac{1}{2})}-u_m^{(i)})
\end{array}\right.\\
    x^{(i+\frac{1}{2})}
= \mathrm{prox}_{\tau \ff}\big(x^{(i)}-\tau  \sum_{m=1}^M L_m^* (2u_m^{(i+\frac{1}{2})}-u_m^{(i)}) \big)\\
x^{(i+1)}=x^{(i)}+\rho^{(i)} (x^{(i+\frac{1}{2})}-x^{(i)}).
\end{array}\right.}

We can remark that when both Technique 1 and Technique 2 are applicable, they give the same result: if $L_m=\mathrm{Id}$, the Chambolle--Pock algorithm \eqref{algoCPp} with $\sum_{m=1}^M \omega_m=1$ and $\sigma=1/\tau$ gives the Douglas--Rachford algorithm \eqref{algorithm_22p}.

As an application of \cref{thm42}, we have:

\begin{theorem}[Chambolle--Pock algorithm \eqref{algoCPp} or \eqref{algoCP2p}]\label{thm82}Let $x^{(0)}\in\mathcal{X}$, $u_1^{(0)}\in\mathcal{U}_1,\ldots,u_M^{(0)}\in\mathcal{U}_M$. Let $\tau >0$ and $\sigma_1>0,\ldots,\sigma_M>0$  be such that 
$\tau\|\sum_{m=1}^M \sigma_m L_m^*L_m\|\leq 1$. Let $(\rho^{(i)})_{i\in\mathbb{N}}$ be a sequence in $[0,2]$ such that $\sum_{i\in \mathbb{N}} \rho^{(i)} (2-\rho^{(i)})=+\infty$.  Then the sequences $(x^{(i)})_{i\in \mathbb{N}}$ and $\big((u_1^{(i)},\ldots,u_M^{(i)})\big)_{i\in \mathbb{N}}$, defined either by the iteration \eqref{algoCPp} or by the iteration \eqref{algoCP2p}, converge weakly to a solution of \eqref{eq22p} and a solution of \eqref{eq23p}, respectively.\end{theorem}

\subsection{The Loris--Verhoeven algorithm}

Let $M\geq 1$ be an integer. Let $\mathcal{X}$ and $\mathcal{U}_m$, $m=1,\ldots,M$, be real Hilbert spaces. Let $\dr{g_m}\in \Gamma_0(\mathcal{U}_m)$, $m=1,\ldots,M$, and let $\h:\mathcal{X}\rightarrow \mathbb{R}$ be a convex and 
differentiable function with  $\beta$-Lipschitz continuous gradient for some real $\beta> 0$. Let $L_m:\mathcal{X}\rightarrow \mathcal{U}_m$, $m=1,\ldots,M$, be bounded linear operators. We consider the convex optimization problem:
\begin{equation}
    \minimize_{x\in\mathcal{X}}\, \dr{\sum_{m=1}^M g_m(L_m x)} + \bl{h(x)}.\label{eq13p}
\end{equation}
The dual problem is:
\begin{equation}
     \minimize_{(u_1,\ldots,u_M)\in\mathcal{U}_1\times\cdots\times\mathcal{U}_M}\, \dr{\sum_{m=1}^M g^*_m(u_m)} + \bl{h^*({\textstyle -\sum_{m=1}^M L_m^*u_m})}.\label{eq21p}
\end{equation}

Let $(\omega_m)_{m=1}^M$ be a sequence of positive weights.
Let $\tau>0$, $\sigma>0$, and 
let $(\rho^{(i)})_{i\in\mathbb{N}}$ be a sequence of nonnegative reals. Let $x^{(0)}\in\mathcal{X}$ and $u_m^{(0)}\in\mathcal{U}_m$ for $m=1,\ldots,M$.
Applying Technique 2 to the Loris--Verhoeven algorithm \eqref{algorithm_1} and setting $\sigma_m=\sigma\omega_m$, we obtain the iteration:
\begin{align}\label{algorithm_1p0}
&\,\mbox{\textbf{Loris--Verhoeven iteration for \eqref{eq13p} and \eqref{eq21p}}: for }i=0,1,\ldots,\notag\\[-1mm]
    &\left\lfloor
    \begin{array}{l}
\mbox{For $m=1,\ldots,M$,}\\
\left\lfloor\begin{array}{l}
u_m^{(i+\frac{1}{2})}
= \mathrm{prox}_{\sigma_m \dr{g^*_m}}\Big(u_m^{(i)}+\sigma_m L_m \big(x^{(i)} - \tau\bl{\nabla h}(x^{(i)})-\tau \sum_{m=1}^M L^*_m u_m^{(i)}\big)\Big)\\
u_m^{(i+1)}=u_m^{(i)}+\rho^{(i)} (u_m^{(i+\frac{1}{2})}-u_m^{(i)})
\end{array}\right.\\
x^{(i+1)}=x^{(i)}-\rho^{(i)}\tau \big(
\bl{\nabla h}(x^{(i)}) + \sum_{m=1}^M L^*_m u_m^{(i+\frac{1}{2})}
\big).
\end{array}\right.\end{align}

If we introduce the variables $\tilde{s}^{(i)}=x^{(i)}+\tau\sum_{m=1}^M L^*_m u_m^{(i)}$ and $a^{(i)}=x^{(i)}-\tau \bl{\nabla h}(x^{(i)})-\tilde{s}^{(i)}=-\tau \bl{\nabla h}(x^{(i)})-\tau\sum_{m=1}^M L^*_m u_m^{(i)}$, for every $i\in\mathbb{N}$, we can rewrite the iteration as:
\ceq{align}{\label{algorithm_1p}
&\,\mbox{\textbf{Loris--Verhoeven iteration for \eqref{eq13p} and \eqref{eq21p}}: for }i=0,1,\ldots,\notag\\[-1mm]
    &\left\lfloor
    \begin{array}{l}
a^{(i)}=x^{(i)}-\tau \bl{\nabla h}(x^{(i)})-\tilde{s}^{(i)}\\
\tilde{s}^{(i+1)}=\tilde{s}^{(i)}+\rho^{(i)} a^{(i)}\\
 \mbox{For $m=1,\ldots,M$,}\\
\Big\lfloor\;
u_m^{(i+1)}=u_m^{(i)}+\rho^{(i)} \Big(\mathrm{prox}_{\sigma_m \dr{g^*_m}}\big(u_m^{(i)}+\sigma_m L_m(x^{(i)}+a^{(i)}) \big)-u_m^{(i)}\Big)\\
x^{(i+1)}=\tilde{s}^{(i+1)}-\tau \sum_{m=1}^M L_m^*u_m^{(i+1)}.
\end{array}\right.}

As an application of \cref{thm31,thm32}, we have:

\begin{theorem}[Loris--Verhoeven algorithm \eqref{algorithm_1p}]\label{thm83}Let $x^{(0)}\in\mathcal{X}$, $u_1^{(0)}\in\mathcal{U}_1,\ldots,
u_M^{(0)}\in\mathcal{U}_M$. Set $\tilde{s}^{(0)}=x^{(0)}+\tau\sum_{m=1}^M L^*_m u_m^{(0)}$. Let $\tau \in (0, 2/\beta)$ and $\sigma_1>0,\ldots,\sigma_M>0$ be such that $\tau\|\sum_{m=1}^M \sigma_m L_m^*L_m\|<1$. Set $\delta=2-\tau\beta/2$. Let $(\rho^{(i)})_{i\in\mathbb{N}}$ be a sequence in $[0,\delta]$ such that $\sum_{i\in \mathbb{N}} \rho^{(i)} (\delta-\rho^{(i)})=+\infty$. 
Then the sequences $(x^{(i)})_{i\in \mathbb{N}}$ and $\big((u_1^{(i)},\ldots,u_M^{(i)})\big)_{i\in \mathbb{N}}$ defined by the iteration \eqref{algorithm_1p} converge weakly to a solution of \eqref{eq13p} and a solution of \eqref{eq21p}, respectively.\end{theorem}

\begin{theorem}[Loris--Verhoeven algorithm \eqref{algorithm_1p}]\label{thm84}
Suppose that $\mathcal{X}$ and the $\mathcal{U}_m$ are all of finite dimension.
Let $x^{(0)}\in\mathcal{X}$, $u_1^{(0)}\in\mathcal{U}_1,\ldots,u_M^{(0)}\in\mathcal{U}_M$. Set $\tilde{s}^{(0)}=x^{(0)}+\tau\sum_{m=1}^M L^*_m u_m^{(0)}$. Let $\tau \in (0, 2/\beta)$ and $\sigma_1>0,\ldots,\sigma_M>0$ be such that $\tau\|\sum_{m=1}^M \sigma_m L_m^*L_m\|\leq 1$. Set $\delta=2-\tau\beta/2$. Let $(\rho^{(i)})_{i\in\mathbb{N}}$ be a sequence in $[0,\delta]$ such that $\sum_{i\in \mathbb{N}} \rho^{(i)} (\delta-\rho^{(i)})=+\infty$.
Then the sequences $(x^{(i)})_{i\in \mathbb{N}}$ and $\big((u_1^{(i)},\ldots,u_M^{(i)})\big)_{i\in \mathbb{N}}$ defined by the iteration \eqref{algorithm_1p} converge to a solution of \eqref{eq13p} and a solution of \eqref{eq21p}, respectively.\end{theorem}

When $\h$ is quadratic, as an application of \cref{thm33}, we have:

\begin{theorem}[Loris--Verhoeven algorithm \eqref{algorithm_1p}, quadratic case]\label{thm85}Suppose that $\h$ is quadratic. Let $x^{(0)}\in\mathcal{X}$, $u_1^{(0)}\in\mathcal{U}_1,\ldots,u_M^{(0)}\in\mathcal{U}_M$. Set $\tilde{s}^{(0)}=x^{(0)}+\tau\sum_{m=1}^M L^*_m u_m^{(0)}$. Let $\tau \in (0, 1/\beta]$ and $\sigma_1>0,\ldots,\sigma_M>0$ be such that $\tau\|\sum_{m=1}^M \sigma_m L_m^*L_m\|\leq 1$. Let $(\rho^{(i)})_{i\in\mathbb{N}}$ be a sequence in $[0,2]$ such that $\sum_{i\in \mathbb{N}} \rho^{(i)} (2-\rho^{(i)})=+\infty$.
Then the sequences $(x^{(i)})_{i\in \mathbb{N}}$ and $\big((u_1^{(i)},\ldots,u_M^{(i)})\big)_{i\in \mathbb{N}}$ defined by the iteration \eqref{algorithm_1p} converge weakly to a solution of \eqref{eq13p} and a solution of \eqref{eq21p}, respectively.\end{theorem}
 
\subsection{The Condat--V\~u algorithm}

Let $M\geq 1$ be an integer. Let $\mathcal{X}$ and $\mathcal{U}_m$, $m=1,\ldots,M$, be real Hilbert spaces. Let $\ff\in\Gamma_0(\mathcal{X})$, let $\dr{g_m}\in \Gamma_0(\mathcal{U}_m)$, $m=1,\ldots,M$, and let $\h:\mathcal{X}\rightarrow \mathbb{R}$ be a convex and 
differentiable function with  $\beta$-Lipschitz continuous gradient for some real $\beta> 0$. Let $L_m:\mathcal{X}\rightarrow \mathcal{U}_m$, $m=1,\ldots,M$, be bounded linear operators. We consider the convex optimization problem:
\begin{equation}
    \minimize_{x\in\mathcal{X}}\, \dg{f(x)}+\dr{\sum_{m=1}^M g_m(L_m x)} + \bl{h(x)}.\label{eq18p}
\end{equation}
The dual problem is:
\begin{equation}
     \minimize_{(u_1,\ldots,u_M)\in\mathcal{U}_1\times\cdots\times\mathcal{U}_M}\, (\ff+\h)^*({\textstyle -\sum_{m=1}^M L_m^*u_m})+
     \dr{\sum_{m=1}^M g^*_m(u_m)}.\label{eq30p}
\end{equation}

Let $(\omega_m)_{m=1}^M$ be a sequence of positive weights.
Let $\tau>0$, $\sigma>0$, and 
let $(\rho^{(i)})_{i\in\mathbb{N}}$ be a sequence of nonnegative reals. Let $x^{(0)}\in\mathcal{X}$ and let $u_m^{(0)}\in\mathcal{U}_m$ for $m=1,\ldots,M$.
Applying Technique 2 to the Condat--V\~u algorithm \eqref{algocv} or  \eqref{algocv2} and setting $\sigma_m=\sigma\omega_m$, we obtain the iterations:
\ceq{align}{\label{algocvp}
&\,\mbox{\textbf{Condat--V\~u iteration, form I, for \eqref{eq18p} and \eqref{eq30p}}: for }i=0,1,\ldots,\notag\\[-1mm]
    &\left\lfloor
    \begin{array}{l}
x^{(i+\frac{1}{2})}= \mathrm{prox}_{\tau \ff}\big(x^{(i)}-\tau\bl{\nabla h}(x^{(i)})-\tau \sum_{m=1}^M L_m^*u_m^{(i)} \big)\\
 x^{(i+1)}=x^{(i)}+\rho^{(i)}(x^{(i+\frac{1}{2})}-x^{(i)})\\
     \mbox{For $m=1,\ldots,M$,}\\
\left\lfloor\begin{array}{l}
    u_m^{(i+\frac{1}{2})}=\mathrm{prox}_{\sigma_m \dr{g^*_m}}\big(u_m^{(i)}+\sigma_m L_m (2x^{(i+\frac{1}{2})}-x^{(i)})\big)\\
    u_m^{(i+1)}=u_m^{(i)}+\rho^{(i)} (u_m^{(i+\frac{1}{2})}-u_m^{(i)}).
    \end{array}\right.
\end{array}\right.}%

\ceq{align}{\label{algocv2p}
&\,\mbox{\textbf{Condat--V\~u iteration, form II, for \eqref{eq18p} and \eqref{eq30p}}: for }i=0,1,\ldots,\notag\\[-1mm]
    &\left\lfloor
    \begin{array}{l}
     \mbox{For $m=1,\ldots,M$,}\\
\left\lfloor\begin{array}{l}
    u_m^{(i+\frac{1}{2})}=\mathrm{prox}_{\sigma_m \dr{g^*_m}}\big(u^{(i)}+\sigma_m L_m x^{(i)}\big)\\
    u_m^{(i+1)}=u_m^{(i)}+\rho^{(i)} (u_m^{(i+\frac{1}{2})}-u_m^{(i)})
    \end{array}\right.\\
x^{(i+\frac{1}{2})}= \mathrm{prox}_{\tau \ff}\big(x^{(i)}-\tau\bl{\nabla h}(x^{(i)})-\tau \sum_{m=1}^M L_m^*(2u_m^{(i+\frac{1}{2})}-u_m^{(i)}) \big)\\
    x^{(i+1)}=x^{(i)}+\rho^{(i)}(x^{(i+\frac{1}{2})}-x^{(i)}).
\end{array}\right.}

If $\h=0$, the Condat--V\~u algorithm reverts to the Chambolle--Pock algorithm \eqref{algoCPp}--\eqref{algoCP2p}.

As an application of \cref{thm61}, we have:

\begin{theorem}[Condat--V\~u algorithm \eqref{algocvp} or \eqref{algocv2p}]\label{thm86}Let $x^{(0)}\in\mathcal{X}$, $u_1^{(0)}\in\mathcal{U}_1,\ldots,u_M^{(0)}\in\mathcal{U}_M$. Let $\tau >0$ 
and $\sigma_1>0,\ldots,\sigma_M>0$ be such that 
$\tau\big(\|\sum_{m=1}^M \sigma_m L_m^*L_m\|+\frac{\beta}{2}\big)<1$. 
Set $\delta=2-\frac{\beta}{2}(\frac{1}{\tau}-\|\sum_{m=1}^M \sigma_m L_m^*L_m\|)^{-1}>1$.
Let $(\rho^{(i)})_{i\in\mathbb{N}}$ be a sequence in $[0,\delta]$ such that $\sum_{i\in \mathbb{N}} \rho^{(i)} (\delta-\rho^{(i)})=+\infty$. Then the sequences $(x^{(i)})_{i\in \mathbb{N}}$ and $\big((u_1^{(i)},\ldots,u_M^{(i)})\big)_{i\in \mathbb{N}}$, defined either by the iteration \eqref{algocvp} or by the iteration \eqref{algocv2p}, converge weakly to a solution of \eqref{eq18p} and a solution of \eqref{eq30p}, respectively.\end{theorem}

When $\h$ is quadratic, as an application of \cref{thm62}, we have:

\begin{theorem}[Condat--V\~u algorithm \eqref{algocvp} or \eqref{algocv2p}, quadratic case]\label{thm87}Suppose that $\h$ is quadratic. Let $x^{(0)}\in\mathcal{X}$, $u_1^{(0)}\in\mathcal{U}_1,\ldots,u_M^{(0)}\in\mathcal{U}_M$. Let $\tau >0$ and $\sigma_1>0,\ldots,\sigma_M>0$  be such that $\tau\|\sum_{m=1}^M \sigma_m L_m^*L_m\|<1$ and $\tau \|\bl{Q}+\sum_{m=1}^M \sigma_m L_m^*L_m\|\leq 1$.
Let $(\rho^{(i)})_{i\in\mathbb{N}}$ be a sequence in $[0,2]$ such that $\sum_{i\in \mathbb{N}} \rho^{(i)} (2-\rho^{(i)})=+\infty$. Then the sequences $(x^{(i)})_{i\in \mathbb{N}}$ and $\big((u_1^{(i)},\ldots,u_M^{(i)})\big)_{i\in \mathbb{N}}$, defined  either by the iteration \eqref{algocvp} or by the iteration \eqref{algocv2p}, converge weakly to a solution of \eqref{eq18p} and a solution of \eqref{eq30p}, respectively.\end{theorem}

\subsection{The Davis--Yin algorithm}

Let $M\geq 1$ be an integer. Let $\mathcal{X}$ be a real Hilbert space. Let $\ff\in\Gamma_0(\mathcal{X})$, let $\dr{g_m}\in \Gamma_0(\mathcal{X})$, $m=1,\ldots,M$, and let $\h:\mathcal{X}\rightarrow \mathbb{R}$ be a convex and  
differentiable function with  $\beta$-Lipschitz continuous gradient for some real $\beta> 0$. We consider the convex optimization problem:
\begin{equation}
\minimize_{x\in\mathcal{X}} \dg{f(x)}+\dr{\sum_{m=1}^M g_m(x)}+\bl{h(x)}.\label{eqdypp}
\end{equation}

Let $(\omega_m)_{m=1}^M$ be a sequence of positive weights whose sum is $1$. Let $\tau>0$ and 
let $(\rho^{(i)})_{i\in\mathbb{N}}$ be a sequence of nonnegative reals. Let $s_m^{(0)}\in\mathcal{X}$ for $m=1,\ldots,M$. By applying Technique 1 to the Davis--Yin algorithm \eqref{eqdy}, we obtain the iteration:
\ceq{align}{\label{algorithm_222p}
&\,\mbox{\textbf{Davis--Yin iteration for \eqref{eqdypp}}: for }i=0,1,\ldots,\notag\\[-1mm]
    &\left\lfloor
    \begin{array}{l}
   x^{(i+\frac{1}{2})}=\mathrm{prox}_{\tau \ff }\big(\sum_{m=1}^M \omega_m s^{(i)}_m\big)\\
   \mbox{For $m=1,\ldots,M$,}\\
   \big\lfloor\;\,
s_m^{(i+1)}=s_m^{(i)}+\rho^{(i)}\Big( \mathrm{prox}_{ \frac{\tau}{\omega_m} \dr{g_m}}\big(2x^{(i+\frac{1}{2})} - s_m^{(i)}-\tau\bl{\nabla h}(x^{(i+\frac{1}{2})} )\big)-x^{(i+\frac{1}{2})}\Big).
\end{array}\right.}

As an application of \cref{thm73}, we have:

\begin{theorem}[Davis--Yin algorithm \eqref{algorithm_222p}]\label{thm88}Let $s_m^{(0)}\in\mathcal{X}$, for $m=1,\ldots,M$. Let $\tau \in (0, 2/\beta)$ and let $\omega_1>0,\ldots,\omega_M>0$ be reals, whose sum is $1$. Set $\delta=2-\tau\beta/2$. Let $(\rho^{(i)})_{i\in\mathbb{N}}$ be a sequence in $[0,\delta]$ such that $\sum_{i\in \mathbb{N}} \rho^{(i)} (\delta-\rho^{(i)})=+\infty$. Then the sequence $(x^{(i+\frac{1}{2})})_{i\in \mathbb{N}}$ defined by the iteration \eqref{algorithm_222p} converges weakly to a solution of \eqref{eqdypp}.\end{theorem}

If $\h=0$, we recover the Douglas--Rachford algorithm \eqref{algorithm_22p}. On the other hand, if $\ff=0$, we recover the Generalized forward-backward algorithm~\cite{rag13,rag19}.

\subsection{The PD3O algorithm}

Let $M\geq 1$ be an integer. Let $\mathcal{X}$ and $\mathcal{U}_m$, $m=1,\ldots,M$, be real Hilbert spaces. Let $\ff\in\Gamma_0(\mathcal{X})$, let $\dr{g_m}\in \Gamma_0(\mathcal{U}_m)$, $m=1,\ldots,M$, and let $\h:\mathcal{X}\rightarrow \mathbb{R}$ be a convex and 
differentiable function with  $\beta$-Lipschitz continuous gradient for some real $\beta> 0$. Let $L_m:\mathcal{X}\rightarrow \mathcal{U}_m$, $m=1,\ldots,M$, be bounded linear operators. We consider the convex optimization problem:
\begin{equation}
    \minimize_{x\in\mathcal{X}}\, \dg{f(x)}+\dr{\sum_{m=1}^M g_m(L_m x)} + \bl{h(x)}.\label{eq18p2}
\end{equation}
The dual problem is:
\begin{equation}
     \minimize_{(u_1,\ldots,u_M)\in\mathcal{U}_1\times\cdots\times\mathcal{U}_M}\, (\ff+\h)^*({\textstyle -\sum_{m=1}^M L_m^*u_m})+
     \dr{\sum_{m=1}^M g^*_m(u_m)}.\label{eq30p2}
\end{equation}

Let $(\omega_m)_{m=1}^M$ be a sequence of positive weights.
Let $\tau>0$, $\sigma>0$, and 
let $(\rho^{(i)})_{i\in\mathbb{N}}$ be a sequence of nonnegative reals. Let $s^{(0)}\in\mathcal{X}$ and let $u_m^{(0)}\in\mathcal{U}_m$ for $m=1,\ldots,M$.
Applying Technique 2 to the PD3O algorithm \eqref{eqpd3o1} and setting $\sigma_m=\sigma\omega_m$, we obtain the iteration:
\begin{align}\label{eqpd3o1p0}
&\,\mbox{\textbf{PD3O iteration for \eqref{eq18p2} and \eqref{eq30p2}}: for }i=0,1,\ldots,\notag\\[-1mm]
    &\left\lfloor
    \begin{array}{l}
x^{(i+\frac{1}{2})}=\mathrm{prox}_{\tau \ff}(s^{(i)})\\
 \mbox{For $m=1,\ldots,M$,}\\
\!\left\lfloor\begin{array}{l}
\!u_m^{(i+\frac{1}{2})}
\!=\! \mathrm{prox}_{\sigma_m \dr{g^*_m}}\Big(u_m^{(i)}\!+\!\sigma_m L_m\big(2x^{(i+\frac{1}{2})}\!-\!s^{(i)}
\!-\!\tau \bl{\nabla h}(x^{(i+\frac{1}{2})})\!-\!\tau \sum_{m=1}^M L_m^*u_m^{(i)}\big) \Big)\\
\!u_m^{(i+1)}=u_m^{(i)}+\rho^{(i)} (u_m^{(i+\frac{1}{2})}-u_m^{(i)}).
\end{array}\right.\\
s^{(i+1)}=s^{(i)}+\rho^{(i)} \big(x^{(i+\frac{1}{2})}-s^{(i)}  -\tau \bl{\nabla h}(x^{(i+\frac{1}{2})})-\tau \sum_{m=1}^M L_m^*u_m^{(i+\frac{1}{2})}\big).
\end{array}\right.\end{align}

If we introduce the variables $\tilde{s}^{(i)}=s^{(i)}+\tau\sum_{m=1}^M L^*_m u_m^{(i)}$ and $a^{(i)}=x^{(i+\frac{1}{2})}-\tau \bl{\nabla h}(x^{(i+\frac{1}{2})})-\tilde{s}^{(i)}$, for every $i\in\mathbb{N}$, we can rewrite the iteration as:
\ceq{align}{\label{eqpd3o1p}
&\,\mbox{\textbf{PD3O iteration for \eqref{eq18p2} and \eqref{eq30p2}}: for }i=0,1,\ldots,\notag\\[-1mm]
    &\left\lfloor
    \begin{array}{l}
x^{(i+\frac{1}{2})}=\mathrm{prox}_{\tau \ff}\big(\tilde{s}^{(i)}-\tau \sum_{m=1}^M L_m^*u_m^{(i)}\big)\\
a^{(i)}=x^{(i+\frac{1}{2})}-\tau \bl{\nabla h}(x^{(i+\frac{1}{2})})-\tilde{s}^{(i)}\\
\tilde{s}^{(i+1)}=\tilde{s}^{(i)}+\rho^{(i)} a^{(i)}\\
 \mbox{For $m=1,\ldots,M$,}\\
\big\lfloor\;
u_m^{(i+1)}=u_m^{(i)}+\rho^{(i)} \Big(\mathrm{prox}_{\sigma_m \dr{g^*_m}}\big(u_m^{(i)}+\sigma_m L_m(x^{(i+\frac{1}{2})}+a^{(i)}) \big)-u_m^{(i)}\Big).
\end{array}\right.}

If $\ff=0$, the PD3O algorithm \eqref{eqpd3o1p} reverts to the Loris--Verhoeven algorithm \eqref{algorithm_1p}. If $\h=0$, it reverts to the Chambolle--Pock algorithm \eqref{algoCPp}.

We can check that if $L_m \equiv \mathrm{Id}$ and $\sum_m \omega_m=1$, by setting $\sigma=1/\tau$, $\sigma_m=\omega_m\sigma$, 
$s_m^{(i)}=s^{(i)}-\frac{\tau}{\omega_m} u_m^{(i)}+\tau \sum_{m=1}^M u_m^{(i)}$ 
in the PD3O algorithm \eqref{eqpd3o1p0}, we recover the Davis--Yin algorithm \eqref{algorithm_222p}. 

As an application of \cref{thm72}, we have:

\begin{theorem}[PD3O algorithm \eqref{eqpd3o1p}]\label{thm89}
Suppose that $\mathcal{X}$ and the $\mathcal{U}_m$ are all of finite dimension.
Let $\tilde{s}^{(0)}\in\mathcal{X}$ and $u_1^{(0)}\in\mathcal{U}_1,\ldots,u_M^{(0)}\in\mathcal{U}_M$. Let $\tau \in (0, 2/\beta)$, $\sigma_1>0,\ldots,\sigma_M>0$  be such that $\tau\|\sum_{m=1}^M \sigma_m L_m^*L_m\|\leq 1$. Set $\delta=2-\tau\beta/2$. Let $(\rho^{(i)})_{i\in\mathbb{N}}$ be a sequence in $[0,\delta]$ such that $\sum_{i\in \mathbb{N}} \rho^{(i)} (\delta-\rho^{(i)})=+\infty$.
Then the sequences $(x^{(i+\frac{1}{2})})_{i\in \mathbb{N}}$ and $\big((u_1^{(i)},\ldots,u_M^{(i)})\big)_{i\in \mathbb{N}}$ defined by the iteration \eqref{eqpd3o1p} converge to a solution of \eqref{eq18p2} and a solution of \eqref{eq30p2}, respectively.\end{theorem}

\section{Conclusion}

We have made a small tour of proximal splitting algorithms and we have shown how a principled analysis in terms of primal-dual monotone inclusions makes it possible to derive existing and new algorithms, with general convergence guarantees. Several connections between apparently distinct algorithms have been established, and we have unleashed the relaxation potential of a large class of proximal splitting algorithms.\footnote{In practice, we recommend that everyone tries overrelaxation as follows:  if, for instance, convergence is guaranteed for $\rho^{(i)}\in (0,2)$, first tune the other parameters, like $\tau$, with $\rho^{(i)}=1$. Once this is done, try $\rho^{(i)}=1.9$;  this will accelerate convergence in 
\br{many} cases.}
Finally, we have derived parallel versions of the algorithms to minimize the sum of an arbitrary number of functions. This selected overview of proximal splitting methods is by no means exhaustive; for instance, we did not cover the forward-backward-forward splitting method and its applications~\cite{tse00,bri11,com12},\cite[section 26.6]{bau17} or Dykstra-like algorithms~\cite{bau08,cha15}, \cite[section 30.2]{bau17}. Moreover, research into proximal splitting techniques is still an active topic, with recent contributions like the asymmetric forward-backward-adjoint algorithm~\cite{lat17}, the forward-backward-half forward algorithm~\cite{bri18}, the  forward-reflected-backward algorithm~\cite{mal18}, or the Forward-Douglas--Rachford-Forward algorithm~\cite{ryu20c}. 
There are many avenues for extensions, 
 including the use of variable stepsizes or metrics~\cite{com14,com142,ped18,mal182}, the study of convergence rates and development of accelerated variants~\cite{bec092,dro14,dos15,lor15,dav16,cha162,bec17,teb18,ryu20,con20,kim21}, and the design of randomized versions~\cite{cev14,ric14,com15,cha18,ped19,gor20,sal20}. It will be very interesting to see how proximal splitting methods are enriched by these modern notions.

\bibliographystyle{siamplain}
\bibliography{IEEEabrv,biblio}

\end{document}